\documentclass[11pt]{amsart}
\usepackage[T1]{fontenc}
\usepackage{amssymb,amsopn,amsmath,amsthm,authblk,hyperref,cleveref,lipsum,scrextend,mathrsfs,array,graphicx,bbm,enumerate,wasysym,yfonts,ulem,dsfont,graphicx,enumitem}
\usepackage{stmaryrd}
\usepackage{subcaption}
\usepackage{verbatim}
\usepackage[all]{xy}
\usepackage{tikz}
\usepackage[mathlines]{lineno}
\usepackage{fullpage}
\usepackage[foot]{amsaddr}
\usepackage{xifthen}
\usepackage{xstring}
\usepackage{xparse}
\usepackage{float}
\usepackage{xcolor,todonotes}
\usepackage[left=2cm,right=2cm,top=2cm,foot=2.5cm]{geometry}
\newtheorem{thm}{Theorem}[section]
\newtheorem{lemma}[thm]{Lemma}
\newtheorem{rem}[thm]{Remark}
\newtheorem{defn}[thm]{Definition}
\newtheorem{claim}[thm]{Claim}
\newtheorem{prop}[thm]{Proposition}

\newcommand{\QQ}{\mathcal{Q}}

\newcommand{\Z}{\mathbb Z}
\newcommand{\N}{\mathbb N}
\newcommand{\R}{\mathbb R}
\newcommand{\E}{\mathbb E}

\renewcommand{\S}{\mathbf S}
\renewcommand{\P}{\mathbb P}

\renewcommand{\1}{\mathbf 1}
\newcommand{\mm}{\mathfrak{m}}

\newcommand{\F}{\mathcal{F}}
\newcommand{\FF}{\mathscr F}

\newcommand{\SM}{\mathsf{Q}}
\newcommand{\SMM}{\widetilde{\SM}}

\newcommand{\type}[1]{\text{Type }#1}
\newcommand{\m}{\mathfrak{m}}

\newcommand{\q}{\mathfrak q}
\renewcommand{\r}{\mathfrak r}
\newcommand{\p}{\mathfrak p}
\newcommand{\Mm}{\widetilde{\m}}
\newcommand{\Mq}{\widetilde{\mathfrak{q}}}
\newcommand{\Mp}{\widetilde{\mathfrak{p}}}

\newcommand{\Ham}{\mathcal{H}}

\newcommand{\ske}{\mathfrak s}
\newcommand{\PM}{\mathcal{M}}
\newcommand{\Qm}{\mathcal Q}
\newcommand{\Bc}{\mathcal B}

\newcommand{\MQm}{\widetilde{\Qm}}
\newcommand{\BM}{\mathfrak{Q}}
\newcommand{\BMM}{\widetilde{\BM}}

\newcommand{\tr}{\mathbf t}
\newcommand{\vroot}{v_r}
\newcommand{\vertex}{\mathfrak{v}}

\newcommand{\Lalg}{\subset}
\newcommand{\MLalg}{\prec}
\newcommand{\Dec}{\phi}
\newcommand{\MDec}{\widetilde{\Dec}}
\newcommand{\Msigma}{\widetilde{\sigma}}
\newcommand{\Maprox}[3]{\left\langle #1 \right\rangle_{#3}^{#2}}

\newcommand\encircle[1]{%
  \begin{tikzpicture}
    \node[draw,circle,inner sep=0.6pt] {#1};
  \end{tikzpicture}}

\newcommand{\glue}{\, \encircle{\scalebox{.3}{g} } \,}

\def\cro#1{\llbracket#1\rrbracket}

\DeclareMathOperator\supp{supp}
\DeclareMathOperator{\Leb}{Leb}

\makeatletter
\g@addto@macro{\endabstract}{\@setabstract}
\newcommand{\authorfootnotes}{\renewcommand\thefootnote{\@fnsymbol\c@footnote}}
\makeatother

\let \epsilon \varepsilon

\numberwithin{equation}{section}

\crefformat{footnote}{#2\footnotemark[#1]#3}

\title{The Markov Property on Decorated Metric Planar Maps}

\author{Pablo Araya$^\dagger$}
\address{$^\dagger$ Universidad de Chile,  Centro de Modelamiento Matemático (AFB170001), UMI-CNRS 2807, Beauchef 851, Santiago, Chile.}

\author{Luis Fredes$^*$}
\address{$^*$Univ. Bordeaux, CNRS, Bordeaux INP, IMB, UMR 5251, F-33400 Talence, France.}
\author{Avelio Sepúlveda$^\dagger$}

\begin{document}
	\maketitle
	\begin{abstract} 
    We introduce decorated metric planar maps: planar maps whose edges $e$ are replaced by segments $[0,w_e]$, together with a continuous function on the resulting metric space, the decoration. We consider laws under which, conditionally on the metric map, the decoration follows a statistical physics model such as the Ising model or the Gaussian free field: the values at the vertices interact through a quadratic Hamiltonian and, given these values, the decoration along each edge is a Brownian bridge. Working with 4-regular maps with a boundary (the duals of quadrangulations with a boundary), we define a Markov property that allows explorations to stop in the middle of an edge. We show that, under mild regularity assumptions, the rerooting-invariant families of laws satisfying this Markov property are exactly the Boltzmann decorated metric maps. We also prove a strong Markov property and show that the random submaps inducing a Markovian decomposition are exactly the stopping maps for a suitable filtration. Finally, we obtain analogous results for classical decorated planar maps.
	\end{abstract}
 

\section{Introduction} 
Random planar maps have attracted considerable interest in the probabilistic and combinatorial communities, particularly as discrete models of random two-dimensional surfaces. As it is natural from a probabilistic point of view, not only uniformly chosen models are interesting. In fact, maps decorated by a statistical physics model, i.e. where the probability of choosing a given map is proportional to the partition function of that statistical physics model, are natural and have been studied not only in the mathematical community but also in physics \cite{kazakov1986ising,boulatov1987ising,kostov1993n}.

In this paper, we are motivated by two main ideas: generalising planar maps to the context where edges have length, and identifying the probability measures on decorated maps which inherit the Markov property of the decorating model. The idea of adding length to a graph and naturally extending the statistical physics model there has appeared at different points in time; however, it was Lupu's work on the metric graph Gaussian free field \cite{Lupu} that started an avalanche of results regarding these models, as the fact that they both satisfy the Markov property and the intermediate value property allows \cite{lupu2018random,ding2020percolation,drewitz2022cluster,drewitz2023critical,drewitz2026critical}. 

One motivation for considering these models comes from the expected connections between decorated random planar maps and Liouville Quantum Gravity surfaces \cite{DKRV,DMS}. These continuous random surfaces are constructed from the exponential of a version of the Gaussian free field and exhibit remarkable integrable properties \cite{KRV,GKRV,NQSZ}.

The result of this paper are motivated by the results in \cite{araya2025characterisation} for the uniformly chosen quadrangulations, and as in that paper we are only going to work with $4$-regular metric maps with a boundary. These are the dual maps to quadrangulations with a boundary. The reason to work in this dual representation is that in our context it makes sense that decorations live on vertices instead of on faces, as in this first case the extension of the model to the edges is more natural.

Random quadrangulations have been extensively studied from both probabilistic and combinatorial perspectives. Explicit formulae for their enumeration are known \cite{tutte} and can also be obtained through bijective methods \cite{SCH98,BG}. Moreover, their local limits \cite{KR05,CM} and scaling limits, including the Brownian map and the Brownian disk \cite{LeGall,Miermont,BetM,CLG,BMR16}, are well understood. Many geometric properties of these limiting objects are also known \cite{LeGall2,LeGall3,GN,CC19,LGR1,LGR2}. Several of these important results rely on Markovian decompositions of random planar maps, commonly referred to in the community as the peeling procedure \cite{Cu}. In particular, Boltzmann measures play a central role in this framework due to their Markovian structure.

We also consider decorated 4-regular maps without edge lengths. An important example of a decoration is an Ising model. Ising-decorated random planar maps (triangulations) are considerably less understood than their undecorated counterparts, although explicit partition functions \cite{BM,BM2} and local limit results \cite{AMS,CT} are known for several models. Nevertheless, understanding their scaling limits as a metric space remains a major challenge. As in the metric setting, Markovian decompositions play an important role in studying these models. Let us mention that the decorations considered in this paper differ from those commonly used in bijective constructions of planar maps \cite{Bet,BetM,BM2}, since we focus on decorations that satisfy a joint Markov property together with the underlying map.

\subsection{Results for metric maps}

A planar map is a rooted graph embedded on the sphere.  We say that it is a metric planar map if it is the metric space that can be constructed by adding strictly positive length to each edge. We call a decoration a continuous function on a metric planar map.

In this work, we focus on decorated $4$-regular metric maps with a boundary. That is to say the law on $(\tilde m, \phi)$, where $\tilde m$ is a metric map constructed from a map in which each vertex has 4 neighbours, except possibly one, which we call the exterior vertex. This setup is dual to the one we used in \cite{araya2025characterisation}, that is the reason why we denote a $4$-regular map the dual representation of a quadrangulation.

\begin{figure}[h!]
    \centering
    \includegraphics[scale=0.8]{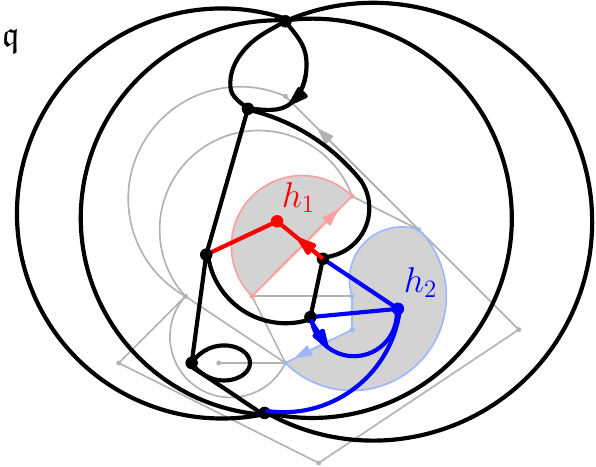}
    \caption{A 4-regular map $\q$ with a boundary and two hole vertices $h_1$ and $h_2$, represented by the shaded regions. The root edge of $\q$ is indicated by a black arrow, while the distinguished oriented edges associated with $h_1$ and $h_2$ are colored in red and blue, respectively. The underlying primal representation is shown in light grey.}
    \label{fig:holes}
\end{figure}

\begin{figure}[h!]
    \centering
    \includegraphics[scale = 0.7]{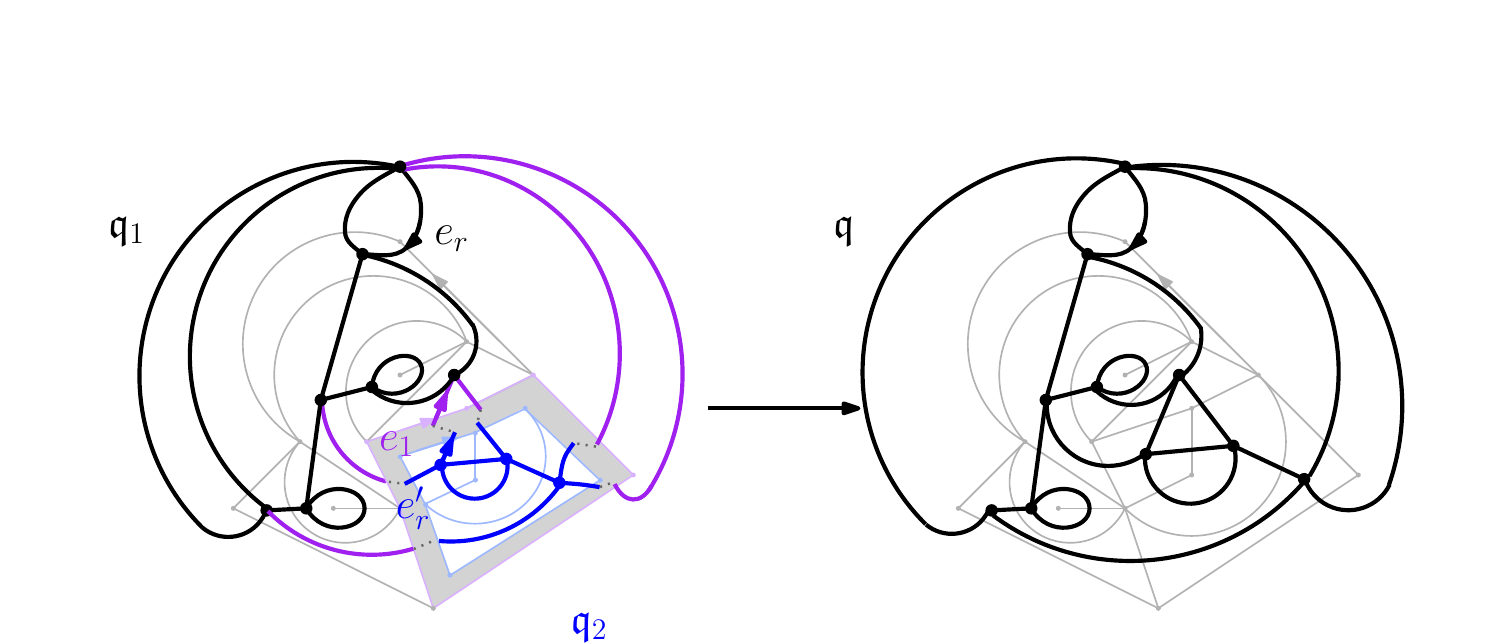}
    \caption{An example of a map $\q$ obtained from the gluing between $\q_1$ and $\q_2$. On the left, $\q_1$ is shown in black and $\q_2$ in blue, with their corresponding primal representations in the background. The grey region represents a hole vertex of $\q_1$, whose active boundary is colored in purple. The gluing identifies the marked edge $e_1$ of this boundary with the root edge $e_r$ of $\q_2$, and the remaining boundary edges according to their orientation. The resulting map $\q$ is shown on the right.}
    \label{fig:gluing}
\end{figure}

We can construct a 4-regular map $\BM$ with boundary from a 4-regular map with hole vertices\footnote{Here by a $4$-regular map with hole vertices we mean $4$-regular maps with boundary that have marked vertices, or hole vertices, whose degree is not necesarilly $4$ (see Section \ref{sec:2}).} $\q$ as follows. For each hole vertex $h$ of $\q$ we take a 4-regular map with boundary $\BM^\q_h$ with boundary size equal to the length of the hole vertex. Then, we glue to each hole vertex its respective 4-regular map to obtain $\BM$ (for a precise definition, see Section \ref{sec:3} and Figure \ref{fig:gluing}). As in \cite{araya2025characterisation}, we call this operation $\q\glue (\BM^\q_h)_h$. When there exist maps $(\BM^\q_h)_h$ such that $\BM= \q\glue (\BM^\q_h)_h$ we say that $\q$ is contained in $\BM$, or $\q\Lalg \BM$.

We are interested in probability measures on decorated metric maps whose decorations follows interactions between vertices that are neighbours. More precisely, given a measure $\mu$ on $\R$, we consider models in which the values of the decoration at the vertices interact through a quadratic Hamiltonian, while, conditionally on these values, the decorations along different edges are independent Brownian bridges. This construction includes metric versions of several classical spin systems. To be more precise, conditionally on a metric graph $\tilde q$ the law of the decoration should be
\begin{align}\label{e.density_metric_graph_1}
    &\P\left(\MDec\in d\Msigma \, \Bigl \vert \, \BMM=\Mq\right)\propto \left(\prod_{ij \in E(\q)}\hat{\P}_{w_{ij}\beta^{-1}}^{\sigma_i,\sigma_j}\left(d\Msigma|_{ij}\right)\right) \prod_{v\in V_i(\Mq)}\mu(d\sigma_v),
 \end{align}
 where $\beta>0$ is the inverse-temperature of the model, $\mu$ is its base measure, $w_{ij}$ is the length of the edge $ij$,  $V_i(\Mq)$ denotes the interior vertices and $\widehat{\P}_{w}^{a,b}$ is the unnormalised $n$-dimensional Brownian bridge measure of duration $w$ from $a$ to $b$. Its total mass is given by
\begin{align*}
    \frac{1}{(2\pi w)^{n/2}}
\exp\left(-\frac{(a-b)^2}{2w}\right).
\end{align*}
 This context is general enough to contain the case of Ising, when $\mu=\delta_{1}+\delta_{-1}$, the Gaussian free field when $\mu(dx)=dx$ and spin $O(N)$ models when $\mu$ is the uniform measure on $\S^{N-1}$.

As we have already discussed, we are interested in the Markov property. This is not such an easy property to describe. To state it, we consider families of probability measures on metric decorated planar $4$-regular maps with boundary $(\P^{\ell,b})_{\ell,b}$ indexed by $\ell$ the half-degree (or half-perimeter) of the exterior vertex and with boundary conditions $b$. The weak Markov property tells us that if we take $\BM \sim (\P^{\ell,b})_{\ell,b}$ and conditioned on the event that $\tilde \q \subseteq \BM$ and $\Dec\mid_{\q}=\sigma\mid _\q$ the law of $\BM$ inside each hole $h$ of $\q$ is independent and has the law of $\P^{\ell_{h},b_h}$, where $\ell_h$ is half the degree of $h$ and $b_h$ is the limit of the decoration at each edge of $h$ (see Definition \ref{def:markov_property_metric}).  Importantly here, $\tilde \q$ can `stop mid-edge' meaning that it can only contain a subset of an edge instead of the whole edge itself.

Our main result, Theorem \ref{t.characterization_metric}, is a characterisation of all families of (reasonable) probability measures that satisfy the Markov property and are invariant under rerooting. They are the Boltzmann decorated metric maps and are defined as follows
\begin{align}\label{e.density_metric_graph_2}
    &\P^{\ell,b}_{\lambda,q,\beta}\left(\BMM\in d\Mq,\MDec\in d\Msigma\right)\propto q^{|V_i(\Mq)|}\left(\prod_{ij\in E(\q)}\exp\left(-\lambda w_{ij} \right )\hat{\P}_{w_{ij}\beta^{-1}}^{\sigma_i,\sigma_j}\left(d\Msigma|_{ij}\right)\right) \prod_{v\in V_i(\Mq)}\mu(d\sigma_v),
 \end{align}

Let us remark that Boltzmann decorated metric maps are not obtained merely by assigning exponential weights to the edge lengths. The Brownian bridge contribution also introduces a nontrivial dependence on these lengths through the Gaussian heat kernel, which means that the edge lengths is not an exponential random variable.

The Markov property of random planar maps has traditionally been studied through the peeling procedure, in which the map is progressively revealed one face at a time (or equivalently, one vertex at a time in the dual representation). An alternative approach was introduced in \cite{araya2025characterisation}, based on the notion of stopping maps. This framework provides a more general class of random submaps that induce a Markovian decomposition, including those that cannot be constructed through the classical peeling procedure. In the metric setting, an exploration may stop in the interior of an edge, rather than only at its endpoints. In consequence particular care is required in choosing the indexing family for the corresponding filtrations, in order to allow such partial-edge explorations.

Our second main result concerns the characterisation of random metric submaps that induce a Markovian decomposition. In Theorem \ref{t.metric_strong_Markov}, we establish a strong Markov property for stopping metric maps: conditionally on the explored metric submap, the unexplored regions are independent Boltzmann metric maps with the corresponding boundary conditions. Conversely, we show that random metric submaps inducing that type of decomposition can be realised as stopping metric maps with respect to a suitable filtration.

\subsection{Non-metric setting} We also study the corresponding problem for decorated planar maps without edge lengths. In this setting, a decorated map is a pair $(\m,\Dec)$, where $\Dec$ assigns a value in $\mathbb R^n$ to each vertex of $\m$. We consider probability measures for which, conditionally on the underlying map $\m=\q$, the decoration has law

\begin{align}
\P(\Dec\in d\sigma|\m=\q)\propto \exp\left(-\frac{\beta}{2}\sum_{i\sim j}||\sigma_i-\sigma_j||^2\right)\prod_{v\in V(\mathfrak q)}\mu(d\sigma_v),
\end{align}

where the sum runs over adjacent vertices and the boundary values are prescribed. For a fixed graph, these distributions belong to the class of generalised Ising systems, studied at least since \cite{New}. Depending on the choice of $\mu$, this framework extends Ising-decorated maps, Gaussian free field decorations, and $O(n)$-spin systems.

As in the metric setting, we define a Markov property for families of probability measures indexed by the perimeter and boundary conditions. We prove in Theorem \ref{t.uniqueness_decorated} that the only rerooting-invariant families satisfying this property are the spin-decorated Boltzmann maps. Furthermore, Theorems \ref{t.strong_Markov_decorated} and \ref{t.local_implies_stopping_dec} provide the corresponding characterisation of random submaps inducing Markovian decompositions.

In particular, our results are motivated by the continuum framework of \cite{AG}, where a Markov property is obtained when cutting a Liouville Quantum Gravity surface by an $\mathrm{SLE}_\kappa$ curve whose central charge is not the natural one. In our discrete setting, a related situation arises, for instance, by considering a map decorated by two Ising systems and a stopping map whose exploration depends only on the first decoration. We believe that the discrete and continuum perspectives may provide useful insights into each other.

\subsection{Conjectural Scaling limits}
One of the main motivations of this work is to understand non-uniform models of decorated planar maps. In particular, when the decoration is given by a critical statistical physics model, one expects the underlying map to converge, after suitable rescaling, to a Liouville Quantum Gravity surface whose parameter is determined by the universality class of the decoration; for instance, $\gamma=\sqrt{3}$ for the Ising model and $\gamma=2$ for the Gaussian free field. If this conjectural picture holds, the stopping sets introduced in this paper should provide a discrete framework for understanding results such as those of \cite{AG}, where an LQG surface is cut by an SLE$_\kappa$ curve with a value of $\kappa$ different from the one naturally associated with the LQG parameter.

A particularly interesting case is that of metric maps decorated by the Gaussian free field (GFF). The GFF on the metric map has been shown to be an integrable model (see for example \cite{Lupu,ALS2,LW2}, we believe it is the same case for the GFF-decorated metric map. In fact in an upcoming work \cite{AFS27}, we will show that the metric map decorated with a GFF is integrable, in the sense that one can find both its local limits as well as understand its partition function and find its critical parameter. This should not be the case when one is working with the non metric decorated case.

Furthermore, we expect a richer phenomenon. Given a decorated metric map, define its first passage set $\mathbb{A}_b$ as the connected component containing the boundary of the set $\{x:\MDec(x)\leq b\}$. In the deterministic metric graph associated with $\eta \mathbb{Z}^2\cap D$, it was shown in \cite{ALS2} that this construction converges to the first passage set of the continuum domain introduced in \cite{aru2017first}. We expect an analogous convergence result to hold for GFF-decorated metric maps. Since the boundary condition of the remaining, unexplored field is constant, its law should depend only on the length of the boundary. Consequently, the evolution of this boundary length should generate a process closely related to the growth-fragmentation processes arising in planar maps \cite{Ber,BBCK}. More precisely, we expect this process to converge, up to a time change, to the growth-fragmentation process associated with the $O(2)$-decorated map and CLE$_4$ on a critical Liouville Quantum Gravity environment \cite{AHPS, Kam}.

\subsection{Ideas of proof}
As we have discussed, there are two main types of theorems in this work. The first one characterises laws over decorated metric maps (resp. decorated maps) that satisfy the Markov property and the second one characterises the (random) decorated metric submaps that satisfy the Markov property.

\subsubsection*{Characterisation of laws}
Let us first describe the proof strategy in the case of metric maps, the non-metric case is somehow simpler except at one interesting and non trivial point (see Claim \ref{c.start_induction_trees}). We start by introducing an intermediate condition. Restrict your attention to the measure obtained by conditioning on a the set of lengths of the edges of the metric map. We require that conditionally on this, the law of a given map is proportional to the partition function of the decoration.

To show that the intermediate condition implies that the measure is Boltzmann, we locally modify the metric map by adding a vertex without altering the energy of the original decoration. We show that the ratio between the densities of the original metric map and its modified version depends only on the length and decoration of the boundary. By expressing this ratio in two different ways, we deduce that it must be determined by a global parameter~$q$, corresponding to the Boltzmann weight.

Finally, we show that the Markov property and the rerooting invariance imply the intermediate step. To do this, we use different types of peeling steps to obtain the same decorated metric map. In this case, we use induction to reduce the number of vertices and conclude.

\subsubsection*{Characterisation of submaps}
 We again split the proof into two parts. First, we define stopping metric maps and show that they satisfy the Markov property. To do this, we follow \cite{araya2025characterisation} and define the notion of filtrations and stopping maps. However, some modifications are necessary for the metric setting. Filtrations need to be indexed by properly chosen equivalence classes of possible submaps rather than just submaps. The reason is that we need to ensure that the event $\{\Mq\subseteq \BMM\}$ has positive probability. Once the right definitions are in place, many of the proofs closely mirror classical arguments. Thus, the main difficulty is properly choosing the right indexation for the filtrations, which require some work to be done in the metric case.

Second, we show that every random metric submap that induces a Markovian decomposition must in fact be a stopping metric map. This is done by starting with a Markovian coupling and explicitly constructing a filtration $\FF$ such that the Boltzmann decorated metric map is an $\FF$-Boltzmann decorated metric map and the metric submap is an $\FF$-stopping metric map. The key observation is that the decomposition of a Boltzmann metric map along a deterministic metric submap $\Mq$ can be performed in successive steps: first revealing part of $\Mq$ and then decomposing the remaining unexplored metric map using the portion of $\Mq$ that has not been revealed yet.

 \medskip
 
The paper is organised as follows. In Section \ref{sec:3} we introduce metric maps and Boltzmann decorated metric maps, establish their weak and strong Markov properties, and show that the only metric submaps admitting a Markovian decomposition are the stopping metric maps with respect to a suitable filtration. In section \ref{sec:4} we characterise probability measures on the space of metric maps that satisfy the Markov property. Finally, in Section \ref{sec:5}, we turn to the simpler, non-metric setting by considering only the combinatorial skeletons of metric maps, thus disregarding edge lengths. We then adapt the techniques developed in the previous sections to establish analogous results.

\subsection*{Acknowledgements}
All authors are partially supported by MATH-AMSUD AMSUD250037. P.A. is supported by ANID-Subdireccion de Capital Humano/ Doctorado Nacional/ 2023-21231096. L.F is partially supported by the Project PEPS JCJC 2025–UMR 5251 (IMB), INSMI. The research of A.S.\,is supported by Centro de Modelamiento Matem\'{a}tico Basal Funds FB210005 from ANID-Chile, FONDECYT regular 1240884 and ERC 101043450 Vortex, and was supported by FONDECYT iniciaci\'on de investigaci\'on N$^o$ 11200085.

\subsection*{AI Declaration}
This paper together with \cite{araya2025characterisation} were originally only one paper that was uploaded to ArXiv on May 2025, so before the AI craziness of 2026. Thanks to the comments of some referees we decided to separate the original paper into two parts. Thus, this work benefited of AI only at the level of language checking and a referee like revision. More precisely, no mathematical results were obtained through interaction with AI, but we hope clarity was improved with our interaction with it.

\section{Preliminaries}\label{sec:2}

In this section, we introduce the necessary framework for the results of this paper. First, we discuss planar maps and their Markov property. In the second part, we study the conditional laws of Brownian bridges, which will play a key role in the definition of metric maps in later sections.

\subsection{Planar maps}\label{planar_maps}

A rooted planar map $\m$ is a finite connected graph properly embedded in the sphere with a distinguished oriented edge $e_r$ called its \textit{root}, whose starting vertex is called the \textit{root vertex}. In the literature the notion of rooted planar map is present in a different way (see for example \cite{CLG2}) where the roles of faces and vertices are interchanged. Working with the dual of a classic rooted planar map provides a more natural description of metric maps and their associated exploration so from here onwards we are going to present everything in the dual.

We denote
\begin{itemize}
    \item $V(\m);$ the set of vertices of $\m$,
    \item $E(\m);$ the set of edges of $\m$,
    \item $\overrightarrow{E}(\m);$ the set of oriented edges (or half-edges) of $\m$,
    \item $F(\m);$ the set of faces of $\m$.
\end{itemize}

Since we work in the dual representation, the root vertex $v_r$ corresponds to the boundary face in the usual primal representation. We define the perimeter of $\m$, denoted by $\text{Per}(\m)$, as the degree of $v_r$. In particular, when $\deg(v_r)=2\ell$, we say that $\m$ has half-perimeter $\ell$.

An important class of maps is that of maps with holes. Since we work with dual representation, we are going to call the holes hole vertices. We say that a map has \textit{hole vertices} if it is provided with a sequence of distinguished vertices $v_1,\dots,v_n$, each one called a hole vertex and with a marked outgoing oriented edge $e_1,\dots,e_n$. We denote by $H(\m)$ the set of hole vertices of $\m$, and by $\PM_H$ the space of rooted planar maps with hole vertices. 

Additionally, we define the \textit{active boundary} of $\m$, denoted by $\text{Active}(\m)$, as the set of edges incident to the hole vertices.

From here onward, and for the sake of simplification\footnote{All the results of this paper should extend easily in the context of more general maps.}, we work only with 4-regular maps, that is, rooted maps whose interior vertices have degree four. The distinguished hole vertices are allowed to have arbitrary degree. Additionally, we restrict to 4-regular maps $\q$ such that $\q$ forms a connected graph and every vertex of $\q$ is incident to a face. This condition is important because it implies the property of being rigid\footnote{This means that a triangulation is connected and no triangulation includes two distinct copies of the original map that coincides in the root.} (see Lemma 4.8 in \cite{angel2003uniform}) for 4-regular maps, which is key in order to have a Markov property. For $\ell,v\in\N$, we denote by $\Qm^{\ell,v}$ the set of such maps having one hole vertex of degree $2\ell$ and $v$ interior vertices. We further define
$$\Qm^\ell=\bigcup_v\Qm^{\ell,v}, \qquad \Qm=\bigcup_\ell\Qm^\ell, $$
and denote by $\Qm_H$ the corresponding space of maps with hole vertices.

Now, we define a way to glue maps. Let $\q_1$ be a rooted planar map with a distinguished hole vertex $v$, and let $\q_2$ be a rooted planar map with hole vertices such that the degree of its root vertex is equal to $\deg(v)$. We glue $\q_2$ into $v$ by identifying the half-edges incident to $v$ with those incident to the root vertex of $\q_2$, according to their cyclic order, starting from the distinguished oriented edge of $v$ and the root edge of $\q_2$. The resulting rooted planar map, denoted by $\q_1\glue\q_2$, inherits its root from $\q_1$.

The resulting map has vertex and edge sets obtained from $V(\q_1)\cup V(\q_2)$ and $E(\q_1)\cup E(\q_2)$ by identifying the corresponding vertices and edges. Its hole vertices are those of $\q_2$, and it has the same root as $\q_1$. More generally, if $\q_1$ has hole vertices $H(\q_1)$, we may glue a family of maps $(\q_v)_{v\in H(\q_1)}$, where the root vertex of each $\q_v$ has degree $\deg(v)$. We denote the resulting map by $\q_1 \glue (\q_v)_{v\in H(\q_1)}.$

\color{blue}
\color{black}

Thanks to this operation, we can define an order relation on maps. We say that a map with hole vertices $\q_1$ is a \textit{submap} of a map with hole vertices $\q_2$ with the same root of $\q_1$, denoted by $\q_1\Lalg \q_2$, if there exists a collection of maps with hole vertices $(\q_v)$, indexed by the hole vertices of $\q_1$, such that
\begin{itemize}
    \item $\text{Per}(\q_v) = \deg(v)$, for every $v\in H(\q_1)$,
    \item and, $\q_2 = \q_1\glue (\q_v)_{v\in H(\q_1)}$.
\end{itemize}

\begin{rem}\label{rem_bef_pee}
    The collection of maps $(\q_v)_{v\in H(\q_1)}$ that needs to be glued to a map $\q_1$ to obtain another map $\q_2$ is uniquely defined\footnote{The property of being rigid assures the uniqueness of the collection $(\q_v)_{v\in H(\q_1)}$ (see \cite{angel2003uniform})} for any map $\q_1 \Lalg \q_2$.
\end{rem}
This relation has a minimal element called the \textit{cemetery}, denoted by $\dagger$, which is the empty map. This map satisfies that $\dagger\Lalg\q$, for any $\q\in\Qm_H$.

\begin{figure}[h!]
    \centering
    \includegraphics[width=5in]{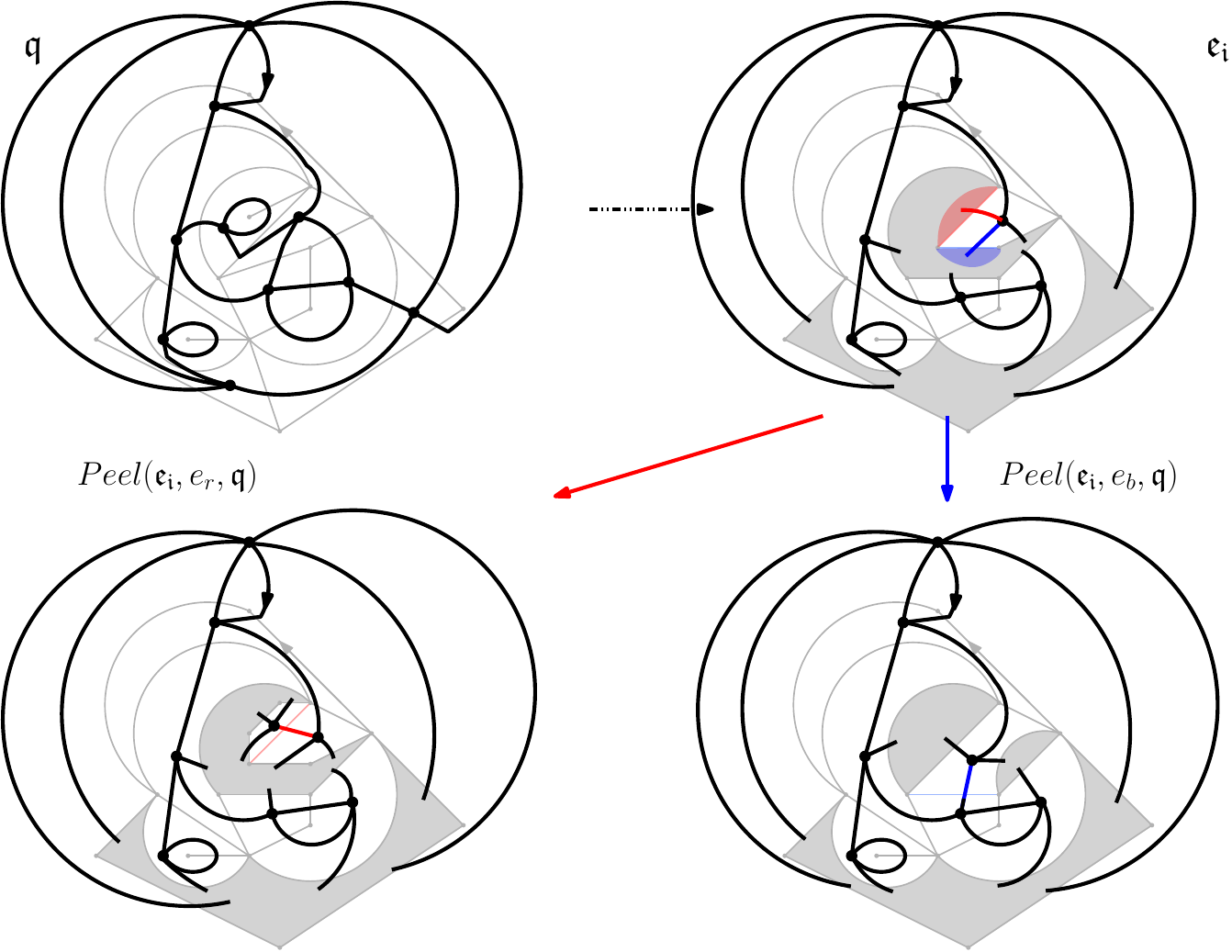}
    \caption{An illustration of the peeling procedure in the dual representation. The top-left panel shows the original map $\q$, while the top-right panel shows the partially explored map $\mathfrak{e}_i$ after several peeling steps. The discovered part is shown in white, and the grey regions represent the unexplored hole vertices, whose boundaries form the active boundary. The bottom panels illustrate two possible peeling steps: a Type 1 transition along the red edge $e_r$ (left) and a Type 2 transition along the blue edge $e_b$ (right).}
    \label{fig:peel}
\end{figure}

The previous order relation allow us to explore rooted planar maps in an algorithmic manner\footnote{Recall that we are working with finite planar maps.}, called the \textit{peeling exploration}. Starting from the root vertex $\mathfrak{e}_0$, the algorithm discovers a vertex of the map $\q$ in each iteration, giving as a result, a finite collection of maps $(\mathfrak{e}_i)_{i}$ with hole vertices, such that
\begin{align}
    \mathfrak{e}_0\Lalg\mathfrak{e}_1\Lalg \dots \Lalg\mathfrak{e}_n = \q.
\end{align}
This exploration depends on a function $\mathcal{A}$, called the \textit{peeling algorithm}, which takes a map with hole vertices $\mathfrak{e}$ and gives an edge $e$ from its active boundary. Let $e=\mathcal{A}(\mathfrak e)$ be the selected edge. Denote by $v_e$ the vertex of $\q$ incident to $e$ that does not belong to the current explored map. Then $v_e$ can be of two types.
\begin{itemize}
    \item\textbf{Peeling step of type 1:} The vertex $v_e$ has not been discovered yet. Then, $\mathfrak{e}_{i+1}$ is obtained by gluing $v_e$ with $\mathfrak{e}_i$ on $e$. The new active boundary set is obtained from the previous one by suppressing the selected edge $e$ and adding the newly exposed edges. In what follows, we denote by \type{$1$}  the peeling step of type 1.
    \item\textbf{Peeling step of type 2:} The vertex $v_e$ has already been discovered. In this case the selected hole vertex splits into two hole vertices joined by the selected edge. If their degrees are $2\ell_1$ and $2\ell_2$, then $ 2\ell_1+2\ell_2+2=\deg(v)$, where $v$ denotes the original hole vertex. Whenever one of the new hole vertices has degree zero, it is regarded as having disappeared. We denote by \type{$2,\ell_1,\ell_2$}  the peeling step of type 2 with $\text{Per}(h_1) = 2\ell_1$ and $\text{Per}(h_2) = 2\ell_2$.
\end{itemize}
 We denote by $Peel(\mathfrak{e_i},e,\q)$ to the resulting map of doing a peeling step on $e\in \text{Active}(\q) $ from $\mathfrak{e_i}$  (See Figure \ref{fig:peel}).  
 
We now state some important properties that may be satisfied by the law of a random planar 4-regular map.

\begin{enumerate}
    \item\label{Invariance-root} \textbf{(Invariance under rerooting)}
    $(\P^{\ell})_{\ell\in\N^{*}}$ is \textit{invariant under rerooting} if for any deterministic rooted 4-regular map $\q_1$ and any copy $\q_2$ of $\q_1$ differing only by the position of the root edge, we have that
    \begin{align*}
            \P^\ell(\BM=\q_1) = \P^\ell(\BM=\q_2), \mbox{ for any }\ell\in\N^*.
    \end{align*}
    \item\label{uniform-fix-faces} \textbf{(Uniformly distributed on $\Qm^{\ell,f}$)} $(\P^{\ell})_{\ell\in\N^{*}}$ is \textit{uniformly distributed in $\Qm^{\ell,v}$} if for any $\q_1,\q_2 \in \Qm^{\ell,f}$ 
    \begin{align*}
        \P^\ell\left(\BM= \q_1\right) = \P^\ell\left(\BM= \q_2\right), \mbox{ for any }\ell\in\N^*.
    \end{align*}

    \item\textbf{(Markov Property)} $(\P^{\ell})_{\ell\in\N^*}$ satisfies the \textit{Markov property}, if, for any deterministic 4-regular map with hole vertices $\q$, we can describe the conditional law  $\P^{\ell}(\cdot \mid \{\q\Lalg \BM\})$ as
    \begin{align*}
        \BM = \q\glue (\BM^{\q}_h)_{h\in H(\q)}.
    \end{align*}
Here, $(\BM^{\q}_h)_{h\in H(\q)}$ is a collection of independent maps with law $\P^{\ell(h)}$ for any $h\in H(\q)$.
\end{enumerate}

A well known probability law that satisfies all the previous properties is the Boltzmann map. 
\begin{defn}
For $q> 0$ and $\ell \in \N$, we say that $\BM$ is a \textit{$q$-Boltzmann map} with half-perimeter $\ell$ if the law of $\BM$ is supported on $\QQ^{\ell}$ and is given by 
\begin{align*}
    \P^{\ell}(\BM = \q) =\frac{1}{W_q^{\ell}}q^{|V(\q)|},\qquad \forall \q\in \QQ^{\ell},
\end{align*}
where the value of normalisation constant $W_{q}^{\ell}$ is given by
\begin{align*}
    W_{q}^{\ell} = \sum_{\q\in\Qm^\ell}q^{|V(\q)|}.
\end{align*}
\end{defn}
\begin{rem}
    Note that this probability measure only makes sense as long as $q$ is a positive constant such that $W_q^{\ell} < \infty$, which is true for $q\leq \frac{1}{12}$ (see for example \cite{Cu}). 
\end{rem}
We will only mention the parameter $q$ when necessary, but in general we omit it. 

The Boltzmann maps satisfies an even stronger version of the Markov property. In order to state the result we need the notion of filtrations and stopping maps.
\begin{defn}\label{def:filtration}
    We say that a collection of $\sigma-$algebras $\FF = (\F_{\q})_{\q\in\Qm_H}$ is a \textit{filtration} if it satisfies the property of \textit{Monotonicity},
\begin{align}
    \tag{Monotonicity} \mbox{For every }\q_1,\q_2\in\Qm_H\mbox{ such that }\q_1\Lalg \q_2\mbox{ one has } \F_{\q_1}\subseteq\mathcal{F}_{\q_2}.\label{filtration-monotonicity}
\end{align}
When working with a probability measure $\P$, the filtrations that we will use will also be complete for it, i.e
\begin{align}
    \tag{Completeness} \mbox{ }\F_{\q}\mbox{ is complete with respect to }\P,\mbox{ for every }\q\in\Qm_H.\label{filtration-completesness}
\end{align}
\end{defn}
For the rest of this section, $\BM$ is going to be a Boltzmann map. 
\begin{defn}\label{def_Boltzmann}
    We say that a Boltzmann map $\BM$ is an $\FF$-\textit{Boltzmann map} if it satisfies the following properties
\begin{itemize}[label={--}]
    \item\label{Adaptability}\textit{Adaptability:} for any $\q$ in $\Qm_H$, the event $\{\q\Lalg \BM\}$ is $\F_{\q}-$measurable.
    \item\label{Indep-increments}\textit{Independent increments:} conditionally on $\F_{\q}$ and on the event $\{\q\Lalg \BM\}$, the law of $(\BM^\q_v)_{v\in H(\q)}$ is that of a collection of independent $q$ Boltzmann maps with boundary given by $v$ for every $v\in H(\q)$.
\end{itemize}
\end{defn}
An special property of random submaps is that the $\sigma$-algebra generated by them can be characterized as follows.
\begin{prop}\label{prop:sigma_algebra}
     Let $\SM$ be a random map and for $\q\in \Qm_H$ consider the events $\Qm_{ \Lalg \q} := \{\SM\Lalg \q\}$. Then,
     \begin{align}\label{event_eq}
         \sigma\left(\SM\right) = \sigma\left(\Qm_{\Lalg \q}:\q\in\Qm_H\right).
     \end{align}
     Recall that as $\SM$ takes only countably many values, its $\sigma$-algebra is generated by the events $\{\SM=\p\}$ for $\p\in\Qm_H$.
\end{prop}

Now, we can define the random submaps that are going to satisfy the Markov property.
\begin{defn}\label{def_stop}
We say that a random planar map $\SM$ is an $\FF$-\textit{stopping map} for an $\FF$-Boltzmann map $\BM$ if
\begin{enumerate}
    \item $\P$-almost surely $\SM\Lalg \BM$,
    \item and, for any $\q\in\Qm_H$ we have that  $\left\{\SM\Lalg \q\right\}\in\F_{\q}$.
\end{enumerate}
We define the $\sigma$-algebra associated with an $\FF$-stopping map $\SM$ as follows
\begin{align*}
    \F_{\SM} := \left\{\left.\Theta\in\bigvee_{\q\in\Qm_H}\F_{\q}\right| \Theta\cap \{\SM\Lalg\q\}\in\F_{\q},\mbox{ for any }\q\in\Qm_H\right\}.
\end{align*}
\end{defn}
Next, we state the stronger version of the Markov property.
\begin{thm}[Strong Markov property] \label{t.strong_Markov}
    Let $\BM$ be an $\FF$-Boltzmann map and $\SM$ an $\FF$-stopping map for $\BM$. Then, $\BM$ can be decomposed as follows
    \begin{align*}
        \BM = \SM \glue (\BM^{\SM}_h)_{h\in H(\SM)},
    \end{align*}
    where, conditional on $\F_{\SM}$, $(\BM^{\SM}_h)_{h\in H(\SM)}$ is a collection of independent $q$-Boltzmann maps with boundary given by $h$ for every $h\in H(\SM)$.
\end{thm}
\subsection{On unnormalised Brownian Bridges} 
For $u,v\in\R$ and $w>0$, we defined the following measure, denoted by $\hat{\P}^{u,v}_{w}$
\begin{align}\label{e.unnormalised_Brownian_Measure}
    \hat{\P}^{u,v}_{w}\left[g\left((P_t)_{t\in[0,w]}\right)\right] = \frac{1}{\sqrt{2\pi w}}\exp\left(-\frac{1}{2w}(u-v)^2\right)\E^{u,v}_{w}\left[g\left((P_t)_{t\in[0,w]}\right)\right],
\end{align}
where $\E^{u,v}_{w}$ denotes the expected value of the unnormalised Brownian Bridge\footnote{We call it unnormalised Brownian Bridge as the expectation in play makes $\hat{\P}^{u,v}_{w}$ a positive measure but not necessarily a probability measure.} conditioned on starting at $u$ and finish on time $w$ at $v$.
The following lemma is a classical property of the Brownian bridges that allows us to decompose the law of a Brownian bridge at a point in the middle of its trajectory.
\begin{lemma}\label{lem-decomp}
Take $u,v\in \R$. Then, for all $w_1,w_2>0$
\begin{align*}
   \int_{\R}\hat{\P}^{u,z}_{w_1}\times \hat{\P}^{z , v}_{w_2}dz  =\hat{\P}^{u,v}_{w}, 
\end{align*}
where $w = w_1 + w_2$. 
\end{lemma}
\begin{proof}
    Take $f$ and $g$ two bounded measurable real-valued functions. Then,
    \begin{align*}
    \E^{u,v}_w\left[f((P_t)_{t\in[0,w_1]})g((P_t)_{t\in[w_1,w]})\right] &= \E^{u,v}_w\left[\E^{u,v}_w\left[f((P_t)_{t\in[0,w_1]})g((P_t)_{t\in[w_1,w]})|P_{w_1}\right]\right]\\
    &= \E^{u,v}_w\left[\E^{u,P_{w_1}}_{w_1}\left[f((P_t)_{t\in[0,w_1]})\right]\E^{P_{w_1},v}_{w_2}\left[g((P_t)_{t\in[0,w_2]})\right]\right],
    \end{align*}
where we used the Markov property to separate the conditional expected value. Then,
\begin{align*}
     &\E^{u,v}_w\left[f((P_t)_{t\in[0,w_1]})g((P_t)_{t\in[w_1,w]})\right] \\
     = &\int_{\R}\sqrt{\frac{w}{2\pi w_1w_2}}\exp\left(-\frac{w}{2w_1w_2}\left(z - \left(\frac{w_2}{w}u + \frac{w_1}{w}v\right)\right)^2\right)\E^{u,z}_{w_1}\left[f((P_t)_{t\in[0,w_1]})\right]\E^{z,v}_{w_2}\left[g((P_t)_{t\in[0,w_2]})\right]dz.
\end{align*}
Notice that we can rewrite the factor inside of the exponential as follows
\begin{align*}
    \frac{w}{2w_1w_2}\left(z - \left(\frac{w_2}{w}u + \frac{w_1}{w}v\right)\right)^2 = \frac{1 }{2 w_1 }(z-u)^2 +\frac{1 }{2w_2 }(z-v)^2 - \frac{1}{2w}\left(u-v\right)^2.
\end{align*}
Replacing this factor we obtain the following
\begin{align*}
    \E^{u,v}_w\left[f((P_t)_{t\in[0,w_1]})g((P_t)_{t\in[w_1,w]})\right]
    &= Z^{u,v}\int_\R \hat{\E}^{u,z}_{w_1}\left[f((P_t)_{t\in[0,w_1]})\right]\hat{\E}^{z,v}_{w_2}\left[g((P_t)_{t\in[0,w_2]})\right]dz,
\end{align*}
where the constant $Z^{u,v}$ is equal to 
\begin{align*}
    Z^{u,v}=\sqrt{2\pi w}\exp\left(\frac{1}{2w}(u-v)^2\right).
\end{align*}
 This proves the lemma.
\end{proof}

\begin{rem}
This measure may be degenerate when $w=0$, which forces $u=v$, but in our model $w$ will be an absolutely continuous random variable, meaning that a.s. it is bigger than $0$, so this problem is not present in our model.
\end{rem}

\begin{rem}\label{r.decomposition_BB}
Lemma \ref{lem-decomp} can also be written as
\begin{align*}
\hat \P^{u,v}_w[ dx] =  dx_{w_1} \hat \P^{u,x_{w_1}}_{w_1}(dx|_{[0,w_1]}) \hat \P^{x_{w_1},v}_{w_2}(dx|_{[w_1,w]}),
\end{align*}
where the density associated with the measure $\hat{\P}^{u,v}_w$ can be decomposed as the multiplication of the unnormalise densities associated with stop in the middle of the trajectory itself.
\end{rem}

\section{Metric Decorated maps: Markov property and associated measures}\label{sec:3}
In this section, we introduce a new model, which we call \textit{metric decorated 4-regular maps}. We modify the model of the previous section by assigning a positive length to each edge of a 4-regular map, so that every edge is naturally identified with an interval of the corresponding length. We first introduce metric 4-regular maps and the operations needed to explore them. We then enrich these metric 4-regular maps with a spin decoration. More precisely, conditionally on the underlying metric 4-regular map, the decoration is distributed according to a spin system defined on it. The canonical example is the Ising model, in which the spins take values in $\{-1,1\}$ at the vertices and are extended to take values in $\R$ along the metric edges. We then introduce decorated metric Boltzmann 4-regular maps with boundaries and show that they satisfy the weak Markov property. Finally, we develop the appropriate framework for the strong Markov property, allowing in particular for explorations that may stop in the middle of an edge. 

To define this model it is useful to slightly change our definition of a 4-regular map with a boundary $\q$. We add what we call phantom exterior vertices $V_e(\q)$. To define them, let us first define $V_i(\q)$ as the interior vertices of $\q$, and for each edge that connects an interior vertex with the root vertex we create a phantom root vertex. Thus, $V_e$ can be identified with $\{0,1,..,2\ell-1\}$ by following the direction of the root edge. We define $V(\q)$ as the union of $V_e(\q)$ with $V_i(\q)$. Notice that each phantom exterior vertex is incident to exactly one edge, which connects it to its corresponding interior vertex.

\begin{figure}[h!]
    \centering
    \includegraphics[scale=1]{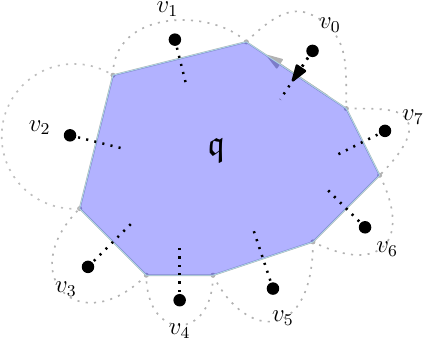}
	\caption{An illustration of the phantom vertices associated with the boundary of a map $\q$ of half-perimeter $\ell$. For each boundary position $k\in\{0,\dots,2\ell-1\}$, we introduce a phantom vertex $v_k$. The dotted segments indicate the correspondence between these vertices and the boundary of $\q$. These segments are auxiliary and are not actual edges of the map.}
    \label{fig:phantom_vertices}
\end{figure}

\subsection{Metric maps}\label{sec:metric}
We call \textit{metric map} the resulting metric space obtained by taking rooted map $\q$ and replacing each edge $e$ by a continuous line segment isometric to an interval $[0,w_e]\subseteq \R$. Metric graphs are classical objects in the theory of quantum graphs (see \cite{Kuc04,BK13}), and were later famously used by Lupu in \cite{Lupu} in the context of the Gaussian free field. Let us note that we do not require that the embedding of the map in $\mathbb{C}$ respects the length of the edges.

For metric maps it is useful to modify the concept of a hole vertex: what for planar maps was represented by a single hole vertex in a metric map will be represented by several boundary points, each lying at the tip of a partially explored edge belonging to the same unexplored hole.

We say that a metric map $\widetilde{\m}$ has \textit{hole} vertices if it is equipped with a finite collection
\begin{align*}
H(\widetilde{\m})=\{h_1,\ldots,h_n\}, \qquad h_i=\{v_{i,1},\ldots,v_{i,k_i}\}
\end{align*}
where each $h_i$ is a finite collection of distinguished points $v_{i,1},\dots,v_{i,k_i}\in \Mm$. Furthermore, we ask the following
\begin{enumerate}
\item[$i)$] All vertices of a hole $h_i$ have degree $1$. That is to say, they are locally homeomorphic to $[0,1)$.
    \item[$ii)$] The embedding can be continuously deformed such that for each hole the vertices can be joined without having edges crossing.
    \item[$iii)$] Each hole is maximal, meaning that there are no vertices $v$ such that $h_i\cup \{v\}$ satisfies both $i)$ and $ii)$.
\end{enumerate}
 The vertices belonging to the same $h_i$ are interpreted as the different boundary points corresponding to a single hole vertex of the underlying combinatorial map. We further distinguish one of the active oriented edges incident to each $h_i$, which determines the marked edge of the corresponding hole vertex. We call \textit{Active boundary} of $\Mm$, denoted by $\text{Active}(\Mm)$, the set of edges adjacent to any hole vertex of $\Mm$.

In this work we concentrate on metric 4-regular maps with boundary. These are metric maps that  can be obtained by starting with a planar map in which all vertices have degree 4 except for $v_r$, the root vertex. We denote the space of metric 4-regular maps with boundary as $\MQm$. Also, for $\ell,v\in\N$ we define $\MQm^{\ell,v}$ as the set of metric 4-regular maps with boundary of half-perimeter $\ell$ and $v$ interior vertices, $\MQm^{\ell} =\cup_v \MQm^{\ell,v}$. We denote by $\MQm_H$ the set of metric 4-regular maps with hole vertices; note that in this case every hole vertex of the original planar map becomes a set of $n$ points, where $n$ is the degree of the original hole vertex.

We define the skeleton function $\ske:\MQm\to\Qm$, which takes a metric 4-regular map $\Mq$ and returns the (minimal) map structure of the previous map, i.e. it forgets the lengths of the edges, and all vertices but the root vertex have degree 4. In the case of metric 4-regular maps with hole vertices the skeleton $\ske(\Mq)$ is obtained as before except that all points of the hole vertex are identified all the vertices belonging to the same hole $h_i$ and embedding it using the continuous deformation from (2). In particular, forgetting the lengths of edges one recovers a map with hole vertices, as introduced in \Cref{sec:2}.

We define a way of gluing metric maps as follows. Let $\Mq_1$ be a metric map with one hole vertex $v_h$, and $\Mq_2$ be a metric map such that $\text{Per}(\Mq_2) = \deg(v_h)$. We denote $\Mq_1\glue\Mq_2$ as the resulting metric 4-regular map that satisfies that 
\begin{itemize}
    \item $\ske(\Mq_1\glue\Mq_2) = \ske(\Mq_1)\glue\ske(\Mq_2)$,
    \item and the lengths are given by 
    \begin{align*}
        w_e = \begin{cases}w_{e_1}^{\Mq_1}, & \text{if }e \text{ comes from }e_1\Lalg\Mq_1,\\ w_{e_2}^{\Mq_2}, & \text{if }e \text{ comes from }e_2\Lalg \Mq_2, \\ w_{e_1}^{\Mq_1}+w_{e_2}^{\Mq_2}, &  \text{if }e_1 \text{ and }e_2\text{ are concatenated through their active endpoints.}\end{cases}
    \end{align*}
\end{itemize}
For every edge $e$ that comes from the identification of two active edges $e_1\in E(\Mq_1)$ and $e_2\in E(\Mq_2)$, their lengths are added only when their active endpoints are identified in the gluing. In this case,
\begin{align*}
    w_e=w_{e_1}^{\Mq_1}+w_{e_2}^{\Mq_2}.
\end{align*}
If the gluing identifies an active endpoint of either edge with the non-active endpoint of the other, the resulting identification is made before suppressing active vertices, and the metric structure is inherited from the resulting glued graph.

We can generalize the gluing procedure for a metric 4-regular map $\Mq$ with more than one hole vertex, where, in each hole vertex, we can glue a metric 4-regular map doing the same procedure as before.

We say that a metric 4-regular map with hole vertices $\Mq_1$ is an \textit{active submap} of $\Mq_2$, denoted by $\Mq_1\MLalg\Mq_2$, if
\begin{itemize}
    \item $\mathfrak{s}(\Mq_1)\Lalg\mathfrak{s}(\Mq_2)$, 
     \item and, for every edge $e\in \text{Active}(\Mq_1)$, $w_e^{\Mq_1}\leq w_e^{\Mq_2}$ with equality allowed only if $e\in\text{Active}(\Mq_2)$. In particular, if an active edge reaches its full length in the reference map, its endpoint is considered as discovered and the edge is no longer active.
\end{itemize}

Also, if for every edge $e\in E(\Mq_1)\setminus \text{Active}(\Mq_1)$, $w_e^{\Mq_1} = w_e^{\Mq_2}$ then we say that $\Mq_1$ is a \textit{submap} of $\Mq_2$, which we denote by $\Mq_1\Lalg\Mq_2$.

\begin{rem}\label{rem:eq}
Note that the relation $\Lalg$ is an order relation for metric maps but the relation $\MLalg$ is not. However, $\MLalg$ induces an order relation for the equivalence classes defined as follows: a metric 4-regular map with hole vertices $\Mq_1$ is \textit{equivalent} to $\Mq_2$, denoted by $\Mq_1\sim\Mq_2$, if
\begin{itemize}
    \item $\mathfrak{s}(\Mq_1) = \mathfrak{s}(\Mq_2)$,
    \item and, for every edge $e$ that is active for $\Mq_1$ and $\Mq_2$ we have that $w_e^{\Mq_1} = w_e^{\Mq_2}$.
\end{itemize} 
The relation $\sim$ is an equivalence relation and we denote by $[\Mq]_{\sim}$ the equivalence class of $\Mq$. In this context, the relation $\MLalg$ means the following 
\begin{align*}
    [\Mq_1]_{\sim} \MLalg [\Mq_2]_{\sim} \Longleftrightarrow (\forall \Mp_1\in [\Mq_1]_{\sim}) (\exists \Mp_2\in[\Mq_2]_{\sim}) \ \Mp_1\MLalg\Mp_2.  
\end{align*}
Under this definition, the relation $\MLalg$ is an order relation on the equivalence classes of the relation $\sim$.
\end{rem} 

The reason why we introduced two different orderings for metric maps is that the event $\{\Mq\Lalg\BMM\}$ is an event of null probability for Boltzmann metric maps but the event $\{\Mq\MLalg\BMM\}$ has positive probability. This is key to making the weak Markov property work. 

We now introduce a notion of monotone convergence for a sequence of metric 4-regular maps with hole vertices. Roughly speaking, such a sequence decreases to a metric 4-regular map with hole vertices if all maps share the same combinatorial structure and the length of each active edge decreases to the corresponding length in the limiting map.
\begin{defn}
Let $\Mq$ be a metric 4-regular map with hole vertices and $\left(\Mq_n\right)_{n\in\N}$ be a sequence of metric 4-regular maps with hole vertices. We say that $\left[\Mq_n\right]_{\sim}\searrow\left[\Mq\right]_{\sim}$ if 
\begin{itemize}
    \item $\ske(\Mq_n) = \ske(\Mq)$ for every $n\in\N$,
    \item and, $w_e^{\Mq_{n}}\searrow w_e^{\Mq}$ for every $e\in \text{Active}(\Mq)$.
\end{itemize}
\end{defn}
A natural example of a sequence that converges to a fixed metric map is the following.
\begin{defn}
Given a metric 4-regular map $\Mq$ and $\varepsilon>0$, we denote by $\Mq_{\varepsilon}$ the metric map satisfying
\begin{itemize}
    \item $\ske(\Mq_{\varepsilon}) = \ske(\Mq) $,
    \item and, $w_e^{\Mq_{\varepsilon}} = w_e^{\Mq} + \varepsilon$ for every edge $e\in E(\Mq)$.
\end{itemize}
\end{defn}
\begin{rem}
Taking $\varepsilon_n = 2^{-n}$ we have that $[\Mq_{\varepsilon_n}]\searrow[\Mq]$ when $n\to\infty$.
\end{rem}

\subsection{Decorated Metric maps}
In this subsection we extend the notion of metric maps by adding a decoration over them. For that consider a pair $(\Mq,\Msigma)$, where the first coordinate is a metric 4-regular map $\Mq\in\MQm_H$ and the second coordinate is a continuous function\footnote{Notice that we work only with decoration that belongs to $\R$ as it simplifies many of the computations, however almost all the results in this section remain true when the decoration belongs to $\R^n$.} $\MDec:\Mq\backslash \{\vroot\}\to\R$.

Let $b:\{0,\dots,2\ell-1\}\to \R$ be a function. We say that $\Msigma$ has \textit{boundary condition $b$} if the values on $V_e(\Mq)$ (the phantom vertices of $\Mq$) are equal to $b$ (after identification of $V_e(\Mq)$ with $\{0,..,2\ell-1\}$). Notice that if a metric 4-regular map has half-perimeter $\ell$ then $b\in \R^{2\ell}$. In order to simplify the notation, we denote $\Dec$ as the decoration $\MDec$ restricted to the interior vertices of the metric $4$-regular map $\Mq$, i.e. $\Dec = \MDec|_{V(\Mq)}$. We are going to work with a measure $\mu$ with $\supp(\mu)\subseteq \R$ and we denote by $\Bc$ the space of all boundary conditions for $\MQm$ that take values in $\supp(\mu)$. Depending on the $\supp(\mu)$ we are going to obtain different decorations such as the Ising model, the GFF, among others.

Let us now discuss the topology of decorated metric maps. This is important, as we need certain continuity properties of the measures we are studying. Take two decorated metric maps $(\tilde \q_1, \tilde \sigma_1)$ and $(\tilde \q_2, \tilde \sigma_2)$. We define the following distance between them
\begin{align*}
d\left((\tilde \q_1, \tilde \sigma_1),(\tilde \q_2, \tilde \sigma_2)\right)=\begin{cases}
\infty & \text{ if }\mathfrak{s}(\Mq_1) \neq \mathfrak{s}(\Mq_2),\\
\sum_{e \in E(\Mq_1)} d_{sko}(\tilde \sigma_1\mid_e,\tilde \sigma_2 \mid_e), & \text{ if }\mathfrak{s}(\Mq_1) = \mathfrak{s}(\Mq_2),
\end{cases}
\end{align*}
where Skorohod's distance is defined as in (12.16) of \cite{Bi}
\begin{align*}
    d_{sko}(X,Y) = \inf_{f} \left \{\|X-Y\circ f\|_\infty \vee \sup_{0\leq s\leq t\leq w_x}\left |\log \frac{f(t)-f(s) }{t-s } \right |\right\}.
\end{align*}
Here the infimum runs over bijective functions $f$ that are increasing going from $[0,w_x]$, the domain of $X$ (the length of the edge $x$), to $[0,w_y]$, the domain of $Y$. Note that the space generated by this distance is, in fact, a Polish space.

Now, we have all the necessary ingredients to define the main object of study of this section, the \textit{spin-decorated metric maps} (SDMM). To define them, fix $\beta >0$ the inverse temperature, and $\mu$ a measure on $\R$.

\begin{defn}
We say that a random pair $(\BMM,\MDec)$ is an SDMM of half-perimeter $\ell$ and boundary condition $b$ if its probability measure satisfies the following
\begin{align}
\P^{\ell,b}\left(\MDec\in d\Msigma\left|\BMM=\Mq \right . \right) \propto \prod_{ij\in E(\Mq)}\hat \P_{w_{ij}\beta^{-1}}^{\sigma_i,\sigma_j}\left((P_t)_{t\in[0,w_{ij}\beta^{-1}]} \in d\Msigma|_{ij}\right)\prod_{v\in V_i(\Mq)}\mu(d\sigma_v),
\end{align}
where $\sigma$ restricted to the phantom vertices is equal to $b$; and where $\hat \P^{c,d}_{w}$ is defined as in Lemma \ref{lem-decomp}: the unnormalised probability measure of a Brownian Bridge of length $w$ started at $c$ and finishing at $d$. Here the edge is oriented going from  $i$ to $j$. \footnote{Here we assume we have numbered the vertices and we point the edges from smallest to biggest endpoint of each edge.}
\end{defn}

Let us give an equivalent representation of an SDMM.
\begin{rem}
Let $(\BMM,\MDec)$ be an SDMM. Then for any $(\tilde q, \sigma)$ 
\begin{align} \label{def_cond_metric}
    \P^{\ell,b}\left(\Dec\in d\sigma\left|\BMM=\Mq \right.\right) &=  \frac{1}{Z^{\beta,b}(\Mq)}\exp\left(-\Ham^{\beta,b}(\Mq,\sigma)\right)\prod_{v\in V_i(\Mq)}\mu(d\sigma_v),
\end{align}
where 
\begin{align}\label{Zq_metric}
&\Ham^{\beta,b}(\Mq,\sigma):= \frac{\beta}{2}\displaystyle\sum_{\substack{i,j\in V(\Mq)\\i\sim j}}\frac{\|\sigma_i-\sigma_j\|^2}{w_{ij}} \ \ \text{ and }
    &Z^{\beta,b}\left(\Mq\right):= \int_{\R^{|V_i(\q)|}}\exp\left(-\Ham^{\beta,b}(\Mq,\sigma)\right)\prod_{v\in V_i(\Mq)}\mu(d\sigma_v).
\end{align}
Furthermore, given $\Dec$ and $\widetilde \BM$, the law of $(\widetilde \phi|_{ij})_{ij\in E(\widetilde \BM)}$ is that of independent Brownian bridges of length $w_{ij}$ starting at $\Dec_i$ and ending at $\Dec_j$. 
\end{rem}

A special case of SDMM are Boltzmann decorated metric maps. They are going to be the focal point of the first part of this section.
\begin{defn}\label{def:BDMP}
We say that $(\BMM,\MDec)$ is a \textit{$(q,\lambda,\beta)$-Boltzmann decorated metric map} if its probability measure satisfies the following

\begin{align}\label{e.Boltzmann-metric}
    &\P^{\ell,b}_{\lambda,q,\beta}\left(\BMM\in d\Mq,\MDec\in d\Msigma\right)\propto q^{|V_i(\Mq)|}\prod_{ij\in E(\q)}\exp\left(-\lambda w_{ij} \right )\hat{\P}_{w_{ij}\beta^{-1}}^{\sigma_i,\sigma_j}\left(d\Msigma|_{ij}\right) \prod_{v\in V_i(\Mq)}\mu(d\sigma_v),
 \end{align}

where the value of the normalisation constant $W_{\lambda,q,\beta}^{\ell, b}$ is given by 
\begin{align}
    W_{\lambda,q,\beta}^{\ell, b} = \sum_{\q}q^{|V_i(\q)|} \int_{\R^{V_i(\Mq)}}  \int_{\R_+^{E(\q)}}     \left( \prod_{ij\in E(\Mq)}\frac{\exp\left (-\frac{\beta}{2} \frac{(\sigma_i-\sigma_j)^2}{w_{ij}} -\lambda w_{ij} \right )  }{ \sqrt{2\pi w_{ij}/\beta}} \right ) \prod_{ij \in E(\Mq)} dw_{ij}  \prod_{v\in V_i(\Mq)} \mu(d\sigma_v),
\end{align}
where the sum is taken over all the maps with half-perimeter $\ell$.
\end{defn}
An attentive reader may realize that the parameters $\lambda$ and $\beta$ play a similar role, as the longer an edge is, the colder the system gets. The truth is slightly more complicated, and is explained in the following remark.
\begin{rem}
The formula \eqref{e.Boltzmann-metric} shows the law decomposed as : firstly a weight according to the number of vertices, secondly the law of lengths and finally the law of decoration along the edges (which depends on the lengths). Two cases of interest arise:
\begin{itemize}
    \item \underline{$\beta>0$} If $(\BMM, \MDec)$ is a $(q,\lambda,\beta)$-Boltzmann decorated metric map, then define $\BMM/\beta$ the metric map where all distances have been divided by $\beta$, $\Dec_\beta(x)= \Dec(\beta x)$, where $x \in \BMM/\beta$. Then $(\BMM/\beta, \MDec_\beta)$ is a $(q\beta^2,\lambda\beta,1)$-Boltzmann decorated metric map. Thus, from here onward we always assume that $\beta=1$, when $\beta>0$.
    \item \underline{$\beta=0$} The decomposition formula \eqref{e.Boltzmann-metric} for $\beta=0$ gives that the spins along the edge are independent (infinite length Brownian bridges give that the spin at two points are not correlated), which can be reduced to undecorated maps with edges of exponential random lengths.  This model is in itself interesting and has been previously studied \cite{Cu,CLG2}. Our techniques also work for this case, as it is simpler.
\end{itemize}
These cases let us focus on the behavior of $q$ and $\lambda$.
\end{rem}

Furthermore, as expected Boltzmann maps are SDMM.
\begin{rem}
The Boltzmann decorated metric maps are indeed an SDMM as
\begin{align*}
    \P^{\ell,b}\left(\MDec\in d\Msigma\left|\BMM=\Mq \right.\right) 
    &\propto  q^{|V_i(\Mq)|}\prod_{ij\in E(\q)}\exp\left(-\lambda w_{ij} \right ) \hat{\P}_{w_{ij}\beta^{-1}}^{\sigma_i,\sigma_j}\left( d\Msigma|_{ij}\right)  dw_{ij}\prod_{v\in V_i(\Mq)}\mu(d\sigma_v) \\
    &\propto\prod_{ij\in E(\Mq)} \hat{\P}_{w_{ij}\beta^{-1}}^{\sigma_i,\sigma_j}\left( d\Msigma|_{ij}\right)\prod_{v\in V_i(\Mq)}\mu(d\sigma_v),
\end{align*}
since over this conditioning $q^{|V_i(\Mq)|}\displaystyle\prod_{ij\in E(\q)}\exp\left(-\lambda w_{ij} \right )$ is a constant factor.
\end{rem}
And of course, Boltzmann metric maps are not trivial.
\begin{rem}
    If the measure $\mu$ is bounded, then there always exists a $q>0$  and $\lambda>0$ s.t. for any $\ell$ and $b$, $W_{q,\lambda,\beta}^{\ell,b}<\infty$. This is because $e^{-\lambda w}/\sqrt{w}$ is integrable in $\R^+$. The case where the measure $\mu$ is Lebesgue, can be treated as in the decorated case as the partition function for the GFF is explicit and increases exponentially with the number of vertices. For the case when $\mu$ is the sum of Dirac's in $\Z$, a similar bound as in the decorated  case allows us to conclude.
\end{rem}

When the decoration belongs to $\R^n$, some modifications need to be added.
\begin{rem}
When the decoration belongs to $\R^n$, we need to use $n$-dimensional Brownian bridges instead of the $1$-dimensional ones. These measures are normalized by a factor $\beta^{d/2}/(2\pi w_{ij})^{d/2}$. This is a priori non-integrable in $w$, thus we need to further ask a condition on the measure $\mu$. It is enough that for any point in the support of the measure, the mass $\mu(B(x,r))$ is smaller than a constant times $r^{d/2-1+\epsilon}$ for some $\epsilon>0$. 
\end{rem}

Now, we state the analogous properties of the invariance under rerooting and Uniform distribution on $\MQm^{\ell,f}$.
\begin{enumerate}
    \item\label{Invariance-root-dec-met} \textbf{(Invariance under rerooting)}
    The SDMM $(\BMM,\phi)\sim \P^{\ell,b}$ with half-perimeter $\ell$ and boundary condition $b$ is \textit{invariant under rerooting} if after rerooting it, its law is that of $\P^{\ell,b_s}$ where $b_s$ is the proper shift of $b$.
    
    \item\label{gibbs-fix-faces-dec-met} \textbf{(Gibbs distribution on $\MQm^{\ell,v}$)} The SDMM $(\BMM,\Dec)\sim\P^{\ell,b}$ has \textit{Gibbs distribution on $\MQm^{\ell,v}$} if for any $\Mq_1,\Mq_2 \in \MQm^{\ell,v}$ such that there exists a bijection $g$ between $E(\Mq_1)$ and $E(\Mq_2)$ that satisfies that $w_e = w_{g(e)}$ for any $e\in E(\Mq_1)$, then\footnote{We say that if $\tilde \q_1$ and $\tilde \q_2$ satisfy these conditions, we say that they are comparable.}
    \begin{align*}
        \frac{\P^{\ell,b}\left(\BMM\in d\Mq_1\right)}{\P^{\ell,b}\left(\BMM\in d\Mq_2\right)} = \frac{Z^b\left(\Mq_1\right)}{Z^b\left(\Mq_2\right)}.
    \end{align*}
\end{enumerate}

A canonical example of an SDMM that satisfies these properties is the Boltzmann decorated metric map.

\begin{prop}\label{prop:boltzman_invariance}
    Boltzmann decorated metric maps are invariant under rerooting (\ref{Invariance-root-dec-met}) and have Gibbs distribution over $\MQm^{\ell,f}$ (\ref{gibbs-fix-faces-dec-met}). 
\end{prop}
\begin{proof}
$\ $
\begin{itemize}
\item\textbf{Invariance under rerooting:}
Let $(\Mq,\Msigma)$ be a decorated metric map with boundary condition $b$, and denote by $(\Mq^s,\Msigma^s)$ its rerooting when shifting the root $s$ edges along the boundary. Denote by $b_s$ the corresponding shift of $b$. Then, 
\begin{align*}
\P^{\ell,b}_{\lambda,q}(\BMM\in d\Mq,\MDec\in d\Msigma) \propto q^{|V_i(\Mq)|} \prod_{ij\in E(\Mq)} \exp(-\lambda w_{ij}) \hat{\P}^{\sigma_i,\sigma_j}_{w_{ij}} (d\Msigma|_{ij})dw_{ij} \prod_{v\in V_i(\Mq)}\mu(d\sigma_v).    
\end{align*}
Rerooting only relabels the vertices and edges. In consequence, the right-hand side is invariant if we change $b$ to $b_s$ (using the reversibility of the Brownian bridge whenever an edge orientation is reversed). Therefore,
\begin{align*}
W^{\ell,b}_{\lambda,q} = W^{\ell,b_s}_{\lambda,q}, 
\end{align*}
and consequently the rerooted pair $(\BMM^{s},\MDec^{s})$ has law $\P^{\ell,b_s}_{\lambda,q}$.

    \item\textbf{Gibbs distribution on $\MQm^{\ell,v}$:} Take $\Mq_1,\Mq_2\in\MQm^{\ell,v}$ as in the hypothesis. For $j=1,2$, we have the following
    \begin{align*}
        \P^{\ell,b}\left(\BMM\in d\Mq_j\right) &= q^{|V_i(\Mq_j)|}\left(\prod_{e\in E(\Mq_j)}\frac{e^{-\lambda w_e}}{\sqrt{2\pi w_e}}\right)Z^{b}(\Mq_j).
    \end{align*}
From here we conclude the desired equality since $\Mq_1$ and $\Mq_2$ have the same number of vertices and there exists a bijection $g:E(\Mq_1)\to E(\Mq_2)$ such that $w_e = w_{g(e)}$, which implies that 
\begin{align*}
    \prod_{e\in E(\Mq_1)}\frac{e^{-\lambda w_e}}{\sqrt{2\pi w_e}} = \prod_{e\in E(\Mq_2)}\frac{e^{-\lambda w_e}}{\sqrt{2\pi w_e}}.
\end{align*}
\end{itemize}
\end{proof}

Let us now show that Boltzmann metric maps are continuous with respect to their boundary conditions. 
\begin{prop}\label{p.continuity_metric}
    Take  $(\BMM,\MDec)$ a $(q,\lambda)$-Boltzmann decorated metric map and assume there is $\delta>0$ such that $\sup_{b'\in B(b,\delta)}W^{\ell,b'}_{q,\lambda}<\infty$.    Then $ \P_{\lambda,q}^{\ell, b_n}\to \P_{\lambda,q}^{\ell, b}$ as $b_n\to b$.
\end{prop}
\begin{proof}
Note that 
\begin{align*}
&\left |\exp\left (-\frac{ (x-y)^2}{2w_{ij}}\right )-\exp\left (-\frac{ (x'-y)^2}{2w_{ij}}\right )\right  |\\
\leq &\frac{|x-x'|}{ w_{ij}}\left|\frac{ (x+x')}{2} - y\right| \exp\left (-\frac{ (x-y)^2\wedge (x'-y)^2}{2w_{ij}}\right) .
\end{align*}

Using this inequality in the edges of $\q$ that intersect the boundary, we can use dominated convergence to see that $W^{\ell,b_n}_{q,\lambda}- W^{\ell,b}_{q,\lambda}\to 0$. This directly implies that the marginal law of $\ske(\BMM)$ under $\P^{\ell,b_n}$ converges to that under $\P^{\ell,b_n}$. Furthermore, one can check that both the length, and the decoration on vertices are also converging. Finally, the trickiest topology is that of the decoration over edges, this is converging as the law of Brownian bridge is continuous for the weak convergence (of the Skorohod topology) in the length, initial point and end point.
\end{proof}

\subsection{Markov property: stopping maps} In this section, we present the Markov property for the metric decorated maps. We begin by stating the weak Markov property for the Boltzmann decorated metric map, and then, we describe the modification of the definition of stopping maps in \cite{araya2025characterisation} for the metric maps in order to have the Markov property.
\subsubsection{Weak Markov property.}

\begin{thm}\label{t.metric_weak_markov}
    Let $(\BMM,\MDec)$ be a Boltzmann decorated metric map and $\Mq\in\MQm_{H}$. Then, conditionally on $\{\Mq\Lalg\BMM\}$ and $\MDec|_{\Mq}$, $(\BMM,\MDec)$ can be decomposed as follows
    \begin{align} \label{e.decomposition_metric}
        \BMM = \Mq\glue (\BMM_h^{\Mq})_{h\in H(\Mq)},\  \text { and } \ 
        \MDec = \MDec|_{\Mq} + \MDec^{\Mq},
    \end{align}
    where, $(\BMM_{h}^{\Mq},\MDec_{h}^{\Mq})$ is a collection of independent Boltzmann decorated metric maps with boundary $\ell(h)$ and boundary condition $\Dec|_{h}$ for every $h\in H(\Mq)$.
\end{thm}

To prove this theorem, we will make use of the following result concerning conditional laws

\begin{lemma}[Theorem 8.5 of \cite{kallenberg2001foundations}]
\label{l.Ley_condicional}
Let $\mathbf X$ and $\mathbf Z$ be two Polish spaces, $\mu$ and $\nu$ two measures on $\mathbf X$ and $\mathbf Z$ respectively and $F:\mathbf X\times \mathbf Z\to \R^+$ a measurable function. Now, take a pair $(X,Z)$ of random variables whose law is proportional to $ F(x,z)\mu(dx)\nu(dz)$. Then, a regular conditional law of $X$ given $Z$ is 
\begin{align*}
 \P(X\in dx\mid Z=z) \propto F(X,z)\mu(dx).
\end{align*}
\end{lemma}

By this notation we mean that the law of $X$ given $Z$ is absolutely continuous with respect to $\mu$ and its Radon-Nikodym derivative is proportional to $F(X,z)$.

Notice that $\Mq$ can have a small part of an edge, in that case, we use the value associated with the tip of $\Mq$ as boundary condition. Here we use the relation $\MLalg$, to work on a positive probability event, which allow us to use Lemma \ref{l.Ley_condicional}.

\begin{proof}[Proof of Theorem \ref{t.metric_weak_markov}]
Let $\Mp\in\MQm_H$ such that $\Mq\Lalg\Mp$ and $\Msigma$ a decoration over $\Mp$.  Denote by $(\Mp_h)_{h\in H(\Mq)}$ the collection of maps such that $\Mp = \Mq\glue(\Mp_h)_{h\in H(\Mq)}$. Let us write the density of the decorated metric map $(\BMM,\MDec)$ under the event $\{\Mq\MLalg\BMM\}$ as follows
\begin{align}
\label{e.law_weak_Markov}
    &\P^{\ell,b}\left[\BMM\in d\Mp, \MDec \in d\Msigma \mid\Mq\MLalg\BMM\right]
    \\
    \nonumber\propto \ &\left(\prod_{h\in H(\Mq)}q^{|V_i(\Mp_h)|} \prod_{e\in E(\Mp_h)}e^{-\lambda w_e}\right)\left(\prod_{e\in E(\Mq)}e^{-\lambda w_e}\right)\left(\prod_{e\in \text{Active}(\Mq)}\1_{w_{e}^{\Mp}-w_{e}^{\Mq}\geq 0}\right)\\
    \nonumber&\left[\sum_{\mm\supseteq\q}\delta_{\mm}(\q\glue(\p_h)_{h\in H(\Mq)}) \prod_{v\in V_i(\mm)}\mu(d\sigma_v) \prod_{ ij\in E(\mm)}\hat{\P}^{\sigma_i,\sigma_j}_{w_{ij}}\left( d\Msigma|_{ij}\right)dw_{ij}\right].
\end{align}

Note that, in each edge of the active boundary of $\Mq$, we can decompose the probability of the Brownian bridge into two parts: the explored and unexplored parts of the edge. We now use Lemma \ref{lem-decomp} and call $k$ the point in the edge $ij$ that is at distance $w_{ij}^{\Mq}$ from $i$, that is to say $w_{ik}=w_{ij}^{\Mq}$. Then, recalling the notation of Remark \ref{r.decomposition_BB}
\begin{align}\label{eq:pivot}
    \hat{\P}^{\sigma_i,\sigma_j}_{w_{ij}}\left( d\Msigma|_{ij}\right)\1_{w_{ij}\geq \tilde w_{ij}^{\Mq}}dw_{ij}
    = d\sigma_{k}  \hat{\P}^{\sigma_i,\sigma_k}_{\tilde w_{ik}}\left( d\Msigma|_{ik}\right) \hat{\P}^{\sigma_k,\sigma_j}_{w_{kj}}\left( d\Msigma|_{kj}\right)\1_{w_{kj}\geq 0}dw_{kj}, 
\end{align}
where $k\in \text{Active}(\Mq)$ such that $k\in e_{ij}$. For a map $\mm\supset\q$, we denote by $\mm^{\text{ext}(\q)}$ a copy of $\mm$ in which we add a vertex at the tip of each active edge of $\q$. Here $\mm^{\text{ext}(\q)}$ lets us divide the behavior of the Brownian bridges at the newly added vertices as in \Cref{eq:pivot}. Applying this to \eqref{e.law_weak_Markov} we get
\begin{align*}
    &\P^{\ell,b}\left[\BMM\in d\Mp, \MDec \in d\Msigma\left|\Mq\MLalg\BMM\right.\right] \\
    \propto &\left(\prod_{h\in H(\Mq)}q^{|V_i(\Mp_h)|} \prod_{e\in E(\Mp_h)}e^{-\lambda w_e}\right)\left(\prod_{e\in E(\Mq)}e^{-\lambda w_e}\right)\left(\prod_{k\in \text{Active}(\Mq)}d\sigma_k\right) \\
    &\left[ \sum_{\mm\supseteq\q} \delta_{\mm}(\q\glue(\p_h)_{h\in H(\Mq)})\left( \prod_{v\in V_i(\mm)}\mu(d\sigma_v)\right) \prod_{ij\in E(\mm^{\text{ext}(\q)})}\hat{\P}^{\sigma_i,\sigma_j}_{w_{ij}}\left( d\Msigma|_{ij} \right)dw_{ij}\right]\\
    \propto & \left(\prod_{ij\in  E(\q)}\hat{\P}^{\sigma_i,\sigma_j}_{w_{ij}}\left( d\Msigma|_{ij} \right)e^{-\lambda w_{ij}}dw_{ij}
\right) \left(\prod_{h\in H(\Mq)}q^{|V_i(\Mp_h)|} \prod_{e\in E(\Mp_h)}e^{-\lambda w_e}\right )\\
    &\left[ \left(\sum_{\mm_h}\delta_{\mm_h}(\p_h)\left(\prod_{v\in \text{Active}(h)}d\sigma_v\right) \left(\prod_{v\in V_i(\mm_h)}\mu(d\sigma_v)\right) \prod_{ij\in E(\mm_h)}\hat{\P}^{\sigma_i,\sigma_j}_{w_{ij}}\left( d\Msigma|_{ij} \right) dw_{ij}\right)\right],
\end{align*}

We now want to use Lemma \ref{l.Ley_condicional}. To do this, define the random variables $X = (\BMM_h^{\Mq},\MDec_h^{\Mq})_{h\in H(\Mq)}$ and $Z = \left((w_e)_{e\in E(\Mq)} , \MDec_{\Mq}\right)$ and the measures
\begin{align*}
    &\mu_X =\prod_{h\in H(\Mq)}\left(\sum_{\mm_h}\delta_{\mm_h}(\p_h) \left(\prod_{v\in \text{Active}(h)}d\sigma_v\right) \left(\prod_{v\in V_i(\mm_h)}\mu(d\sigma_v)\right)\prod_{ij\in E(\mm_h)}\hat{\P}^{\sigma_i,\sigma_j}_{w_{ij}}\left( d\Msigma|_{ij} \right)dw_{ij}\right),\\
    &\mu_Z = \left(\prod_{v\in V_i(\Mq)}\mu(d\sigma_v)\right)\left(\prod_{ij\in E(\Mq)}\hat{\P}^{\sigma_i,\sigma_j}_{w_{ij}}\left( d\Msigma|_{ij} \right)dw_{ij}\right).
\end{align*}
Furthermore, define 
\begin{align*}
    F(X,Z) = \left(\prod_{h\in H(\q)}q^{|V_i(\p_h)|}\prod_{e\in E(\Mp_h)}e
    ^{-\lambda w_e}\right)\prod_{e\in E(\Mq)}e^{-\lambda w_e}.
\end{align*} 
Then, Lemma \ref{l.Ley_condicional} implies that
\begin{align*}
    &\P^{\ell,b}\left[\BMM \in d\Mp, \MDec \in d\Msigma\left|\Mq\Lalg\BMM, \MDec_{\Mq} \in d\Msigma_{\Mq}  \right.\right]\\
    \propto & \prod_{h\in H(\q)}\left(q^{|V_i(\p_h)|}\prod_{e\in E(\Mp_h)}e^{-\lambda w_e}\sum_{\mm_h}\delta_{\mm_h}(\p_h)\left(\prod_{v\in V_i(\mm_h)}\mu(d\sigma_v)\right) \prod_{ij\in E(\mm_h)}\hat{\P}^{\sigma_i,\sigma_j}_{w_{ij}}\left(  d\Msigma|_{ij} \right)dw_{ij} \right).
\end{align*}

This directly implies the independence between the metric maps associated with each hole and that in each hole they have the distribution of a Boltzmann decorated metric map.
\end{proof}

\subsubsection{Filtration and stopping maps.}
We need a new notion of filtrations for metric maps. They are going to be collections of $\sigma$-algebras indexed by the equivalence classes of metric maps with respect to $\sim$ but we are going to denote it just by $[\Mq]$.

\begin{defn}
We say that $\FF=(\F_{[\Mq]})_{\Mq\in\MQm_H}$ is a \textit{filtration} if it satisfies the properties \ref{filtration-monotonicity}\footnote{with respect to the relation $\MLalg$.}, \ref{filtration-completesness} and \textit{Right Continuity},
\begin{align}
    \tag{Right Continuity}\label{filtration-continuity} \mbox{for any }[\Mq_n]\searrow[\Mq]\mbox{, we have that }\displaystyle\bigcap_{n\in\mathbb{N}}\F_{[\Mq_n]} = \F_{[\Mq]}.
\end{align}
\end{defn}

For the rest of this section, $(\BMM,\MDec)$ is going to be a Boltzmann decorated metric 4-regular map. 

\begin{defn}
We say that a Boltzmann decorated metric 4-regular map $(\BMM,\Dec)$ is an \textit{$\FF$-Boltzmann decorated metric 4-regular map} if it satisfies the following properties
\begin{itemize}[label={--}]
    \item\label{Adaptability-dec}\textit{Adaptability:} the event $\{\Mq\MLalg\BMM\}$ and the function $\Dec|_{\Mq}\1_{\Mq\MLalg\BMM}$ are $\F_{[\Mq]}$-measurable.
    \item\label{Indep-increments-dec}\textit{Independent increments:} conditionally on $\F_{[\Mq]}$ and $\MDec|_{\Mq}$, the law of $(\BMM_{h}^{\Mq},\Dec_{h}^{\Mq})_{h\in H(\Mq)}$ is that of a collection of independent $q$-Boltzmann decorated metric maps with boundary $h$ and boundary condition $\MDec|_{h}$ for every $h\in H(\Mq)$.
\end{itemize}
\end{defn}
\begin{rem}
    An example of filtration is the \textit{natural filtration} associated with a Boltzmann decorated metric map $(\BMM,\Dec)$ defined as follows,

\begin{align*}
    \F_{[\Mq]} := \overline{\bigcap_{\varepsilon >0}\bigvee_{\Mp\Lalg \Mq_{\varepsilon}} \sigma\left( \{\Mp \MLalg \BMM\} , \Dec|_{\Mp}\1_{\Mp\MLalg\BMM} \right)}^{\P}.
\end{align*}

The intersection with respect to $\varepsilon>0$ is made in order to satisfy \ref{filtration-continuity}. The previous $\sigma$-algebra is well defined, since the events only depend on the lengths of the active edges, which are the same regardless of the representative chosen from the equivalence class.
\end{rem}

With this, we are in condition to define the random metric submaps that are going to satisfy the Markov property. Unfortunately, we need to modify the definition of the previous sections. 
\begin{defn}
We say that a random metric map with hole vertices $\SMM$ is an \textit{$\FF-$stopping metric map} if
\begin{enumerate}
    \item $\P$-almost surely $\SMM\Lalg\BMM$,
    \item and, for any $\Mq\in\MQm_H$, we have that $\left\{\SMM\MLalg\Mq \right\}\cap \left\{\Mq\MLalg \BMM\right\} \in\mathcal{F}_{[\Mq]}.$
\end{enumerate}
\end{defn}

Let us note that for a random metric map with hole vertices $\SMM$, the events $\MQm_{\MLalg\Mq}$ characterise the $\sigma$-algebra generated by $\SMM$.
\begin{prop}
Let $\SMM$ be a random metric map and for $\Mq\in\MQm_H$ consider the events $\MQm_{\MLalg\Mq} :=\{\SMM\MLalg \Mq\}$. Then,
\begin{align}
    \sigma\left(\left[\SMM\right]\right) = \sigma\left(\MQm_{\MLalg\Mq} : \Mq\in\MQm_H\right).
\end{align}
\end{prop}
\begin{proof}
First denote by $\mathcal{G} = \sigma\left(\MQm_{\MLalg\Mq} : \Mq\in\MQm_H\right).$

\begin{itemize}
    \item For $\subseteq$: First notice that $\ske(\SMM)$ is $\mathcal{G}$-measurable. This comes from the fact that the events $\{\ske(\SMM)\Lalg\q\}\in\mathcal{G}$ for every $\q\in\mathcal{Q}_H$ and these events characterize the $\sigma$-algebra generated by $\ske(\SMM)$ (Proposition \ref{prop:sigma_algebra}). Next, we are left to prove that the event
    \begin{align*}
        \{\ske(\SMM)= \ske(\Mq),w_{e}^{\SMM}\leq r\}\in\mathcal{G}.
    \end{align*}
This comes directly from the fact that $\Mq$ takes a countable collection of skeletons. This implies that the lengths of all active edges are $\mathcal{G}$-measurable. We conclude that $\sigma([\SMM])\subseteq\mathcal{G}$.    
    \item For $\supseteq$: Notice that given an skeleton $\ske(\Mq)$, the equivalence class $[\Mq]$ is determined by the lengths of the active edges. Then, the event $\{\SMM\MLalg\Mq\}$ is fully determined by a countable number of skeletons and a finite number of inequalities of type $w_e^{\SMM}\leq w_e^{\Mq}$. In consequence, the event $\{\SMM\MLalg\Mq\}$ belongs to $\sigma(\left[\SMM\right])$ which implies $\mathcal{G}\subseteq \sigma([\SMM])$.
\end{itemize}
\end{proof}

For the proof of the strong Markov property we need a way to approximate submaps in a continuous way from above. 

\begin{defn}\label{def:aprox_active}
Let $\Mp\Lalg\Mq$ be two metric maps and let $\varepsilon>0$. We define the metric map with hole vertices $\Maprox{\Mp}{\Mq}{\varepsilon}\subseteq \Mq$ as follows. Starting from $\Mp$, for each edge $e\in\text{Active}(\Mp)$, we extend its explored length to the next multiple of $\varepsilon$, or until the endpoint of the corresponding edge in $\Mq$ is reached. More precisely, 
\begin{align*}
    w_e^{\Maprox{\Mp}{\Mq}{\varepsilon}} := \left( \varepsilon\left\lceil\frac{w_e^{\Mp}}{\varepsilon} \right\rceil \right) \wedge w_e^{\Mq}.
\end{align*}
When two extended active edges meet or overlap, we identify their common segment according to the identifications in $\Mq$. Whenever this procedure completely discovers an active edge, we reveal its endpoint and add all the edges incident to the newly discovered vertex that are not already contained in $\Mp$. These new edges are added with length zero and become active edges. In particular, when two active edges meet in the interior of the same explored edge, that edge becomes fully explored without revealing any additional vertex. The extension procedure is applied only to the edges that were already active in $\Mp$.
\end{defn}

\begin{rem}
    Notice that $\Maprox{\Mp}{\q}{\varepsilon}$ is different from $\Mp_{\varepsilon}$ because it only extends the lengths of the active edges and allows us to extend in a way that the submap obtained is still an active submap of $\Mq$. Furthermore, $[\Maprox{\Mp}{\q}{\varepsilon}]\MLalg [\Mq]$ but $[\Mp_{\varepsilon}]\not\MLalg [\Mq]$ in general.
\end{rem}

\begin{rem}
The previous approximation induces a third possible type of peeling step. Take a length parameter $L>0$ as input.  
\begin{itemize}
    \item\textbf{Peeling step of type 3:} We say that we obtain a peeling step of type 3 and parameter $L$, if starting from the root edge we discover the edge up to $L$ starting from the end and we still have not found any other vertex. Note that in this case, the information we have seen contains all the spins associated with this section of the edge. 
\end{itemize}
To obtain this type of peeling step, we explore an edge until we either have explored $L$ units of length or we have hit another vertex. We denote by $Peel(\tilde{\mathfrak{e_i}},L,e,\Mq)$ the resulting map of doing a peeling on $e \in \text{Active}(\Mq)$ of $\Mq$ from $\mathfrak{e_i}$ (See Figure \ref{fig:peeldagger})

\begin{figure}[H]
    \centering
    \includegraphics[scale=0.8]{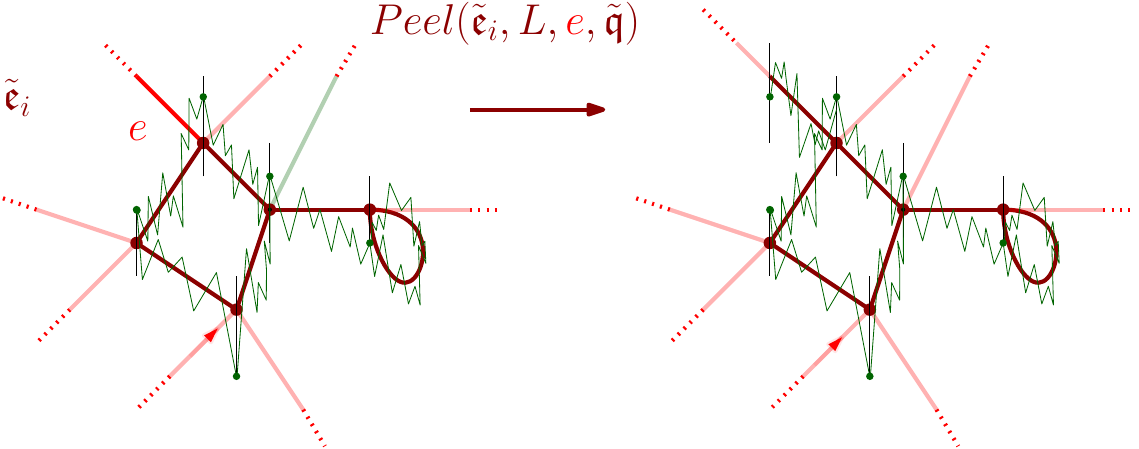}
    \caption{Representation of the peeling step for the metric map $\Mq$, starting from $\tilde{\mathfrak{e}}_i$, when selecting the edge $e$; i.e. when $L$ is smaller than the full length of the edge $e$. Note that this induces a peeling step of type 3. We added a representation of the decoration as an example to also illustrate that the spins are revealed.  }
    \label{fig:peeldagger}
\end{figure}
The Markov property induced by this type of peeling step turns the tip of the edge revealed by the peeling step of type 3 into phantom vertex with its decoration in order to induce a boundary condition. 
\end{rem}

Next, we present some properties associated with the stopping metric maps. 
\begin{prop}
Let $(\BMM,\Dec)$ be an $\FF$-Boltzmann decorated metric map.
\begin{enumerate}
    \item Let $\SMM_n$ be a decreasing sequence of $\FF-$stopping metric maps for $(\BMM,\Dec)$ such that $\SMM_n\searrow\SMM_{\infty}$ where $\SMM_{\infty}$ is metric map with hole vertices. Then $\SMM_{\infty}$ is an $\FF-$stopping metric map for $(\BMM,\Dec)$.
    \item Let $\SMM$ be an $\FF-$stopping metric map for $(\BMM,\Dec)$ and $n\in\N$, then, $\Maprox{\SMM}{\BMM}{2^{-n}}$ is an $\FF-$stopping metric map for $(\BMM,\Dec)$.
\end{enumerate}
\end{prop}
\begin{proof} 
Let $\Mq\in\MQm_H$.
    \begin{enumerate}
        \item Note the following 
        \begin{align*}
            \{\SMM_{\infty}\MLalg \Mq\} \cap \{\Mq\MLalg \BMM\}&= \underbrace{\bigcup_{K\in\N}\bigcap_{k\geq K}\underbrace{\bigcup_{N_0\in\N}\bigcap_{n\geq N_0}\left(\{\SMM_n\MLalg\Mq_{1/k} \}\cap \{\Mq_{1/k}\MLalg \BMM\}\right)}_{\in \F_{[\Mq_{1/k}]}}}_{\displaystyle\bigcap_{k\in \N}\F_{[\Mq_{1/k}]}},
        \end{align*}
where, thanks to the right continuity, we conclude.
\item Let $n\in\N$ and note the following
\begin{align*}
    \left\{ \Maprox{\SMM}{\BMM}{2^{-n}} \MLalg \Mq  \right\} \cap \{\Mq\MLalg\BMM\} &= \left\{ \SMM \MLalg \Maprox{\Mq}{}{-2^{-n}}  \right\} \cap \{\Mq\MLalg\BMM\},\\
\end{align*}
where $\Maprox{\Mq}{}{-2^{-n}}$ is defined as in \ref{def:aprox_active} but we take the floor function instead of the ceil function. In consequence the map obtained is an active metric submap of $\Mq$. Then,
\begin{align*}
    \left\{ \Maprox{\SMM}{\BMM}{2^{-n}} \MLalg \Mq  \right\} \cap \{\Mq\MLalg\BMM\} &= \left[\left\{ \SMM \MLalg \Maprox{\Mq}{}{-2^{-n}}  \right\} \cap \left\{\Maprox{\Mq}{}{-2^{-n}}\MLalg\BMM\right\}\right]\cap \{\Mq\MLalg\BMM\}.
\end{align*}
where, thanks to the property of transitivity, we can add the event $\Maprox{\Mq}{}{-2^{-n}}\MLalg\BMM$. Finally, notice that the first term belongs to $\F_{[\Maprox{\Mq}{}{-2^{-n}}]} \subseteq \F_{[\Mq]}$ and, since $\BMM$ it is an $\FF$-Boltzmann map, then $\{\Mq\MLalg\BMM\}\in\F_{[\Mq]}$. 
    \end{enumerate}
\end{proof}

For the following, we need to define the $\sigma$-algebra associated with the stopping metric maps.
\begin{defn}
The $\sigma$-algebra associated with an $\FF-$stopping metric map $\SMM$ is defined as follows
\begin{align}
    \F_{[\SMM]} := \left\{\left. \Theta \in \bigvee_{\Mq\in\MQm_H} \F_{[\Mq]} \right| \Theta\cap \left\{\SMM\MLalg\Mq\right\} \cap \left\{\Mq \MLalg \BMM\right\} \in \F_{[\Mq]}, \mbox{ for any }\Mq\in\MQm_H\right\}.
\end{align}
\end{defn}

For this $\sigma$-algebra we have the following properties.
\begin{prop}\label{prop:SMM-filtration}
    Let $(\BMM,\Dec)$ be an $\FF$-Boltzmann decorated metric map.
    \begin{enumerate}
        \item If $\SMM$ is an $\FF-$stopping metric map for $(\BMM,\Dec)$, then $\SMM$ is $\F_{[\SMM]}$-measurable.
        \item If $\SMM_1$ and $\SMM_2$ are two $\FF-$stopping metric maps for $(\BMM,\Dec)$ such that $\SMM_1 \MLalg \SMM_2$ almost surely, then $\F_{[\SMM_1]}\subseteq \F_{[\SMM_2]}.$
        \item If $(\SMM_n)_{n\in\N}$ is a sequence of $\FF-$stopping metric maps for $(\BMM,\Dec)$ such that $\SMM_n\searrow \SMM$, then
        \begin{align*}
            \displaystyle\bigcap_{n\in\N}\F_{[\SMM_n]} = \F_{[\SMM]}.
        \end{align*}
        \item For any real integrable random variable $X$, we have that    
        \begin{align*}
            \E\left[X\1_{[\SMM]=[\Mq]}|\F_{[\SMM]}\right] = \E\left[X|\F_{[\Mq]}\right]\1_{[\SMM]=[\Mq]}.
        \end{align*}
    \end{enumerate}
\end{prop}

\begin{proof}$\ $
\begin{enumerate}
    \item It suffices to prove that $\{\SMM\MLalg\Mq\}\in\F_{[\SMM]}$ for all $\Mq\in\MQm_H$. In fact, for all $\Mq_1\in\MQm_H$,
    \begin{align*}
        \{\SMM\MLalg\Mq\}\cap \{ \SMM\MLalg\Mq_1 \} \cap\{ \Mq_1\MLalg \BMM \} =\left [\{\SMM\MLalg\Mq\wedge \Mq_1\}\cap \{ \Mq_1\wedge \Mq\MLalg\BMM \}\right ]\cap\left[\{ \SMM\MLalg\Mq_1 \} \cap\{ \Mq_1\MLalg\BMM \}\right].
    \end{align*}
    where, $\Mq\wedge\Mq_1$ denotes the biggest metric map that is contained in both $\Mq$ and $\Mq_1$ with respect to $\MLalg$. Notice that the first term belongs to $\F_{[\Mq\wedge\Mq_1]}$ so,consequently, it belongs to $\F_{[\Mq_1]}$. The second term, from the definition of $\SMM$, also belongs to $\F_{[\Mq_1]}$.
    
    \item Let $\Theta\in\F_{[\SMM_1]}$. Then,
    \begin{align*}
        \Theta\cap\{ \SMM_2\MLalg\Mq \} \cap\{ \Mq\MLalg \BMM \} \stackrel{a.s.}{=} \left[\Theta\cap\{\SMM_1\MLalg\Mq \}\cap\{ \Mq\MLalg \BMM \} \right]\cap\{ \SMM_2\MLalg\Mq \}\cap \{\SMM_1\MLalg \SMM_2\}.
    \end{align*}
Notice that the first term, from the hypothesis, belongs to $\F_{[\Mq]}$. Also, the second term belongs to $\F_{[\Mq]}$, because $\SMM_2$ is an $\FF$-stopping map. The last term is the complement of a negligible event, and, in consequence, belongs to $\F_{[\Mq]}$.   
    
    \item Notice that the inclusion $\supset$ comes from the monotonicity property shown in (2) for stopping metric maps. On the other hand, if $\Theta\in\displaystyle\bigcap_{n\in\N}\F_{[\SMM_n]}$,
    \begin{align*}
        \Theta\cap\{ \SMM\MLalg\Mq \} \cap\{ \Mq\MLalg \BMM \} &=  \Theta\cap\{ \Mq\MLalg \BMM \}\cap \bigcup_{K\in\N}\bigcap_{k\geq K}\bigcup_{N_0\in\N}\bigcap_{n\geq N_0}\{ \SMM_n\MLalg\Mq_{1/k} \}\\
        &= \bigcup_{K\in\N}\bigcap_{k\geq K}\bigcup_{N_0\in\N}\bigcap_{n\geq N_0}\left[\Theta\cap\{ \Mq_{1/k}\MLalg \BMM \}\cap\{ \SMM_n\MLalg\Mq_{1/k} \}\right].
    \end{align*}
Notice that $\displaystyle\bigcup_{N_0\in\N}\bigcap_{n\geq N_0}\Theta\cap\{ \Mq_{1/k}\MLalg \BMM \}\cap\{ \SMM_n\MLalg\Mq_{1/k}\}\in \F_{[\Mq_{1/k}]}$ for every $k\in\N$, so, in consequence 
\begin{align*}
    \bigcup_{K\in\N}\bigcap_{k\geq K}\bigcup_{N_0\in\N}\bigcap_{n\geq N_0}\Theta\cap\{ \Mq_{1/k}\MLalg \BMM \}\cap\{ \SMM_n\MLalg\Mq_{1/k}\}\in\bigcap_{k\in\N}\F_{[\Mq_{1/k}]} = \F_{[\Mq]},
\end{align*}
where we use the right continuity of the filtration.
\item Let us show that the right-hand side satisfies the integral property of the conditional expectation. To do that take $\Theta\in\F_{[\SMM]}$ and compute
\begin{align*}
    \E\left[X\1_{[\SMM] = [\Mq]}\1_{\Theta}\right] = \E\left[\E\left[X\1_{[\SMM] = [\Mq]}\1_{\Theta}|\F_{[\Mq]}\right]\right] = \E\left[\E\left[X|\F_{[\Mq]}\right]\1_{[\SMM]=[\Mq]}\1_{\Theta}\right].
\end{align*}
\end{enumerate}
\end{proof}
\subsubsection{Strong Markov property.} In this section, we state the strong Markov property for the Boltzmann decorated metric map.
\begin{thm}\label{t.metric_strong_Markov}
    Let $(\BMM,\Dec)$ be a Boltzmann decorated metric map and $\SMM\in\MQm_H$ an $\FF-$stopping metric map for $(\BMM,\Dec)$. Then, $\BMM$ and $\Dec$ can be decomposed as follows,
    \begin{align}
        \BMM = \SMM\glue \left(\BMM_h^{\SMM}\right)_{h\in H\left(\SMM\right)},\  \text { and } \ \Dec = \Dec|_{\SMM} + \Dec^{\SMM},
    \end{align}
    where, conditional on $\F_{[\SMM]}$ and $\Dec_{\SMM}$, $\left(\BMM_h^{\SMM},\Dec^{\SMM}_h\right)_{h\in H\left(\SMM\right)}$ is a collection of independent Boltzmann decorated metric maps with boundary $\ell(h)$ and boundary condition $\Dec|_{h}$ for every $h\in H\left(\SMM\right)$.
\end{thm}

First let us prove the case where $[\SMM]$ takes a countable number of values.
\begin{lemma}\label{lem-countable}
    Theorem \ref{t.metric_strong_Markov} is true when $[\SMM]$ takes a countable values.
\end{lemma}
\begin{proof}
    Denote by $(\BMM_h^{\SMM})_{h\in H(\SMM)}$ the collection of maps $\BMM$ arising from $\SMM\Lalg\BMM$ almost surely. Let $[P]
    $ denote the set of all possible values of $[\SMM]$, which is countable. Therefore, for each possible value, we choose one representative. We denote by $P$ the set of representatives chosen for each possible value of $[\SMM]$. Notice that this does not imply that $\SMM$ itself takes only countably many values, but rather that its equivalence class does. Take $(f_h)_{h\in H(\SMM)}$ a collection of measurable functions taking values in $[0,1]$ (measurable with respect to $H(\SMM)$). Then,
    \begin{align*}
        \E\left[\prod_{h\in H(\SMM)}f_h\left(\BMM_h^{\ \SMM},\MDec_h^{\ \SMM} \right)\right] &= \sum_{\Mq\in P}\E\left[\left.\prod_{h\in H(\SMM)}f_h\left(\BMM_h^{\ \SMM},\MDec_h^{\ \SMM} \right)\1_{\left[\SMM\right] = [\Mq]}\right|\F_{[\SMM]}\right]\\
        &= \sum_{\Mp\in P}\prod_{h\in H(\Mq)}\E\left[\left.f_h\left(\BMM_h^{\ \Mq},\MDec_h^{\ \Mq} \right)\right|\F_{[\Mq]}\right]\1_{\left[\SMM\right] = [\Mq]}.
    \end{align*}
We conclude since $\{[\SMM] = [\Mq]\}$ is $\F_{[\SMM]}$-measurable (from \ref{prop:SMM-filtration} (1)) and that conditionally on $\F_{[\Mq]}$ the metric decorated map $(\BMM_h^{\ \Mq},\MDec_h^{\ \Mq})$ has the law of a Boltzmann decorated metric map for any $\Mq\in P$, this is nothing else than the weak Markov property.
\end{proof}

Now, we can begin with the proof of the strong Markov property.

\begin{proof}[Proof of Theorem \ref{t.metric_strong_Markov}]
    For the general case, let $n\in\N$ and consider the following metric map
    \begin{align*}
        \SMM_n := \Maprox{\SMM}{\BMM}{2^{-n}}.
    \end{align*}
    The sequence of metric maps $\left(\SMM_n\right)_{n\in\N}$ satisfies the following 
    \begin{enumerate}
        \item $\left[\SMM_n\right]$ takes countably many values, for any $n\in\N$,
        \item $\SMM_n$ is an $\FF$-stopping metric map for $\BMM$ for any $n\in\N$,
        \item $\left[\SMM_n\right]\searrow [\SMM]$,
        \item and, for a large enough $n$, $\ske(\SMM_n)= \ske(\SMM)$.
    \end{enumerate}
    By $(1)$ and $(2)$, $\SMM_n$ satisfies the strong Markov property. Thus, we have that
    \begin{align*}
        \BMM &= \SMM_n\glue \left(\BMM_h^{\SMM_n}\right)_{h\in H\left(\SMM_n\right)},\  \text { and } \ \Dec = \Dec_{\SMM_n} + \Dec^{\SMM_n},
    \end{align*}
    where, conditional on $\F_{\left[\SMM_n\right]}$ and $\Dec_{\SMM_n}$, $\left(\BMM_h^{\SMM_n},\Dec^{\SMM_n}_h\right)_{h\in H\left(\SMM_n\right)}$ is a collection of independent Boltzmann decorated metric maps with boundary $\ell(h)$ and boundary condition $\Dec|_{h}$ for every $h\in H\left(\SMM_n\right)$. For $F$ a real bounded measurable function we define the following
    \begin{align*}
        M_n := \E\left[ F\left(\BMM,\Dec\right)|\F_{[\SMM_n]}\right].
    \end{align*}
    The sequence $(M_n)_{n\in\N}$ is a backward martingale with respect to the filtration $\left(\F_{[\SMM_n]}\right)_{n\in\N}$. Also, by Proposition \ref{prop:SMM-filtration}, we have that $\F_{\left[\SMM_n\right]}\searrow \F_{[\SMM]}$. In consequence 
    \begin{align*}
        \lim_{n\to\infty}M_n = \E\left[ F\left(\BMM,\Dec\right)|\F_{[\SMM]}\right].
    \end{align*}
    Now, for $n$ large enough, consider the following function
    \begin{align*}
        F_n(\BMM,\Dec) &= \prod_{h\in H(\SMM)}f_h(\BMM^{\SMM_n}_h,\Dec_h^{\SMM_n}),
    \end{align*}
    where $(f_h)_{h\in H(\SMM)}$ is a sequence of real bounded measurable functions. Then,
    \begin{align*}
        \E\left[ F_n\left(\BMM,\Dec\right)|\F_{[\SMM_n]}\right] &= \prod_{h\in H(\SMM)}\E^{\ell(h),b_{n,h}}[f_h(\BMM^{\SMM_n}_h,\Dec^{\SMM_n}_h)],
    \end{align*}
    where $b_{n,h}$ denotes the boundary condition of the hole $h$ induced by $\SMM_n$. Thanks to Proposition \ref{p.continuity_metric}, we have that
\begin{align*}
    \lim_{n\to\infty}\E^{\ell(h),b_{n,h}}[f_h(\BMM^{\SMM_n}_h,\Dec^{\SMM_n}_h)] = \E^{\ell(h),b_{h}}[f_h(\BMM^{\SMM}_h,\Dec^{\SMM}_h)],
\end{align*}
    where $b_h$ denotes the boundary condition of the hole $h$ induced by $\SMM$. Mixing all together we conclude that the law of $\left(\BMM_h^{\SMM},\Dec^{\SMM}_h\right)_{h\in H\left(\SMM\right)}$ is the law of a collection of independent Boltzmann decorated metric map with boundary $h$ and boundary condition $\Dec|_{h}$ for every $h\in H\left(\SMM\right)$.
\end{proof}

\subsection{Local maps are stopping maps.}
In this subsection, we characterise all random metric submaps of a Boltzmann decorated metric map that satisfy the Markov property. Informally, we show that the only ones are stopping metric maps. To do this, we first need to define the concept of a local metric map, namely, the metric maps that induce a Markovian decomposition.
\begin{defn}
We say that the triplet $\left(\BMM,\MDec,\SMM\right)$ is a \textit{local metric map} if,
\begin{itemize}
    \item $(\BMM,\MDec)$ is a Boltzmann decorated metric map;
    \item $\P$-almost surely $\SMM\Lalg\BMM$;
    \item and, conditionally on $\SMM$ and $\MDec_{\ \SMM}$, $(\BMM_h^{\ \SMM},\MDec^{\ \SMM}_h)_{h\in H\left(\SMM\right)}$ is a collection of independent Boltzmann decorated metric maps with perimeter $\ell(h)$ and boundary condition $\MDec|_h$ for every $h\in H\left(\SMM\right)$. 
\end{itemize}
\end{defn}
The previous definition is the generalization of \textit{local maps} defined in \cite{araya2025characterisation} for planar maps but adapted to metric maps. The next result shows that the only metric maps that induce the Markov property are the stopping metric maps.
\begin{thm}\label{t.local_implies_stopping}
    Let $\left(\BMM,\MDec,\SMM\right)$ be a local metric map. Then, there exists a filtration $\FF = (\F_{[\Mq]})_{\Mq\in\MQm_H}$ such that $(\BMM,\MDec)$ is an $\FF-$Boltzmann decorated metric map and $\SMM$ is an $\FF-$stopping metric map for $(\BMM,\MDec)$.
\end{thm}

Almost all of the proof made in Theorem 3.2 in \cite{araya2025characterisation} for local maps is analogous for this case. So, we are going to mention the main differences and how each lemma changes for this case but we leave the proof to the reader. We start defining the filtration that gives the result. For $\Mq\in\MQm_H$ we define the following $\sigma$-algebra
\begin{align}\label{e.def_filt_metric}
\F_{[\Mq]} = \bigcap_{\varepsilon>0}\ \bigvee_{\Mq^{\varepsilon}\in \Maprox{\Mq}{\BMM}{\varepsilon}}\ \bigvee_{\Mp\MLalg\Mq^{\varepsilon}} \overline{\sigma\left(\{\Mp\MLalg\BMM\},\{\SMM\MLalg \Mp\}, \{\MDec|_{\Mp}\1_{\Mp\MLalg\BMM}\}, \{\SMM^{\wedge \Mp}\1_{\Mp\MLalg\BMM}\} \right)}^{\P}.
\end{align}
Here by $\Mq^{\varepsilon}\in \Maprox{\Mq}{\BMM}{\varepsilon}$ we mean that there is a $\BMM\supset \Mq$ such that $\Mq^{\varepsilon}$ can be obtained by an $\varepsilon$ increasing of $\Mq$ in $\BMM$. Clearly $\SMM$ is an $\FF$-stopping metric map and $(\BMM,\MDec)$ satisfies the property of Adaptability. The only thing that is left to prove is the property of Independent increments. To prove that property we give the analogous lemmas 

\begin{lemma}\label{lem.local_stopping_1}
Let $\left(\BMM,\MDec,\SMM\right)$ be a local map. Then for $\Mq\in\MQm_H$ and a collection $(f_h)_{h\in H(\Mq)}$ of real, measurable, bounded functions, we have
\begin{align*}
    &\E\left[\left.\prod_{h\in H(\Mq)}f_h(\BMM_h^{\ \Mq},\MDec_h^{\ \Mq})\right|\SMM,\MDec|_{\Mq},\{\SMM\MLalg\Mq\MLalg\BMM\}\right]\1_{\SMM\MLalg\Mq\MLalg\BMM} = \prod_{h\in H(\Mq)}\E^{\frac{|h|}{2},\MDec|_{ h}}\left[f_h(\BMM,\MDec)\right ] \1_{\SMM\MLalg\Mq\MLalg\BMM}.
\end{align*}
Here, the pair $(\BMM_h^{\ \Mq},\MDec_h^{\ \Mq})$ is the one that appears from the weak Markov property of $(\BMM,\MDec)$.
\end{lemma}
\begin{lemma}\label{lem.local_stopping_2}
Let $\left(\BMM,\MDec,\SMM\right)$ be a local map and take $\Mq\in\MQm_H$. Additionally, fix $f_{\Mq}, (f_h)_{h\in H(\Mq)}, (g_{\Mp})_{\Mp\MLalg\Mq}$ and $(h_{\Mp})_{\Mp\MLalg\Mq}$ a collection of real, measurable and bounded functions and $\Lambda(\Mq)$ a collection of submaps of $\Mq$. Then, we have the following
\begin{align}\label{e.stopping_2}
    &\E\left[f_{\Mq}(\Mq,\MDec|_{\Mq})\1_{\SMM\MLalg\Mq\MLalg\BMM}\left (\prod_{\Mp\in \Lambda(\Mq)}\1_{\SMM\MLalg\Mp\MLalg\BMM}g_{\Mp}(\MDec|_{\Mp})h_{\Mp}(\SMM^{\wedge\Mp})\right ) \left( \prod_{h\in H(\Mq)}f_h\left(\BMM_h^{\Mq},\Dec_h^{\Mq}\right )\right)\right] \\
   \nonumber &= \E\left[f_{\Mq}(\Mq,\MDec|_{\Mq})\1_{\SMM\MLalg\Mq\MLalg\BMM}\left (\prod_{\Mp\in \Lambda(\Mq)}\1_{\SMM\MLalg\Mp\MLalg\BMM}g_{\Mp}(\MDec|_{\Mp})h_{\Mp}(\SMM^{\wedge\Mp})\right ) \left( \prod_{h\in H(\Mq)}\E^{\frac{ |h| }{2 },\MDec|_{h}}\left[f_h(\BMM,\MDec)\right]\right )\right].
\end{align}
\end{lemma}

\section{Characterisation of Markovian laws on decorated metric maps}\label{sec:4}
We present a characterisation of all SDMM satisfying the Markov property. This section is technical so on a first reading, we recommend that the reader skip and assume as true the Lemmas \ref{l.absolute_continuity_metric} and \ref{l.continuity_3_1} that prove the absolute continuity of the lengths of the edges and each value of the decoration with respect to the Lebesgue measure and Lemmas \ref{l.continuity_3_1}, \ref{l.peleeing_decom} and \ref{l.continuity_peleing_1_metric} which gives the continuity of the peeling transitions with respect to the value of the decoration at the longest of the edges. The only result you need for the following subsection is the part of Lemma \ref{l.peleeing_decom} which shows that the peeling step of type 3 has an exponential form and the peeling steps of type 1 and 2 are independent from the peeling step of type 3.

We need to state formally the Markov property for a general collection of measures $(\P^{\ell,b})_{(\ell,b) \in \N^{*}\times\Bc}$.

\begin{defn}\label{def:markov_property_metric}
We say that $(\P^{\ell,b})_{(\ell,b) \in \N^{*}\times\Bc}$ satisfies the \textit{Markov property}, if, for any deterministic metric 4-regular map with hole vertices $\Mq$, we can describe the conditional law $\P^{\ell,b}(\cdot|\{\Mq\Lalg\BMM\})$ as     \begin{align}
        \BMM = \Mq \glue (\BMM_h^{\Mq})_{h\in H(\Mq)},\  \text { and } \ \MDec = \MDec|_{\Mq} + \MDec^{\Mq}.
    \end{align}
    Here, conditional to $\MDec|_{\Mq}$, $(\BMM^{\Mq}_h,\MDec^{\Mq}_h)_{h\in H(\Mq)}$ is a collection of independent decorated metric maps with law $\P^{\ell(h),\MDec|_{h}}$ for every $h\in H\left(\Mq\right)$.
\end{defn}
                    
Now, we are ready to state the main result of this section that characterises all the decorated metric 4-regular maps having the Markov property. 

\begin{thm}\label{t.characterization_metric}
Take a collection of measures $(\P^{\ell,b})_{(\ell,b) \in \N^{*}\times\Bc}$ on spin-decorated metric 4-regular maps with a half-perimeter $\ell$ and boundary condition equal to $b:\cro{0,2\ell-1}\mapsto \R$ (and $\beta=1$ for simplicity). Furthermore, assume that for all $\ell\in \N$,  $$\lim_{\epsilon \to 0} \inf_{b}\P^{\ell,b}\left(\inf_{e \in E(\BMM)}w_{e}>\varepsilon\right)=1\text{ and }\P^{\ell,b}(\BM \text{ has at least one interior vertex})>0$$ and that the measure $\P^{\ell,b}$ is continuous in $b\in \R^\ell$ for the weak topology of measures. Then, the following are equivalent.
\begin{enumerate}
    \item There are $q,\lambda>0$ such that for any $\ell$ and $b$, $\P^{\ell,b}$ is a $(q,\lambda)$-Boltzmann decorated metric map,
    \item $(\P^{\ell,b})_{(\ell,b) \in \N^{*}\times\Bc}$ has the Markov property, and for each $\ell$ and $b$, the measure $\P^{\ell,b}$ has the Gibbs distribution on $\MQm^{\ell,v}$ for every $v$ and is invariant under rerooting,
    \item $(\P^{\ell,b})_{(\ell,b) \in \N^{*}\times\Bc}$ has the Markov property and for each $\ell$ and $b$ the measure $\P^{\ell,b}$ is invariant under rerooting.
\end{enumerate}
\end{thm}

Before the proof of this result, we need a few lemmas. We start by noting that the set of metric maps with $\ske(\tilde \q)=\q$ is in bijection with $(\R^+)^{\sharp E(\q)}$. We denote this bijection $w(\tilde \q)$. Now, we see that the induced law on $(\R^+)^{\sharp E(\q)}$ is absolutely continuous with respect to Lebesgue.
\begin{lemma}\label{l.absolute_continuity_metric}
Let $(\P^{\ell,b})_{(\ell,b)\in\N\times\Bc}$ be as in Theorem \ref{t.characterization_metric} and suppose that it satisfies the Markov property. Then, on the event that $\ske(\BMM) = \ske(\Mq)$, the law of $w(\BMM)$ is absolutely continuous with respect to $\Leb^{|E(\Mq)|}$.  From now on, we define
\begin{align*}
p^{\ell,b}(\tilde \q):=\frac{\P^{\ell,b}\left(\BMM\in d\Mq\right) }{ \Leb^{|E(\q)|}(dw)}(w(\tilde \q))
\end{align*}
\end{lemma}
\begin{proof}
    We order the edges $(e_{k})$ of $\ske(\Mq)$ in such a way that they could be discovered via a peeling procedure, and define $w_k$ the length of the edge $e_k$ in $\BMM$. It suffices to show that for any $k$, the conditional law of $w_k$ given $(w_j)_{j<k}$ is absolutely continuous with respect to Lebesgue.
    
    Take $k\in \N$ and assume that we have peeled in order the edges of $e_k$; in particular this implies that we are on the event $\q_{k-1}\subseteq \BMM$ for a given deterministic $\q_{k-1}$. Thus, we have that the conditional law of $w_k$ given $\q_{k-1}\subseteq \BMM$ and $\phi|_{\q_{k-1}}$ is the law of the length of a metric chose according to $\P^{\ell_{k-1},b_{k-1}}$ for a given semi-perimeter $\ell_{k-1}$ and (random) boundary condition $b_{k-1}$. This implies that to prove the lemma it is enough to see that for any boundary condition, $w_1$ the length of a given edge from the root is absolutely continuous with respect to Lebesgue, without  any need to further condition on $\ske(\BMM)=\ske(\tilde \q)$ as this is an event of positive probability.
    
    We start by using the following claim.
    \begin{claim}\label{c.constructing_w}
     For any $\delta>0$ there exists a triplet $(z,\tau, \tau')$, such that conditionally on a decoration $z$, the lengths $\tau$ and $\tau'$ are independent and there exists a deterministic set $E$ such that
    \begin{align}
        \label{e.law_as_a_sum}w_1\stackrel{law}{=}\begin{cases} \tau + \tau' & \text{ if } z\in E,\\
        \tau & \text{ if } z \notin E.
        \end{cases}
    \end{align}
    Furthermore, $\P(z\notin E)\leq \delta$, and $\P(\tau\in A\mid z=x)=0$ for any $x\in E$ and any set $A$ of Lebesgue measure $0$.
    \end{claim}
    
    Let us first see how to use the claim to conclude. Take a set $A$ of Lebesgue measure $0$. We bound
    \begin{align}
        \P(w_1\in A) &\leq \delta + \P(\tau + \tau' \in A, z \in E)\leq  \delta + \E\left[\P\left (
        \tau \in A-\tau' \mid z,\tau' \right ) \mathbf 1_{z \in E}\right ]\leq \delta.        
    \end{align}
    Thus, $\P(w_1\in A)=0$ which implies that the law of $w_1$ is absolutely continuous with respect to Lebesgue.
    
    We, now, just need to prove the claim.
    \begin{proof}[Proof of Claim \ref{c.constructing_w}] We need to discover a small part of an edge of the boundary (not necessarily the root itself), let us call it $j$ and its associated boundary $b(j)$.  Fix $\epsilon>0$ and define the following ``stopping time'' in the edge associated with $j$
    \begin{align*}
        \tau_\epsilon :=\inf\{t>0: \MDec(t)=b(j)\pm \epsilon\}\wedge \inf\{t>0: t \text{ is a vertex of }\BMM\},
    \end{align*}
    where we are abusing the notation by calling $t$ the point in the edge $e_k$ that is at distance $t$ from the original phantom vertex $j$. Informally, $\tau_\epsilon$ is continuously exploring $e_k$ until we either hit a vertex or the decoration hits $b(j)\pm\epsilon$. We now define $z=\MDec(\tau_\epsilon)$ and $E=\{b(j)\pm\epsilon\}$. It is clear that \eqref{e.law_as_a_sum} holds, and that conditionally on $z$, $\tau_\epsilon$ is independent of $\tau'$ the remainder of the length of $w_1$. We are left to prove two things
    \begin{itemize}
        \item \textit{It is unlikely that $z\notin E$.} As we know that $\P^{\ell,b}(w_1=0)=0$, we know that there is $L>0$ such that $\P^{\ell,b}(w_1>L)\geq 1-\delta/2$. We will choose $\epsilon$ small enough such that $\P(z\notin E)\leq \delta$. To do this, we call $v$ the vertex that is at the other end of $e_1$ and condition on $w_1$ to see that
        \begin{align*}
            \P^{\ell,b}(\tau_\epsilon\geq L\mid L, \MDec_v) &\leq  \delta/2+ \E^{\ell,b}\left[ \P^{b(j),\MDec_{\phi_1}}_{w_1}\left ((P_t)_{t\in [0,w_1]}\subseteq [b(j)-\epsilon,b(j)+\epsilon]\right ) \1_{w_1\geq L} \right ]\\ 
            &\leq \delta/2 + \frac{2\sqrt{2}\epsilon }{\sqrt{\pi L} }.
        \end{align*}
        Which gives the result, as a.s. the event $z\notin E$ is contain on the event $w_1<L$ union $\tau_\epsilon\geq L$.
        \item \textit{$\tau_\epsilon$ does not put a lot of mass on $0$-measure sets.} Take a set $A$ with $0$ Lebesgue measure. Note that conditionally on $w$, $\MDec(v)$ and $z \in E$ the law of $\tau$ is absolutely continuous with respect to Lebesgue (as it is that of the first hitting time of $b(j)\pm \epsilon$ of a Brownian bridge of length $w$ going from $b(j)$ to $\MDec(v)$ that hits exactly at $z$). Thus, for any $x \in \{b(j)\pm\epsilon\}$
        \begin{align*}
        \P^{\ell,b}\left(\tau\in A\mid z=x \right)=\E\left[\P\left( \tau \in A \mid w_1,z=x,\MDec_v\right ) \right]  = 0.
        \end{align*}

    \end{itemize}
        
    \end{proof}
\end{proof}

Assume that $(\P^{\ell,b})_{(\ell,b)\in\N^{*}\times\Bc}$ is like in the theorem and satisfies the Markov property. Then, if we take a metric 4-regular map $\Mq$ with half-perimeter $\ell$, the law of $\Dec$ is absolutely continuous with respect to $\mu^{|V_i(\Mq)|}\times \text{Leb}^{|E(\Mq)|}$ on the event that the map structure of $\BM$ is equal to the map structure of $\Mq$.
\begin{align}\label{e.def_p_q_s_metric}
    p^{\ell,b}(\Mq,\sigma) &= \frac{d\P\left(\BMM \in d\Mq,\Dec\in d\sigma\right)}{\mu^{|V_i(\q)|}\times \text{Leb}^{|E(\q)|}(d\sigma,dw)}(\sigma,w(\Mq)) =  p^{\ell,b}(\tilde \q)\frac{\exp\left(-\Ham^{b}(\Mq,\sigma)\right)}{Z^{b}(\Mq)}.
\end{align}
We are going to identify peeling steps of different types.
\begin{itemize}
    \item We say that we have a peeling step of type $1$ if $[\vroot,\vertex]\subseteq \BM$, where we denote by $[\vroot,\vertex]$ to the map with hole vertices containing an external vertex of degree $2\ell$ and the root is connected to a vertex $\vertex$ of degree $4$ with all the active edges of length $0$.
    \item The peeling step of type $3$ with length $t>0$ is identified with $[0,t]\subseteq \BMM$, where $[0,t]$ is the map with hole vertices containing an external vertex of degree $2\ell$ and the root edge is active and has length $t$, all the other vertices have length $0$.
\end{itemize}

\begin{lemma}\label{l.continuity_3_1}
Let us work in the context of Lemma \ref{l.absolute_continuity_metric} and fix $t>0$, on the event $[0,t]\subseteq \BMM$, the law of $\MDec_t$ is absolutely continuous with respect to Lebesgue. Additionally,
\begin{align}\label{e.def_p_3_metric}
p_3^{\ell,b} (t,x) :=  \frac{\P^{\ell,b}([0,t]\subseteq \BMM, \MDec_t \in dx)}{\Leb(dx)} (x)
\end{align}
is continuous\footnote{As in the chapter before, $p_3^{\ell,b}$ is defined a.e., thus we are saying that there is a continuous representative in $L^1(\R)$.} in both $t$ and $x$.

Additionally, on the event $[\vroot,\vertex]\subseteq \BMM$, $(w_{\vroot\vertex}, \phi_\vertex)$ is absolutely continuous with respect to Leb$\times \mu$. Then, for a.e. $w$
\begin{align}\label{e.def_p_1_metric}
p_1^{\ell,b}(w,x):= \frac{\P^{\ell,b}\left ([\vroot,\vertex]\subseteq \BM,w_{\vroot,\vertex}\in dw, \MDec_\vertex \in dx \right ) }{\Leb\times \mu (dw, dx) } (w,x),
\end{align}
is continuous in $x$.
\end{lemma}
\begin{proof}We start by proving the statements related to $p_3^{\ell,b}$. As conditionally on the event $\{[0,t]\subseteq \BMM\}$, the probability that the length of the root edge is bigger than $t+\epsilon$ goes to $1$ as $\epsilon \to 0$, and conditionally on that event $\MDec_t$ is absolutely continuous with respect to Lebesgue, we see that \eqref{e.def_p_3_metric} is well defined. Now, fix $b$ and $\ell$; there is a map $\q$ such that $\P^{\ell,b}(\ske(\BMM)=\q)>0$.  Calling $e$ the root edge we have that
\begin{align*}
\P^{\ell,b}(\ske(\BMM)=\q, w_e>2t)>0,
\end{align*}
thus, we can see that a.e. on $x$
\begin{align*}
\frac{d\P^{\ell,b}(\ske(\BMM)=\q, w_e>2t,\MDec_t\in dx)}{\text{Leb}(dx) } (x)= p_3^{\ell,b}(t,x) \P^{\ell,b_x}(\ske(\BMM)=\q,w_e>t),  
\end{align*}
where $b_x$ is equal to $b$ except in $0$ where it takes value $x$. As $\P^{\ell,b_x}(\ske(\BMM)=\q,w_e>t)$ is continuous in $x$ and $t$ positive when $x=b(0)$, we are only left to prove that the left-hand side is also continuous in $x$ and $t$. This is proven by dominated convergence as $\BMM$ is an SDMM and thus in the root edge the decoration behaves like a Brownian bridge, when conditioned on the value at $2t$. Then,
\begin{align*}
\frac{d\P^{\ell,b}(\ske(\BMM)=\q, w_e>2t,\MDec_t\in dx)}{\text{Leb}(dx) } (x) &= \P^{\ell,b}\left(\ske(\BM)=\q \right)  \E^{\ell,b}\left[ \frac{e^{\frac{-(x-l(t))^2 }{2 t } }}{\sqrt{ 2\pi t} }\mathbf 1_{w_e>2t}\mid \ske(\BMM)=\q\right] , 
\end{align*}
where $l(t)$ is the linear interpolation at $t$ between the boundary condition of the root edge and the value $\MDec_{2t}$.
Note that the term inside the expected value is continuous in $x$ and bounded by $(2\pi)^{-1/2}t^{-1} \mathbf 1_{w_e>2t}$. Thus we can conclude continuity on $t$ and $x$ by the bounded convergence theorem.    

Now, we work with $p_1$, again it is direct to see that \eqref{e.def_p_1_metric} is well posed as in any finite skeleton $\q$ with all lengths bigger than $\epsilon>0$, $\MDec_v$ is absolutely continuous with respect to $\mu$. Fix $\ell$, $b$ and $z\in \R$, and define the graph $\r= [\vroot,\vertex]\glue \q$, where $\q$ is such that $\P^{\ell+1,\hat b_z}(\q)>0$ where $\hat b_z$ takes value $z$ in $0$, $1$ and $2$ and values $b(j-2)$ at $j\geq 3$. We note that
\begin{align*}
\frac{d\P^{\ell,b}\left(\ske(\BMM)=\r, w_{\vroot\vertex}\in dw,\MDec_\vertex\in dx,\displaystyle\min_{e\in \BMM} w_{e}>\epsilon\right)}{\text{Leb}\times \mu(dw,dx) } (w,x)= p_1(w,x) \P^{\ell,\hat b_x}\left(\ske(\BMM)=\q,\min_{e\in \BMM} w_{e}>\epsilon \right ).
\end{align*}
As $\P^{\ell,\hat b_x}\left(\ske(\BMM)=\q \right )$ is continuous in $\hat b_x$ and positive when $x=z$, we just need to show that the left-hand side is continuous in $x$. We see that
\begin{align*}
    &\frac{d\P^{\ell,b}(\ske(\BMM)=\r, w_e\in dw,\MDec(v)\in dx)}{\text{Leb}\times \mu (dw,dx) } (w,x) \\
    = &\int\frac{d\P^{\ell,b}(\ske(\BMM)=\r, w_{\vroot\vertex}\in dt, w_{e}\in dt_e \ \forall e\neq \vroot\vertex)}{\text{Leb}^{\sharp E(\r)} (\prod_{e\in E(\r)}dw_e)}\frac{\P^{\ell,b}(\MDec_\vertex\in dx\mid \BMM=\r_{w,(w_e)}) }{ \mu(dx) }\prod_{e\neq \vroot\vertex} dw_e\1_{w_e\geq \epsilon},
\end{align*}
where $\r_{w,(w_e)}$ is the metric map with skeleton $\r$, length $w$ in the edge $\vroot\vertex$ and length $w_e$ in each edge different from $\vroot\vertex$.  As $\frac{\P^{\ell,b}(\MDec_\vertex\in dx\mid \BMM=\r_{w,(w_e)}) }{ \mu(dx) }$ is bounded and continuous in $x$ (as long as $w_e\geq \epsilon$), we conclude by the bounded convergence theorem.
\end{proof}

We are missing one key possibility: the event where the edge of the root is connected directly to the outer vertex twice. This is the so-called a peeling step of type 2. In this case, it is also clear that the length of the edge is absolutely continuous with respect to Lebesgue, as the event of having a peeling step of type $2$ has positive probability. Thus, we can define
\begin{align*}
p_{2}^{\ell,b,\ell_1,\ell_2}(w) &=\frac{d\P^{\ell,b}(\type{2,\ell_1,\ell_2, w_e\in dw})}{\text{Leb}(dw)}(w),
\end{align*}
where $\ell_1+\ell_2=\ell-1$.

Now we can compute the relationship between all the types of peeling steps and compute their rates.
\begin{lemma}\label{l.peleeing_decom}
Let us work in the context of Lemma \ref{l.absolute_continuity_metric}. Then, there exist positive constants $r^{\ell,b},s^{\ell,b,\ell_1,\ell_2}$ and $\lambda^{\ell,b}$ such that for a.e. $t$
\begin{align}
\label{e.formula_p3}p_3^{\ell,b}(t)=p_{3}^{\ell,b}(t,b(0)) &= \frac{1}{\sqrt{2\pi t}}e^{-\lambda^{\ell,b}t}\\
 \label{e.formula_p1}    p_{1}^{\ell,b}(t,b(0)) &= p_3^{\ell,b}(t) r^{\ell,b},\\
    p_{2}^{\ell,b,\ell_1,\ell_2}(t) &\label{e.formula_p2}=p_3^{\ell,b}(t) s^{\ell,b,\ell_1,\ell_2}\exp\left(-\frac{(b(0)-b(k))^2}{2t}\right),
\end{align}
where $b(k)$ is the boundary value of the vertex discovered by the peeling step of type 2 at the other side of the root edge. In particular, $p_1^{\ell,b}(t,b(0))$ and $p_{2}^{\ell,b,\ell_1,\ell_2}(t)$  are continuous in $t$.
\end{lemma}
\begin{proof}
Fix $t>0$ and $t_1+t_2=t$. Take $\q$ a map such that $\P^{\ell,b}(\ske(\BMM)=\q)>0$. For a sequence of lengths  $w=(w_e)_{e\in E(\q)}$ we define $\q_{w}$ as the metric map with skeleton $\q$ and length $w_e$ for the edge $e$. Denote by $1$ the root edge, and assume that $w_1>t$. Define $\hat w_e=w_e$ for all $e\neq 1$ and $\hat w_1= w_1-t$. We have that
\begin{align}\label{e.1_for_3}
    \frac{d\P^{\ell,b}( \ske(\BMM)=\q,w(\BMM) \in dw, \Dec_t \in dx )}{d\text{Leb}^{\sharp E(\q)}\times \text{Leb}}(w,b(0))=p_3^{\ell,b}(t,b(0))\frac{d\P^{\ell,b}( \ske(\BMM)=\q,w(\BMM) \in dw) }{d\text{Leb}^{\sharp E(\q)}}(\hat w,b(0))>0.
\end{align}
Now, instead of looking at the values at $t$, we can also look the value at $t_1$ and use the fact that $\BMM$ is a spin-decorated map to see that
\begin{align}\label{e.2_for_3}
    &\frac{d\P^{\ell,b}( \ske(\BMM)=\q,w(\BMM) \in dw,\phi_t \in dx,\Dec_{t_1}\in dy )}{d\text{Leb}^{\sharp E(\q)}\times \text{Leb}}(w,b(0),y)\\ \nonumber= &\frac{d\P^{\ell,b}( \ske(\BMM)=\q,w(\BMM) \in dw,\Dec_t \in dx,\Dec_{t_1}\in dy )}{d\text{Leb}^{\sharp E(\q)}\times \text{Leb}}(w,b(0),b(0))\exp\left(-\frac{(b(0)-y)^2}{2t_1}-\frac{(y-b(0))^2}{2t_2}\right)
    \\\nonumber=&
    \ p_3^{\ell,b}(t_1,b(0))p_3^{\ell,b}(t_2,b(0))\exp\left(-\frac{(b(0)-y)^2}{2t_1}-\frac{(y-b(0))^2}{2t_2}\right)\frac{d\P^{\ell,b}( \ske(\BMM)=\q,w(\BMM) \in dw) }{d\text{Leb}^{\sharp E(\q)}}(\hat w,b(0)).
\end{align}
Putting together \eqref{e.1_for_3} and \eqref{e.2_for_3}, we have that \begin{align*}
    p_3^{\ell,b}(t) &= p_3^{\ell,b}(t_1)p_3^{\ell,b}(t_2) \int_{\R}\exp\left(-\frac{(b(0)-y)^2}{2t_1}-\frac{(y-b(0))^2}{2t_2}\right)dy\\
    &= p_3^{\ell,b}(t_1)p_3^{\ell,b}(t_2)\sqrt{\frac{2\pi t_1t_2}{t}},
\end{align*}
where, in the last equality, we used  Lemma \ref{lem-decomp}. Then, $p_3^{\ell,b}(t)$ satisfies
\begin{align*}
   \sqrt{2\pi t} p_3^{\ell,b}(t) &= \left(\sqrt{2\pi t_1}p_3^{\ell,b}(t_1)\right) \left(\sqrt{2\pi t_2}p_3^{\ell,b}(t_2)\right).
\end{align*}
This fact together with its continuity implies that there is $\lambda^{\ell,b}>0$ such that
\begin{align*}
    p_3^{\ell,b}(t) &=\frac{1}{\sqrt{2\pi t}} e^{-\lambda^{\ell,b}t}.
\end{align*}

For $p_1^{\ell,b}$ we do a similar computation to the one of $p_3^{\ell,b}$ to see that a.e. on $t, t_1, t_2$
\begin{align*}
    p_1^{\ell,b}(b(0),t)&= p_{3}^{\ell,b}(t_1)p_1^{\ell,b}(b(0),t_2)\int_{\R}\exp\left(-\frac{(b(0) - y)^2}{2t_1} - \frac{(b(0) - y)^2}{2t_2}\right)d\sigma\\
    &= p_{3}^{\ell,b}(t_1)p_1^{\ell,b}(b(0),t_2)\sqrt{\frac{2\pi t_1t_2}{t}},
\end{align*}
where, in the last equality, we used Lemma \ref{lem-decomp}. Then, we take $t_2\to 0$, and  obtain that 
\begin{align*}
r^{\ell,b}:=\displaystyle\lim_{t_2\to 0}\sqrt{2\pi t_2}p_1^{\ell,b}(b(0),t_2) \in \R^+.
\end{align*}
From here \eqref{e.formula_p1} follows.

For the second type of peeling step, we can do as before to decompose it as follows
\begin{align*}
    p_2^{\ell,b,\ell_1,\ell_2}(t) &= p_3^{\ell,b}(t_1)p_2^{\ell,b,\ell_1,\ell_2}(t_2)\exp\left(\frac{(b(0)-b(k))^{2} }{2 t_2 } \right)\int_{\R}\exp\left(-\frac{(b(0)-y)^2}{2t_1}-\frac{(y - b(k))^2}{2t_2}\right)dy\\
    &= p_3^{\ell,b}(t_1)p_2^{\ell,b,\ell_1,\ell_2}(t_2)\sqrt{\frac{2\pi t_1t_2}{t}}\exp\left(\frac{(b(0)-b(k))^{2} }{2 t_2 }-\frac{(b(0)-b(k))^2}{2t}\right).
\end{align*}
Taking again  $t_2\to 0$, we obtain that
\begin{align*}
s^{\ell,b,\ell_1,\ell_2} := \displaystyle\lim_{t_2\to 0}\sqrt{2\pi t_2} p_2^{\ell,b,\ell_1,\ell_2}(t_2)\exp\left(\frac{(b(0)-b(k))^2}{2t_2}\right) \in \R^+,
\end{align*}
from here \eqref{e.formula_p2} follows.
\end{proof}

The last lemma we need is the continuity of $p_1$.
\begin{lemma}\label{l.continuity_peleing_1_metric}
In the context of Lemma \ref{l.absolute_continuity_metric}, the function $p_1^{\ell,b}(w,x)$ is continuous in $(w,x) \in \R^+\times \R$.
\end{lemma}
\begin{proof}
    Take $\ell$, $b$ and $\q=[\vroot,\vertex]\glue\r$ with $\P^{\ell+1,\hat b}(\r)>0$ as in the proof of Lemma \ref{l.continuity_3_1}. For a decoration $\sigma$ in $\r$, we define $\sigma_x$ the decoration in $\q$ that extends $\sigma$ and takes value $x$ on the vertex $\vroot\vertex$. Now, we see that a.e. in $w$ and $(w_e)_{e\in E(\r)}$
    \begin{align}\label{e.1_for_1}
    &\frac{d\P^{\ell,b}\left(\ske(\BMM)=\q, w_{\vroot\vertex}\in dw, \MDec_{\vertex} \in dx, \displaystyle\min_{e\in \BMM} w_{e}\geq\epsilon\ \forall e\in E(\r), \max |\MDec_u|\leq R \ \forall u \in V(\r)\right)}{d\text{Leb}\times \mu (dw,dx)}(w,x)\\
    \nonumber\\
    &=p_1^{\ell,b}(w,x)\P^{\ell+1,\hat b_x}\left(\ske(\BMM)=\r,  \min_{e\in \BMM} w_{e}\geq\epsilon \ \forall e\in E(\r), \max |\MDec_u|\leq R \ \forall u \in V(\r)\right).
    \end{align}
    Furthermore, we can use the fact that $\BMM$ is a spin-decorated map to see that \eqref{e.1_for_1} is equal to
    \begin{align*}
    & \int p^{\ell,b}(\tilde q_{w,(w_e)},\sigma_x) \1_{\{w_e\geq \epsilon, \ |\sigma_u)|\leq R\} } \prod_{e\in E(\r)} dw_e \prod_{u\in V(\r)} \mu(d\sigma_u)\\
    = &\int p^{\ell,b}(\tilde q_{w,(w_e)},\sigma_{b_0}) \exp\left(\sum_{u\sim v}\frac{(b(0)-\sigma_u)^2-(x-\sigma_u)^2}{2w(uv)} \right ) \1_{\{w_e\geq \epsilon, \ |\sigma_u)|\leq R\} }\prod_{e\in E(\r)}dw_e \prod_{u\in V(\r)}\mu(d\sigma_u)\\
    = & \ p_1^{\ell,b}(w,b(0))\int p^{\ell,\hat b}(\tilde r_{(w_e)},\sigma)  \exp\left(\sum_{u\sim v}\frac{(b(0)-\sigma_u)^2-(x-\sigma_u)^2}{2w(uv)} \right ) \1_{\{w_e\geq \epsilon, \ |\sigma_u)|\leq R\} }\prod_{e\in E(\r)}dw_e \prod_{u\in V(\r)} \mu(d\sigma_u).
    \end{align*}
    We conclude putting both equalities together and using that the exponential function appearing inside the integral is continuous in $x$ and bounded and the integral is not $0$ as long as $R$ is big enough. Note that $p_1^{\ell,b}(w,b(0))$ is continuous in $w$ and the integral does not depend on $w$.
\end{proof}

Now, we begin with the proof of the main result of this section. Part of the proof is similar to that of Theorem 2.19 in \cite{araya2025characterisation}. However, there are both technical difficulties as in the metric setting, the peeling process is continuous, and non technical difficulties, from the fact that we need to understand the role that the boundary values are playing.
\begin{proof}[Proof of Theorem \ref{t.characterization_metric}]
    Thanks to Proposition \ref{prop:boltzman_invariance}, we only need to prove that $(2) \Rightarrow (1)$ and $(3)\Rightarrow(2)$. 
    
    First, fix $w_1,w_2>0$, $\ell\in\Z_+$ and $b\in\Bc$. Let us first prove that $\lambda^{\ell,b}$ does not depend on $b$. By comparing energies, it is possible to see that
\begin{align*}
p_3^{\ell,b}(w_1,x)\,p_2^{\ell,b_x,\ell-1,0}(w_2)=p_3^{\ell,b}(w_1,b(0))\,p_2^{\ell,b,\ell-1,0}(w_2)\exp\big(-\tfrac{(x-b(0))^2}{2w_1}-\tfrac{(x-b(1))^2}{2w_2}+\tfrac{(b(0)-b(1))^2}{2w_2}\big).
\end{align*}
Thus, after using \eqref{e.formula_p2} we see that $e^{w_2(\lambda^{\ell,b_x}-\lambda^{\ell,b})}$
does not depend\footnote{This uses the fact that $s^{\ell,b,\ell-1,0}>0$. This property follows from the fact that the probability of obtaining a given tree is strictly positive, as $\BMM$ is a spin-decorated map together with the fact that it is Gibbs distributed.} on $w_2$. This implies that $\lambda^{\ell,b}$ is only a function of $\ell$.

Now, we prove that $\lambda^{\ell}$ is constant, to do this, assume $\ell\geq 2$. From studying a map where both the boundary edges $0$ and $1$, and $2$ and $3$ are identified, we see that for any $w_1,w_2\in \R$
\begin{align*}
    p_2^{\ell,b,\ell-1,0}(w_1)p_2^{\ell-1, \dot b,\ell-2,0}(w_2)= p_2^{\ell,b_s,\ell-1,0}(w_2)p_2^{\ell-1,\dot b_s,\ell-2,0}(w_1),
\end{align*}
where $b_s$ is the right shift on $b$ so that $2$ plays the role of the root edge and $\dot b$, resp. $\dot b_s$, is the boundary condition that is left after doing a peeling step. As a consequence of the equation, we have that $e^{ w_1(\lambda^{\ell}-\lambda^{\ell-1})}$ does not depend on $w_1$, which implies that $\lambda^{\ell}$ does not depend on $\ell$. This means that for all $t>0$, $p_3^{\ell,b}(t)=p_3(t)$.

    \begin{itemize}
        \item $(2)\Rightarrow (1):$ Take $\Mq\in\MQm_H$ and $\sigma$ a decoration over $V_i(\Mq)$. We define $q_{w_1,w_2}(\ell, b, \Mq, \sigma)$ as follows
        \begin{align*}
        q_{w_1,w_2}(\ell, b, \Mq, \sigma) &:= \frac{p^{\ell,b}\left(\Mq(w_1,w_2),\overline{\sigma}\right)}{p^{\ell,b}\left(\Mq,\sigma\right)},
        \end{align*}
    where, $\Mq(w_1,w_2)$ denotes a copy of $\Mq$ where we add a new internal vertex, two new edges of length $w_1$ and $w_2$ and $\overline{\sigma}$ is an extension of the decoration, which assigns the decoration of the phantom vertex for the new internal vertex as presented in figure \ref{fig:trans}.
    
    \begin{figure}[h!]
        \centering
        \includegraphics[scale = 0.85]{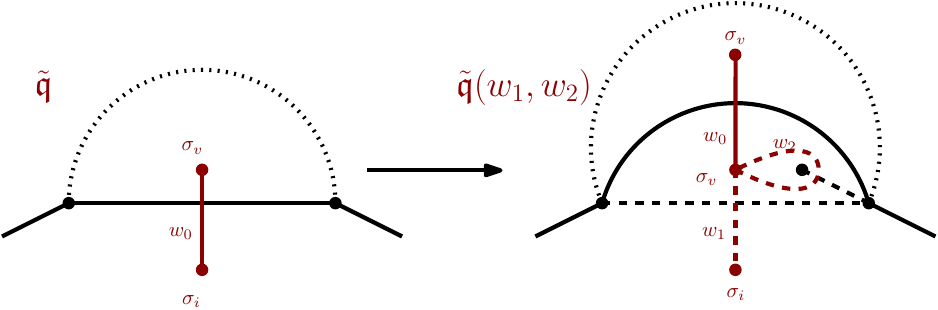}
        \caption{An illustration of the construction used in the proof of Theorem 2.19 in \cite{araya2025characterisation}. The left panel shows an interior vertex with spin $\sigma_i$ connected to a phantom vertex with spin $\sigma_v$. In the right panel, we introduce a new interior vertex and two additional edges. The new vertex is assigned the same spin as the phantom vertex, so that the contribution to the Hamiltonian associated with the original connection between the phantom and interior vertices remains unchanged.}
        \label{fig:trans}
    \end{figure}
    
We are going to prove that $q_{w_1,w_2}(\cdot,\cdot, \cdot,\cdot)$ is constant.
We first prove for a fixed $w_1,w_2>0$, for any $\ell\in\N^{*}$ and $b\in\Bc$, $q_{w_1,w_2}(\ell, b, \cdot, \cdot)$ is constant. We first see that $q_{w_1,w_2}(\ell, b, \Mq, \sigma)$ does not depend on the structure of $\Mq$ and on the number of vertices of $\Mq$ by using the Markov property with a peeling step of type 1 and one of type 2
    \begin{align*}
        q_{w_1,w_2}(\ell, b, \Mq, \sigma) &= \frac{p^{\ell,b}\left(\Mq(w_1,w_2),\overline{\sigma}\right)}{p^{\ell,b}\left(\Mq,\sigma\right)} = \frac{p_1^{\ell,b}(w_1,b(0))p^{\ell+1,\hat{b}}\left(\widehat{\Mq}(w_2),\widehat{\sigma}\right)}{p^{\ell,b}\left(\Mq,\sigma\right)} = p_1^{\ell,b}(b(0),w_1)p_2^{\ell+1,\hat{b},0,\ell}(w_2).
    \end{align*}
Here $\widehat{\Mq}$ is the resulting metric map obtained after doing a peeling of type $1$ and $\widehat{\sigma}$ is the corresponding decoration.
This last term does not depend on $\Mq$ and $\sigma$. Notice that, thanks to Lemma \ref{l.peleeing_decom}, we have the following
\begin{align}\label{eq.q_decomp_metric}
    q_{w_1,w_2}(\ell,b) &= q(\ell,b)p_3(w_1)p_3(w_2).
\end{align}
where $q(\ell,b) := r^{\ell,b}s^{\ell+1,\hat{b},0,\ell}$.

Now, we have that for any $(\Mq_1,\sigma_1)$ and $(\Mq_2,\sigma_2)$ where $\Mq_1$ and $\Mq_2$ are comparables
\begin{align}\label{e.ratio_metric}
    \frac{p^{\ell,b}(\Mq_1,\sigma_1)}{p^{\ell,b}(\Mq_2,\sigma_2)} = \exp\left(\Ham^b(\Mq_2,\sigma_2)-\Ham^b(\Mq_1,\sigma_1)\right).
\end{align} 
This follows from \eqref{e.def_p_q_s_metric} and having Gibbs distribution on $\MQm^{\ell,f}$ (see \eqref{gibbs-fix-faces-dec-met}). We first apply this formula for trees. Note that if $\tilde \tr$ is a tree of size $\ell$, we can always peel it from its leaves and obtain that
\begin{align}\label{e.base_case_metric_quad}
    p^{\ell,b}(\tilde \tr)= e^{-\Ham^{b}(\tilde \tr)}\prod_{k=0}^{\ell-1} s^{\ell-k,b_k,\ell-k-1,0} p_3(w_k),
\end{align}
where $b_k$ is the right boundary condition for the edge peeled at time $b_k$. This together with \eqref{e.ratio_metric} implies that $\prod_{k=0}^{\ell-1} s^{\ell-k,b_k,\ell-k-1,0}$ does not depend on the skeleton of $\tilde \tr$.

This together with Lemma \ref{l.absolute_continuity_metric} applied to trees, implies that if $\tr_1$ and $\tr_2$ are two trees with the same lengths, then we obtain that $p^{\ell,b}$ follows a Boltzmann formula for trees; meaning that the probability in Definition \ref{def:BDMP} is only supported on trees.

Take $\Mq_1$ and $\Mq_2$ two 4-regular maps such that 
\begin{itemize}
    \item $|V_i(\Mq_1)| = |V_i(\Mq_2)|+1$;
    \item $|E(\Mq_1)| = |E(\Mq_2)| + 2$;
    \item and, there are $|E(\Mq_2)|$ edges of $E(\Mq_1)$ that have the same lengths as $E(\Mq_2)$. 
\end{itemize}
Then, using \eqref{e.ratio_metric}
\begin{align}
    \nonumber\frac{p^{\ell,b}(\Mq_1,\sigma_1)}{p^{\ell,b}(\Mq_2,\sigma_2)} &= \frac{p^{\ell,b}(\Mq_1,\sigma_1)}{p^{\ell,b}(\Mq_2(w_1,w_2),\overline{\sigma_2})}\frac{p^{\ell,b}(\Mq_2(w_1,w_2),\overline{\sigma_2})}{p^{\ell,b}(\Mq_2,\sigma_2)}\\  
    &\nonumber= \exp\left(\Ham^b(\Mq_2,\sigma_2)-\Ham^b(\Mq_1,\sigma_1)\right)q_{w_1,w_2}(\ell,b)\\
    \label{e.ratio_metric_2}&= q(\ell,b)\exp\left(\Ham^b(\Mq_2,\sigma_2)-\Ham^b(\Mq_1,\sigma_1)\right)p_3(w_1)p_3(w_2),
\end{align}
where, $w_1$ and $w_2$ are the sizes of the edges that are different between $E(\Mq_1)$ and $E(\Mq_2)$.

 This implies by induction that
\begin{align}\label{e.boltzmann_epsilon}
    p^{\ell,b}(\Mq,\sigma) = \frac{q(\ell,b)^{|V_i(\Mq)|}}{W^{\ell,b}}\exp\left(-\Ham^b(\Mq,\sigma)\right)\prod_{e \in E(\Mq)}p_3(w_{e}),
\end{align}
where $W^{\ell,b}<\infty$. It remains to prove that $q(\ell,b)$ does not depend on $\ell$ and $b$. Denote by $\sigma_x$ to the decoration over $\Mq(w_1,w_2)$ that is equal to $\sigma$ everywhere except that its value in the new vertex is $x\in \text{supp} \ \mu$. Suppose further that the new vertex's neighbours $v$ in $\Mq(w_1,w_2)$ satisfies $\sigma(v) = x$. Then
\begin{align*}
    p^{\ell,b}(\Mq(w_1,w_2),\overline{\sigma}) = p^{\ell,b}(\Mq(w_1,w_2),\sigma_x).
\end{align*}
Define, now, $b_x$ as the boundary condition that takes values $b_x(0) = x$ and everywhere else coincides with $b$. In consequence, and thanks to \eqref{e.boltzmann_epsilon}, we obtain that
\begin{align*}
    q_{w_1,w_2}(\ell,b) &= \frac{p^{\ell,b}(\Mq(w_1,w_2),\sigma_x)}{p^{\ell,b}(\Mq,\sigma)} \\
    &= p_1^{\ell,b}(b_0,w_1)p_2^{\ell+1,\hat{b}_x,0,\ell}(w_2)\frac{p^{\ell,b_x}(\Mq,\sigma) }{ p^{\ell,b}(\Mq,\sigma)}\\
    &=p_1^{\ell,b}(b_0,w_1)p_2^{\ell+1,\hat{b}_x,0,\ell}(w_2) \left(\frac{q(\ell,b_x)}{q(\ell,b)}\right)^{|V_i(\q)|}\frac{W_{q(\ell,b)}^{\ell,b}}{W_{q(\ell,b_x)}^{\ell,b_x}}e^{\frac{1}{2} \frac{(x-b(0))^2}{w_1}}.
\end{align*}
Since the left-hand side does not depend on the number of vertices we conclude that $q(\ell,b) = q(\ell,b_x)$. This means that we can modify just one value of $b$ and keep the value of $q$. By rerooting invariance we conclude that we can change all of them.

Finally, we are going to prove that $q(\ell)$ does not depend on $\ell$. For this we compute
\begin{align*}
    q_{w_1,w_2}(\ell,b) &= \frac{p_1^{\ell,b}(b_0,w_1)p^{\ell+1,\hat{b}}\left(\widehat{\Mq}(w_2),\widehat{\sigma}\right)}{p^{\ell,b}\left(\Mq,\sigma\right)}\\
    &= p_1^{\ell,b}(b_0,w_1)\left(\frac{q(\ell+1)}{q(\ell)}\right)^{|V_i(\Mq)|}\frac{W_{q(\ell,b)}^{\ell,b}}{W_{q(\ell+1,\hat{b})}^{\ell+1,\hat{b}}}p_{3}(w_2).
\end{align*}
Again, as the left-hand side does not depend on the number of vertices we conclude that $q(\ell+1) = q(\ell)$. And, in consequence, $q_{w_1,w_2}$ does not depend on $\ell$.

\item $(3)\Rightarrow (2):$ We need to show that for any $\Mq_1,\Mq_2$ as in \eqref{gibbs-fix-faces-dec-met} we have the following
\begin{align}\label{e.metric_Gibbs}
    \frac{p^{\ell,b}\left(\Mq_1\right)}{p^{\ell,b}\left(\Mq_2\right)} &= \frac{Z^{b}\left(\Mq_1\right)}{Z^{b}\left(\Mq_2\right)}.
\end{align}
In fact, we prove something stronger: if $\Mq_1$ and $\Mq_2$ have the same number of interior vertices
\begin{align}\label{e.change_length_metric}
    \frac{p^{\ell,b}\left (\Mq_1\right)}{Z^{b}\left(\Mq_1\right)\prod_{e \in E( \Mq_1)} p_3(w^1_e) }= \frac{p^{\ell,b}\left (\Mq_2\right)}{Z^{b}\left(\Mq_2\right)\prod_{e \in E( \Mq_2)} p_3(w^2_e) },
\end{align}
where $w_e^i$ is the length of the edge $e$ in $\Mq_i$. We prove this by induction in the number of vertices. We temporarily assume the following claim, which is the base case of the induction. Its proof will be made in the end.
\begin{claim}\label{c.inductive_step_metric}
 If $\Mq_1$ and $\Mq_2$ are dual metric trees with $\ell$ edges, then \eqref{e.change_length_metric} holds.
\end{claim}

For the inductive step, take $n\in\N$ and assume that \eqref{e.metric_Gibbs} and \eqref{e.change_length_metric} hold for metric maps with no more than $n$ interior vertices. Take $\Mq_1$ and $\Mq_2$ with $n+1$ interior vertices. Thanks to re-rooting invariance, let us first assume that the root of $\Mq_1$ and $\Mq_2$ are the boundary of an internal vertex that we call $\mathfrak{v}$. Additionally, define $\Mq_{i}^{*}$ as the result of applying a peeling step on the root to $\Mq_i$ and for $x\in \text{supp}(\mu)$
\begin{align*}
    p^{\ell,b}(\Mq_i,x) := \frac{d\P(\BMM\in d\Mq_i,\Dec(\mathfrak{v})\in dx)}{\mu(dx)}
\end{align*}
is the Radon-Nikodym derivative of the decoration at the vertex $\mathfrak{v}_i$ with respect to the measure $\mu$ on
the 0-measure event that the map is exactly $\Mq_i$. To show that \eqref{e.metric_Gibbs} is true for $\Mq_1$ and $\Mq_2$, we compute
\begin{align*}
    p^{\ell,b}(\Mq_i,b(0)) &= \P^{\ell,b}(\BMM\in d\Mq_i)\frac{1}{Z^{b}(\Mq_i)}\int \exp\left(-\Ham^{b^{*}}(\Mq_i^{*},\sigma)\right)\mu^{n-1}(d\sigma) = \P^{\ell,b}(\BMM\in d\Mq_i)\frac{Z^{b^{*}}(\Mq_i^{*})}{Z^{b}(\Mq_i)},
\end{align*}
where $b^{*}$ is the boundary condition that appears when doing a peeling step on the vertex $\vertex$ and discovering that it has value $b(0)$. Then, we can use the Markov property to peel the vertex $\vertex$
\begin{align*}
    \frac{p^{\ell,b}(\Mq_1,b(0))}{Z^{b^{*}}(\Mq_1^{*})}& = p_1^{\ell,b}(b(0),w_0^1)\left(\prod_{e\in E(\Mq_1^*)}p_3(w^1_e)\right) \frac{p^{\ell+1,b^{*}}(\Mq_1^{*})}{Z^{b^{*}}(\Mq_1^{*})\prod_{e\in E(\Mq_1^*)}p_3(w^1_e)
    } \\
    &= p_1^{\ell,b}(b(0),w_0^1)\left(\prod_{e\in E(\Mq_1^*)}p_3(w^1_e)\right)\frac{p^{\ell+1,b^{*}}(\Mq_2^{*})}{Z^{b^{*}}(\Mq_2^{*})\prod_{e\in E(\Mq_2^*)}p_3(w^2_e)}\\ &=\frac{p_1^{\ell,b}(b(0),w_0^1) \prod_{e\in E(\Mq_1^*)}p_3(w_e^1)}{p_1^{\ell,b}(b(0),w_0^2)\prod_{e\in E(\Mq_2^*)}p_3(w_e^2) } \cdot \frac{p^{\ell,b}(\Mq_2)}{Z^{b^*}(\Mq_2)},
\end{align*}
where we use the induction hypothesis for $\Mq_1^{*}$ and $\Mq_2^{*}$. Here, $w_0^i$, denotes the length of the root edge of $\Mq_i$, while $w_e^i$ denotes the length of the edge $e$. Note that, thanks to Lemma \ref{l.absolute_continuity_metric},
\begin{align*}
    \frac{p_1^{\ell,b}(b(0),w_0^1) }{p_1^{\ell,b}(b(0),w_0^2) } = \frac{p_3(w_0^1) }{p_3(w_0^2) }.
    \end{align*}

Now, in the case that the edge of the root of $\tilde \q_2$ does not connect with an internal vertex but has a face connected to an edge $\mathfrak e$, we proceed as in Theorem \ref{t.uniqueness_decorated} and take $\tilde \q_3$ that has an internal vertex connected to the root and to $\mathfrak e$. Then
\begin{align*}
\frac{p^{\ell,b}\left (\Mq_1\right)}{\prod_{e \in E(\tilde \q_1)} p_3(e) Z^{b}\left(\Mq_1\right)}= \frac{p^{\ell,b}\left (\Mq_3\right)}{\prod_{e \in E(\tilde \q_3)} p_3(e) Z^{b}\left(\Mq_3\right)}= \frac{p^{\ell,b}\left (\Mq_2\right)}{\prod_{e \in E(\tilde \q_2)} p_3(e) Z^{b}\left(\Mq_2\right)}.
\end{align*}
    \end{itemize}

We are only missing the claim.
\begin{proof}[Proof of Claim \ref{c.inductive_step_metric}] Note first that by using the right peeling steps and Lemma \ref{l.peleeing_decom}, there exists a function $f^{\ell,b}$ that takes skeletons of trees such that for any dual metric tree
\begin{align}\label{base_case}
    p^{\ell,b}(\tilde \tr)= f^{\ell,b}(\ske(\tr)) Z^{\ell,b}(\tr)\prod_{ e \in E(\tr)} p_3(w_{e}) .
\end{align}
We need to show that $ f^{\ell,b}(\cdot)$ is constant. We proceed by induction on the number of edges $\ell$ of the tree. Again, the induction step is easier. Take $\ell\geq 2$ and assume that $f^{\ell,b}(\ske(\tr))$ is constant for any $\ell$, $b$ with length smaller than or equal $\ell$. Then, if $\tilde \tr_1$ and $\tilde \tr_2$ are two trees such that the root edge is incident to a leaf, we have that
\begin{align*}
    p^{\ell,b}(\tilde \tr_i)=p_2^{\ell,b,\ell-1,0}(w_0^1) p^{\ell-1,b^*}(\tilde \tr_i^*) = p_3^{\ell,b}(w_0^1)s^{\ell,b,\ell-1,0}p^{\ell-1,b*}(\tilde \tr_i^*).
\end{align*}
This, together with adding a proper $\tilde \tr_3$ when $\tilde \tr_1$ and $\tilde \tr_2$ do not share a leaf allows us to conclude the inductive step.

For the base case denote by $e = (0,j)$ the root edge of length $w_e = w$ of the tree $\tilde \tr$. If we explore a portion $u$ of the edge, thanks to the Markov property, we obtain the following
\begin{align*}
    p^{\ell,b}(\tilde \tr) = p_3^{\ell,b}(u,x)p^{\ell,b_x}(\tilde\tr_u)\sqrt{\frac{2\pi u(w-u)}{w}} \exp\big(\tfrac{(x-b(0))^2}{2u}+\tfrac{(x-b(j))^2}{2(w-u)}-\tfrac{(b(0)-b(j))^2}{2w}\big),
\end{align*}
where $\tilde\tr_u$ is the tree where the root edge has length $w_e-u$ and $b_x$ denotes the boundary condition $b$ where we replace $b(0)$ by $x$. Note that $\ske(\tilde\tr_u) = \ske(\tilde\tr)$. Then, by \ref{base_case} we obtain the following
\begin{align*}
    \frac{f^{\ell,b}(\ske(\tilde\tr))}{f^{\ell,b_x}(\ske(\tilde\tr))} = p_3^{\ell,b}(u,x)e^{\lambda u}\sqrt{2\pi u}\exp\left(\frac{(b(0)-x)^2}{2u}\right),
\end{align*}
where we note that the right hand side of the previous equation does not depends on the skeleton of $\tr$. Now take $\tilde\tr_1$ and $\tilde\tr_2$ two dual metric trees with $\ell(\tilde\tr_1) = \ell(\tilde\tr_2)$. Then
\begin{align*}
    \frac{f^{\ell,b}(\ske(\tilde\tr_1))}{f^{\ell,b}(\ske(\tilde\tr_2))} = \frac{f^{\ell,b_x}(\ske(\tilde\tr_1))}{f^{\ell,b_x}(\ske(\tilde\tr_2))},
\end{align*}
where we conclude that the quotient does not depend on  the boundary condition $b$. Finally take a constant boundary condition to see that, by invariance under rerooting, $f^{\ell,b}(\tilde\tr_1) = f^{\ell,b}(\tilde\tr_2)$ where we conclude.

\end{proof}
\end{proof}

\section{Decorated maps: Markov property and associated measures}\label{sec:5}
In this section, we are interested in the study of spin-decorated maps. That is to say, planar maps that has a decoration on top of it. We modify the model of the previous section by forgetting the lengths of the edges and focusing in the skeleton of metric maps. The decoration now is a function $\sigma$ from the vertices of $\BM$ to $\R^d$, such that conditionally on $\BM$ the law of $\sigma$ is that of a spin system on $\BM$. In this case, the measures we work with are simpler and almost all the results of the metric decorated case are analogous to this model.

We start by defining spin-decorated maps. Then, we introduce decorated Boltzmann 4-regular maps with boundaries and show that they satisfy the weak and strong Markov property for simpler framework than before. Then, once again, we characterise all `reasonable' spin-decorated 4-regular maps satisfying the Markov property. Almost all of the proofs for this case are analogous and simpler than the previous section.

\subsection{Decorated Maps} 
Consider a pair $(\q,\sigma)$, where the first coordinate is a 4-regular map $\q\in\Qm_H$ and the second coordinate is a function $\sigma:V(\q)\to \R^n$. The notion of boundary condition is analogous to this model. Notice that if a 4-regular map has half-perimeter $\ell$ then $b\in (\R^n)^{2\ell}$. The measure we work in this section satisfies that $\supp(\mu)\subseteq \R^n$  and we denote by $\Bc$ the space of all boundary condition for $\Qm$ that take values in $\supp(\mu)$. 

Since in these model there are no continuous edges, we only need to fix the inverse temperature $\beta>0$.

\begin{defn}
We say that a random pair $(\BM,\Dec)$ is a SDM of half-perimeter $\ell$ and boundary condition $b$ if its probability measure satisfies that for $\Dec$ with boundary condition $b$
\begin{align} \label{def_cond}
    \P^{\ell,b}(\Dec\in d\sigma|\BM=\q)& =  \frac{1}{Z^b(\q)}\exp\left(-\Ham^b(\q,\sigma)\right)\prod_{v\in V_i(\q)}\mu(d\sigma_v),
\end{align}
where $\Ham^{b}(\q,\sigma):= \frac{\beta}{2}\displaystyle\sum_{i\sim j}\|\sigma_i-\sigma_j\|^2$ and where
\begin{align}\label{Zq}
    Z^b(\q) = Z_{\beta, \ell,\mu}^b(\q) := \int_{(\R^n)^{|V_i(\q)|}}\exp\left(-\Ham^b(\q,\sigma)\right)\prod_{v\in V_i(\q)}\mu(d\sigma_v).
\end{align}
\end{defn}
A special case of SDM are Boltzmann decorated maps. 
\begin{defn}
    We say that $(\BM,\Dec)$ is a \textit{$q$-Boltzmann decorated map} with half-perimeter $\ell$ and boundary condition $b$ if its density satisfies the following
\begin{align}
    \P^{\ell,b}\left(\BM=\q,\Dec\in d\sigma\right)\propto q^{|V_i(\q)|}\exp\left(-\Ham^{\beta,b}(\q,\sigma)\right)\prod_{v\in V_i(\q)}\mu(d\sigma_v),
\end{align}
where the value of the normalisation constant $W^{\ell,b}_{q,\beta}$ is given by
\begin{align}
    W^{\ell,b}_{q,\beta} &= \sum_{\q\in\Qm^{\ell}} q^{|V_i(\q)|}\int_{(\R^{n})^{|V_i(\q)|}}\exp\left(-\Ham^{\beta,b}(\q,\sigma)\right)\prod_{v\in V_i(\q)}\mu(d\sigma_v).
\end{align}
\end{defn}

\begin{rem}
For most of the measures we are interested in there always exists a $q>0$ s.t. $W_{q}^{\ell,b}<\infty$. This is clearly the case when the measure $\mu$ is finite as one can compare with the case without decoration. When the measure is Lebesgue in $\R^n$, we can check by noting that 

\begin{align*}
Z^b(\q)=c(b,\q)(2\pi/\beta)^{N/2}\text{det}(\Delta_D)^{-1/2},
\end{align*} 

where $\Delta$ is the Laplacian in $\q$ with boundary conditions $0$ (rooted at the boundary) and $N = n\cdot |V_i(\q)|$. This Laplacian is related by the Matrix-tree theorem to the number of spanning trees in $\q$ which grows at most exponentially with the number of interior vertices of $\q$. For the case where the decoration takes values on $\Z$ and $\mu$ is the counting measure, we can manually bound the $Z^b_\Z(\q)\leq c^{|V(q)|}(Z^b_\R(\q) + 1)$. 
\end{rem}

Now, we state the analogous properties of being invariant under rerooting (\ref{Invariance-root-dec-met}) and have Gibbs distribution on $\Qm^{\ell,f}$ (\ref{gibbs-fix-faces-dec-met}) for spin-decorated maps.
\begin{enumerate}
    \item\label{Invariance-root-dec} \textbf{(Invariance under rerooting)}
    The SDM $(\BM,\Dec)\sim \P^{\ell,b}$ with half-perimeter $\ell$ and boundary condition $b$ is \textit{invariant under rerooting} if after rerooting it, its law is that of $\P^{\ell,b_{s}}$ where $b_s$ is the proper shift of $b$.
    
    \item\label{gibbs-fix-faces-dec} \textbf{(Gibbs distribution on $\Qm^{\ell,f}$)} The SDM $(\BM,\Dec)\sim\P^{\ell,b}$ has \textit{Gibbs distribution on $\Qm^{\ell,f}$} if for any $\q_1,\q_2 \in \Qm^{\ell,f}$
    \begin{align*}
        \frac{\P^{\ell,b}\left(\BM= \q_1\right)}{\P^{\ell,b}\left(\BM= \q_2\right)} = \frac{Z^b(\q_1)}{Z^b(\q_2)}.
    \end{align*}
\end{enumerate}

Once again, an example of probability laws that satisfy these two previous properties are the Boltzmann decorated maps. 
\begin{prop}
    Boltzmann decorated maps are invariant under rerooting (\ref{Invariance-root-dec}) and have Gibbs distribution over $\Qm^{\ell,f}$ (\ref{gibbs-fix-faces-dec}). 
\end{prop}
\begin{proof}
The proof of this result is analogous to the proof made in section \ref{sec:3}.    
\end{proof}

\subsection{Markov property: stopping maps} 
In this subsection we state the Markov property for the Boltzmann maps in its two versions: the weak one and the strong one. Also, we show that the only random submaps that satisfy the Markov property for Boltzmann maps are the stopping maps.

First, we state the weak Markov property associated with the Boltzmann decorated maps.
\begin{thm}\label{t.decorated_weak_Markov}
    Let $(\BM,\Dec)$ be a Boltzmann decorated map and $\q\in\Qm_{H}$. Then, conditionally on $\{\q\Lalg\BM\}$ and $\Dec|_{\q}$, $(\BM,\Dec)$ can be decomposed as follows
    \begin{align} \label{e.decomposition}
        \BM = \q\glue (\BM_h^{\q})_{h\in H(\q)},\  \text { and } \ 
        \Dec = \Dec|_{\q} + \Dec^{\q},
    \end{align}
    where, $(\BM_{h}^{\q},\Dec_{h}^{\q})$ is a collection of independent Boltzmann decorated maps with boundary $\ell(h)$ and boundary condition $\Dec|_{h}$ for every $h\in H(\q)$.
\end{thm}
The boundary conditions at hole vertices are defined as follows. From the perspective of each hole vertex, its adjacent vertices are treated as phantom vertices. Each phantom vertex retains the spin of the corresponding original vertex, as determined by the information of $\Dec_\q$, and is connected only to the corresponding interior vertex of the hole.

\begin{proof}[Proof of Theorem \ref{t.decorated_weak_Markov}]
The proof of this theorem is exactly analogous and simpler than the proof made in Theorem \ref{t.metric_weak_markov}. Thus we refer the reader to that proof.
\end{proof}

Since the model in this section is simpler, the definitions that work for this case are simplified. The definition of filtration is the same as presented in \ref{def:filtration}. The definition that needs to be changed is the notion of Adaptability and Independent increments.

\begin{defn}\label{def_dec_Boltzmann}
Let $\FF = (\F_{\q})_{\q\in\Qm_H}$ be a filtration. We say that a Boltzmann decorated map $(\BM,\Dec)$ is an \textit{$\FF-$Boltzmann decorated map} if it satisfies the following properties
\begin{itemize}[label={--}]
    \item\textit{Adaptability:} the event $\{\q\Lalg \BM\}$ and the function $\Dec|_{\q}\1_{\{\q\Lalg \BM\}}$ are $\F_{\q}$-measurable.
    \item\textit{Independent increments:} conditionally on $\F_{\q}$, on the event $\{\q\Lalg\BM\}$, the law of $(\BM_h^{\q},\Dec^{\q}_h)_{h\in H(\q)}$ is that of a collection of independent $q$-Boltzmann decorated maps with boundary $h$ and boundary condition $\Dec|_{h}$ for every $h\in H(\q)$.
\end{itemize}
\end{defn}

\begin{rem}
As before, an example of a filtration for this model is the \textit{natural filtration}, defined as follows
\begin{align}
    \F_{\q} := \overline{\bigvee_{\p\Lalg \q} \sigma\left( \left\{\p \Lalg \BM \right\},\Dec|_{\p}\1_{\{\p\Lalg \BM\}} \right)}^{\P}.
\end{align} 
Of course the weak Markov property implies that $(\BM,\Dec)$ is an $\FF$-Boltzmann decorated map for its natural filtration.
\end{rem}
In this context, we keep the definition of stopping map as well as the definition of the $\sigma$-algebra associated with these random maps given in Definition \ref{def_stop}.

The following result is the strong Markov property for the Boltzmann decorated maps.
\begin{thm}\label{t.strong_Markov_decorated}
   Let $\FF = (\F_{\q})_{\q\in\Qm_H}$ be a filtration. Let $(\BM,\Dec)$ be an $\FF$-Boltzmann decorated map and $\SM\in\Qm_H$ an $\FF-$stopping map for $(\BM,\Dec)$. Then, $\BM$ and $\Dec$ can be decomposed as follows,
    \begin{align*}
        \BM = \SM\glue (\BM^{\SM}_h)_{h\in H(\SM)},\  \text { and } \ \Dec = \Dec|_{\SM} + \Dec^{\SM},
    \end{align*}
    where, conditional on $\F_{\SM}$ and $\Dec|_{\SM}$, $(\BM_h^{\SM},\Dec^{\SM}_h)_{h\in H(\SM)}$ is a collection of independent Boltzmann decorated maps with boundary $\ell(h)$ and boundary condition $\Dec|_h$ for every $h\in H(\SM)$.
\end{thm}
\begin{proof}
The proof of this theorem is simpler than the one of Theorem \ref{t.metric_strong_Markov}. Actually, since now the structure has countable number of values, the proof made in Lemma \ref{lem-countable} is enough for this case.
\end{proof}

A natural question that arises is which random maps allow a Markovian decomposition. Once again, the answer is the notion of local maps.

\begin{defn}
We say that the triplet $\left(\BM,\Dec,\SM\right)$ is a \textit{local map} if,
\begin{itemize}
    \item $(\BM,\Dec)$ is a (decorated) Boltzmann map;
    \item $\P$-almost surely $\SM\Lalg\BM$;
    \item and, conditionally on $\SM$ and $\Dec_{\SM}$, $(\BM_h^{\SM},\Dec^{\SM}_h)_{h\in H(\SM)}$ is a collection of independent Boltzmann decorated maps with perimeter $|h|$ and boundary condition $\Dec|_h$ for every $h\in H(\SM)$. 
\end{itemize}
\end{defn}
The next result is the analogous result of Theorem \ref{t.local_implies_stopping}.
\begin{thm}\label{t.local_implies_stopping_dec}
    Let $\left(\BM,\Dec,\SM\right)$ be a local map. Then, there exists a filtration $\FF = (\FF_{\q})_{\q\in\Qm_H}$ such that $(\BM,\Dec)$ is an $\FF-$Boltzmann decorated map and $\SM$ is an $\FF-$stopping map for $(\BM,\Dec)$.
\end{thm}
Once again, the proof of this result is the same as that presented in the previous section. Therefore we only specify the filtration that satisfies the result. For $\q\in\Qm_H$ we define the following $\sigma$-algebra
\begin{align}\label{e.def_filt}
        \F_{\q} = \overline{\bigvee_{\p\Lalg\q} \sigma\left(\{\p\Lalg\BM\},\{\SM\Lalg \p\}, \Dec|_{\p}\1_{\p\Lalg\BM}, \SM^{\wedge \p}\1_{\p\Lalg\BM} \right)}^{\P}.
\end{align}
Clearly $\SM$ is an $\FF$-stopping map and $(\BM,\Dec)$ satisfies the property of Adaptability. The only thing that is left to prove is the property of Independent increments but this comes directly by Lemmas \ref{lem.local_stopping_1} and \ref{lem.local_stopping_2}.

\subsection{Characterisation of Markovian laws on decorated maps.}

As before, we present a characterisation of all the SDM satisfying the Markov property. In order to do this, we need to state the Markov property for a general collection of measures $(\P^{\ell,b})_{(\ell,b)\in\N^*\times\Bc}$ 
\begin{defn}
We say that $(\P^{\ell,b})_{(\ell,b)\in\N^*\times\Bc}$ satisfies the \textit{Markov property}, if for any deterministic 4-regular map with hole vertices $\q$, we can describe the conditional law $\P^{\ell,b}(\cdot|\{\q\Lalg\BM\})$ as 
    \begin{align}
        \BM = \q \glue (\BM_h^{\q})_{h\in H(\q)},\  \text { and } \ \Dec = \Dec|_{\q} + \Dec^{\q}.
    \end{align}
    Here, conditional to $\Dec|_{\q}$ ,$(\BM^{\q}_h,\Dec^{\q}_h)_{h\in H(\q)}$ is a collection of independent decorated maps with law $\P^{\ell(h),\Dec|_{h}}$ for every $h\in H(\q)$.
\end{defn}
The peeling exploration introduced in the previous section extends naturally as follows: each time that a new vertex is discovered, it also reveals the value of its spin.

Now, we are ready to state the main result of this section that characterises all the decorated 4-regular maps having the Markov property.
\begin{thm}\label{t.uniqueness_decorated}
Take a collection of measures $(\P^{\ell,b})_{(\ell,b)\in \N^*\times \Bc}$ on spin-decorated 4-regular maps with a half-perimeter $\ell$ and boundary condition equal to $b:\cro{0,2\ell-1}\mapsto \supp \mu$. Furthermore, assume that for any $\ell$ and $b$ we have that 
\begin{align*}
    \P^{\ell,b}(\BM \text{ has at least one interior vertex})>0,
\end{align*}
and that the measure $\P^{\ell,b}$ is continuous on $b\in(\supp\mu)^{2\ell}$ for the weak topology of measures. Then, the following are equivalent.
\begin{enumerate}
    \item There is a $q>0$ such that for any $\ell$ and $b$, $\P^{\ell,b}$ is a $q$-Boltzmann decorated map, 
    \item $(\P^{\ell,b})_{(\ell,b)\in \N^*\times \Bc}$ has the Markov property, and for each $\ell\in \N$ and $b\in \Bc$, the measure $\P^{\ell,b}$ has the Gibbs distribution on $\Qm^{\ell,f}$ for every $f$ and is invariant under rerooting.
    \item $(\P^{\ell,b})_{(\ell,b)\in \N^*\times \Bc}$ has the Markov property and for each $\ell$ and $b$ the measure $\P^{\ell,b}$ is invariant under rerooting.
\end{enumerate}
\end{thm}

Before starting the proof, assume that $(\P^{\ell,b})_{(\ell,b)\in \N^*\times \Bc}$ is like in the theorem and satisfies the Markov property. Take $\q$ with half-perimeter $\ell$ and $\sigma:V(\q)\mapsto \supp \mu$ with boundary condition $b$. Note that the law of $\sigma$ on the event that  $\q=\BM$ is absolutely continuous with respect to $\mu^{\sharp V_i(q)}$ and the Radon-Nikodym derivative is given by
\begin{align} \label{e.def_p_q_s}
    p^{\ell,b}(\q,\sigma) := \frac{d \P^{\ell,b}(\BM=\q, \Dec|_{V_i(\q)} \in d\sigma) }{d \mu^{\sharp V_i(\q)}(d\sigma)} (\sigma) = \P^{\ell,b}(\BM=\q)\frac{ \exp(-\Ham^b( \q,\sigma)) }{Z^b(\q)}.
\end{align}
Here the last equality follows from \eqref{def_cond}. Note furthermore that $p^{\ell,b}(\q,\cdot)$ is continuous.
Now, we also need to show that on the event that the first peeling step is of type $1$, the value of the decoration $\Dec$ in the peeled faced is absolutely continuous with respect to $\mu$.

\begin{lemma}\label{l.continuity_decorated}
Let us work in the context of Theorem \ref{t.uniqueness_decorated}, assuming also that $(\P^{\ell,b})$ has the Markov property. On the event where the first peeling step is of type $1$, call $\vertex$ the new vertex discovered. Then, the law of $\Dec_\vertex$ is absolutely continuous with respect to $\mu$. Furthermore,
\begin{align*}
    p^{\ell,b}_1(x):= \frac{d\P^{\ell,b}(\text{Type 1}, \Dec_\vertex\in dx)  }{d\mu(dx) }(x)
\end{align*}
is continuous\footnote{To be more precise, it has a version (i.e. a representative) in $L^1(\mu)$ which is continuous. We always use this version as the equalities we obtain are all almost everywhere equalities on $\mu$ and Lemma \ref{l.continuity_decorated} implies that we can extend these equalities to all points} on the support of $\mu$.
\end{lemma}
\begin{proof}
To show absolute continuity with respect to $\mu$, we do as follows. Take $A$ such that $\mu(A)=0$. Then by the definition of an SDM for all deterministic $\q$ such that the peeling of the root vertex gives us a peeling of type $1$ we have that $\P^{\ell,b}(\Dec_\vertex\in A,  \BM=\q)=0$. We conclude by summing over all such $\q$.

Now, to show continuity we proceed as follows. Assume that the probability of having a peeling step of type 1 is positive. Take $x\in \supp \mu$ and take $\q$ such that $\P^{\ell+1,\hat b_x}(\q)>0$, where $\hat b_x$ is the function that takes values $x$ in $0$, $1$ and $2$ and values $b(j-2)$ for all $j\geq 3$. Define $\r$ as the gluing between the result of a peeling step of type $1$ of a map with boundary size $\ell$ with $\q$. Take $\sigma$ a decoration of $\r$ and  $\sigma^{(x)}$ taking value $x$ in the vertex adjacent to the root edge and $\sigma$ elsewhere, we have that
\begin{align}
p^{\ell,b}(\r,\sigma^{(x)})= p_1^{\ell,b}(x) \P^{\ell+1,\hat b_x}(\BM=\q) \frac{e^{-\Ham^b( \r, \sigma)} }{Z^b(\r) }.
\end{align}
As $\P^{\ell+1,\hat b_x}(\BM=\q)$ and $p^{\ell,b}(\r,\sigma^{(x)})$ are continuous and non-zero in $x$, we see that $p_1^{\ell,b}(x)$ is a.e. continuous on $x$.
\end{proof}

We can now start the proof of the main theorem of this section.

\begin{proof}[Proof of Theorem \ref{t.uniqueness_decorated}]\label{proof.sec4_charac}
This proof follows closely that of Theorem \ref{t.characterization_metric}. In fact, most of the proof is analogous to that of Theorem \ref{t.characterization_metric} and is even simpler in this case. The implication $(2) \Rightarrow (1)$ mimics the same argument as before. For $\ell\in\Z_+$, $b\in \Bc$, $\q\in\Qm_H$ and $\sigma$ a decoration over $V(\q)$, we define $q(\ell,b,\q,\sigma)$ as follows
\begin{align*}
        q(\ell,b,\q,\sigma) := \frac{p^{\ell,b}(\overline{\q},\overline{\sigma})}{p^{\ell,b}( \q,\sigma)},
\end{align*}
where $\overline{\q}$ is the skeleton of the metric map $\overline{q}(w_1,w_2)$ defined in the proof of Theorem \ref{t.characterization_metric}. Following the same procedure we can prove that $q(\cdot,\cdot,\cdot,\cdot)$ is constant, i.e. does not depend on any of the parameters; we leave it to the reader to complete this proof.

For the implication $(3) \Rightarrow (2)$ a new difficulty appears. First, we have to show that for any $\q_1,\q_2 \in \Qm^{\ell,v}$
    \begin{align}\label{e.Gibbs}
    \frac{\P^{\ell,b}(\BM = \q_1) }{\P^{\ell,b}(\BM = \q_2) }=\frac{ Z^b(\q_1)}{Z^b(\q_2) }.
    \end{align}
    We prove this by induction on the number of interior vertices (both maps have to have the same number). However, this proof is trickier than the metric case\footnote{Here is where we use that for any $\ell$ and $b$ we have that $\P^{\ell,b}(\BM \text{ has at least one internal vertex})>0$ (hypothesis of Theorem \ref{t.uniqueness_decorated}), the reason being that if there is probability 0 to peel a new vertex, then the measure will be restricted to trees.}. To initiate the induction, we need the following claim. 
    \begin{claim}\label{c.start_induction_trees}
 Take two trees $\mathbf{t}_1$ and $\mathbf{t}_2$ with $\ell$ edges, we have that
    \begin{align}\label{e.ratio_trees}
        \frac{\P^{\ell,b}(\BM = \mathbf{t}_1)}{\P^{\ell,b}(\BM = \mathbf{t}_2)} = \frac{Z^b(\mathbf{t}_1)}{Z^b(\mathbf{t}_2)}.
    \end{align}    \end{claim}
    As this result is technical, we defer its proof to the end and for now we assume it.
    
   Take $n\in \N$ and assume that \eqref{e.Gibbs} holds for maps with no more than $n$ interior vertices. Choose $\q_1$ with $n+1$ vertices. At least one internal vertex has to share an edge with the boundary. By rerooting invariance, we assume that the root edge is adjacent to an internal vertex we call $\vertex$. Let $\q_2$ be another map with $n+1$ vertices and suppose that that the root edge is adjacent to an internal vertex, that we call $\vertex$ too (we will deal with the other case later). Additionally, define $\q_i^{*}$ as the result of applying a peeling step in the root of $\q_i$ and for $x\in \supp(\mu)$
   \begin{align*}
   p^{\ell,b}(\q_i,x) &:= \frac{d \P(\BM=\q_i, \Dec(\vertex)\in dx) }{\mu(dx) }(x)
   \end{align*}
  is the Radon-Nikodym derivative of the decoration at the vertex $\vertex$ with respect to the measure $\mu$ on the event that the map is exactly $\q_i$.
   
  To show that \eqref{e.Gibbs} is true for $\q_1$ and $\q_2$, we compute   
   \begin{align}\label{eq:quot:RN}
   p^{\ell,b}(\q_i, b(0))&= \P^{\ell,b}(\BM=\q_i)\frac{1 }{Z^b(\q_i) }\int \exp(-\Ham^{b^{*}}(\q_i^{*}, \sigma )) \mu^{n-1}(d\sigma)= \P^{\ell,b}(\BM=\q_i)\frac{Z^{b^*}(\q_i^{*}) }{ Z^{b} (\q_i)},
   \end{align}
   where $b^{*}$ is the boundary condition\footnote{To be more precise, it is exactly equal to $b$ everywhere, except that now $0$ has become three edges with value equal to $b(0)$} that appears when doing a peeling step on the vertex $\vertex$ and discovering that it has value $b(0)$. 
   
  Then, we can use the Markov property to peel the vertex $\vertex$ and see that, thanks to \cref{eq:quot:RN}, $\P^{\ell,b}(\BM=\q_1)/Z^b(\q_1)$ is equal to
   \begin{align*}
    \frac{ p^{\ell,b}(\q_1,  b(0))}{Z^{b^{*}}(\q_1^{*}) } = \frac{p^{\ell,b}_1(b(0))\P^{\ell+1, b^{*}}(\BM=\q_1^{*}) }{Z^{b^*}(\q_1^*) }
    = \frac{p^{\ell,b}_1(b(0))\P^{\ell+1, b^*}(\BM=\q_2^*) }{Z^{b^*}(\q_2^*) }   = \frac{\P^{\ell,b}(\BM=\q_2) }{Z^b(\q_2)}.
   \end{align*}
   Note that in the second equality we used the induction hypothesis as $\q_1^*$ and $\q_2^*$ have $n$ vertices.
   
   We have shown that as long as $\q_1$ and $\q_2$ both have interior vertices that are adjacent to the root edge they satisfy \eqref{e.Gibbs}. To extend this for general $\q_1$ and $\q_2$, we argue as follows. First, we assume by rerooting invariance that $\q_1$ has a vertex adjacent to the root edge.  Then, we take $\mathfrak e$ an edge of the outer boundary of $\q_2$ that is adjacent to an internal vertex (this is possible as $\q_2$ has at least one internal vertex). Now, we construct $\q_3$ a map that has a vertex $\vertex$ such that its boundary $\partial \vertex$ contains both the root edge and $\mathfrak e$. See \cref{fig:q3} for a formal explanation.
   \begin{figure}[h!]
       \centering
       \includegraphics[scale=0.5]{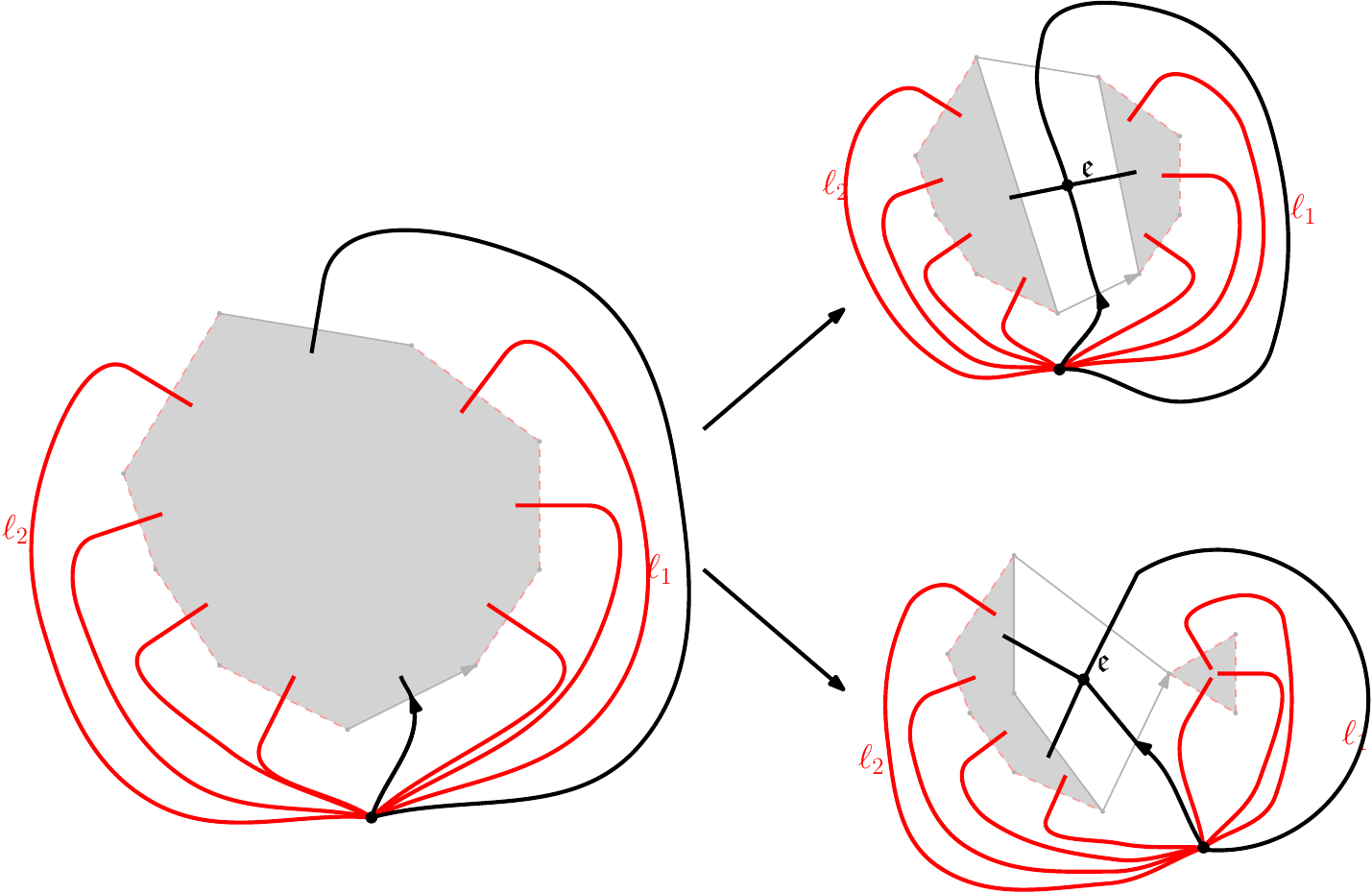}
       \caption{An illustration of the transformation according to the parity of $\ell_1$. Starting from the tip of the root edge and following its orientation along the boundary, we count the edges until reaching $\mathfrak{e}$, obtaining $\ell_1$. If $\ell_1$ is odd, we perform the transformation shown in the upper-right panel, whereas if $\ell_1$ is even, we use the transformation shown in the lower-right panel. The grey regions represent hole vertices with even boundary length, which can be filled with 4-regular maps with a boundary so that the resulting maps have the same number of interior vertices.}
       \label{fig:q3}
   \end{figure}
   Thus,
   \begin{align*}
   \frac{\P^{\ell,b}(\q_1) }{\P^{\ell,b}(\q_2) }= \frac{\P^{\ell,b}(\q_1) }{\P^{\ell,b}(\q_3) }\frac{\P^{\ell,b}(\q_3) }{\P^{\ell,b}(\q_2) }= \frac{Z^b(\q_1) }{ Z^b(\q_3)}\frac{Z^b(\q_3) }{ Z^b(\q_2)}=\frac{Z^b(\q_1) }{ Z^b(\q_2)}.
   \end{align*}

   To finish, we just need to show the claim.
   
   \begin{proof}[Proof of Claim \ref{c.start_induction_trees}] We proceed by induction on $\ell$, this is obviously true for $\ell=1$, however our base case needs to be $\ell=2$. Surprisingly, the induction step is more straightforward than the base case, so we begin with it.
   
   \underline{Induction step:} Assume that \eqref{e.ratio_trees} is true for all trees of size $\ell$ with $\ell\geq 2$. Take $\mathbf t_1$ and $\mathbf t_2$ two trees of size $\ell+1$ that ``share a leaf'', meaning that there is a $k\in \{0,..,2\ell+1\}$ such that the boundary vertex labeled $k$ is a leaf for both trees, in particular this means that the semi-edges $k$ and $k+1$ are identified in both trees. By rerooting invariance, we assume that both of them associate the root edge to next edge counter-clockwise. Then for any boundary condition $b$
   \begin{align*}
       \P^{\ell+1,b}(\BM={\mathbf t}_1) = \P^{\ell+1,b}(\text{type 2},0,\ell) \P^{\ell,b^{*}}(\BM={\mathbf t}_1^*) =\P^{\ell+1,b}(\BM={\mathbf t}_2) \frac{Z^{b^*}( {\mathbf t}_1^*) }{ Z^{b^*} ({\mathbf t}_2^*) },
   \end{align*}
   where $\tr_1^*$ and $b^*$ are the results of applying one peeling step on the root to $\tr_1$ and $b$ respectively. We conclude this case by noting that 
   \begin{align*}
       Z^{b^*} ({\mathbf t}_i^*)=  Z^b ( {\mathbf t}_i)\exp\left(\frac{\beta}{ 2} (b(0)-b(1))^2\right).
   \end{align*}
  
   For the general case, we construct an intermediate tree of size $\ell+1$ having two special leaves, one matching a leaf of $\mathbf{t}_1$ and another matching a leaf of $\mathbf{t}_2$. Such a tree always exists if $\ell+1\geq 3$; this is left to the reader since the argument is similar to that of \Cref{fig:q3}.

   \underline{Base case:} Now, we verify the claim for $\ell=2$ and an arbitrary boundary condition $b\in \Bc$. In this case, we numerate the edges of the boundary from $0$ to $3$. Define $\mathbf t_1$ as the tree where $0$ and 1, and 2 and 3 are glued together and $\mathbf t_2$ as that where $0$ and $3$, and $1$ and $2$ are glued together. Note that $\mathbf t_1$ and $\mathbf t_2$ are the only two trees that can be obtained in this setup. See \Cref{fig:2trees} for an explanation on the resulting $\mathbf{t}_1$ and $\mathbf{t}_2$.
   \begin{figure}[h!]
       \centering
       \includegraphics[scale = 0.65]{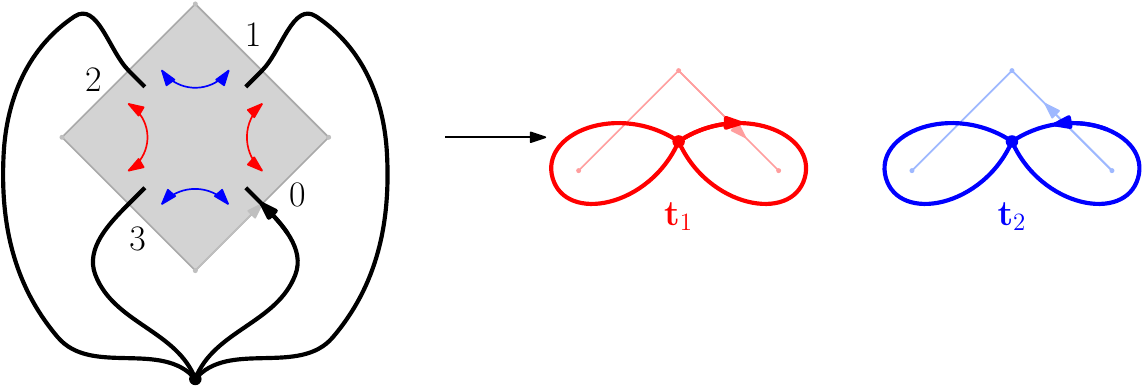}
       \caption{Two possible identifications of the four half-edges incident to the hole vertex, producing the maps $\mathbf{t}_1$ and $\mathbf{t}_2$. The identifications corresponding to $\mathbf{t}_1$ are shown in red, while those corresponding to $\mathbf{t}_2$ are shown in blue. The resulting maps are displayed on the right.}
       \label{fig:2trees}
   \end{figure}
   
  Now, we need to define $\tr_i^\boxplus$. Let $\boxplus$ be the 4-regular map with hole vertices presented in \Cref{fig:t1square} and let $\tr_i^\boxplus= \boxplus \glue \tr_i$ and $\sigma^\boxplus$ its corresponding decoration where each vertex that is adjacent to the root vertex receives spin equal to the spin in the phantom vertex. 
   
    \begin{figure}[h!]
       \centering
       \includegraphics[scale = 0.75]{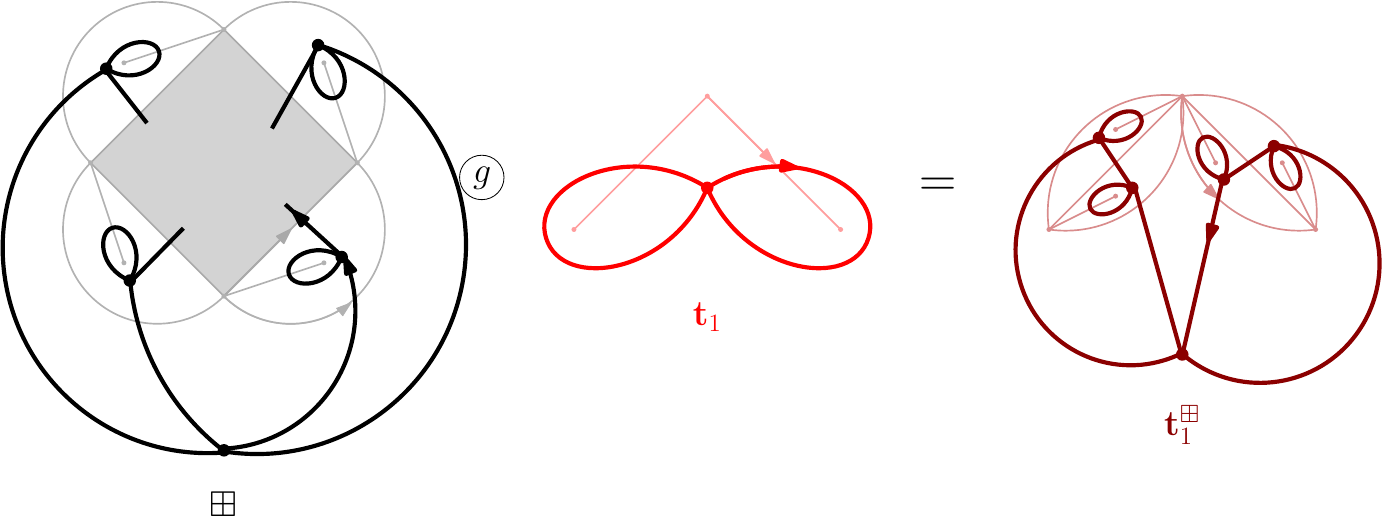}
       \caption{An example of the gluing construction. The left panel shows a map with a hole vertex, represented by the grey region. The middle panel shows the auxiliary map $\tr_1$, which is glued into this hole to obtain the map displayed on the right. In the decorated setting, each ordinary interior vertex of the resulting map carries an associated spin.}
       \label{fig:t1square}
   \end{figure}

   Define $s$ the counterclockwise shift $s(i)=i+1\mod 4$ and note $q_4=q(2,b)q(2,b\circ s)q(2,b\circ s^2) q(2,b\circ s^3)$. Thanks to the Markov property and by iterating the construction use to construct $\overline{\q}$ we have that
   \begin{align}
       p^{2,b}(\tr_i^{\boxplus},\sigma^\boxplus)= q_4\P^{2,b}(\BM=\tr_i).
   \end{align}
   
   Now, define $b_s=b\circ s$ and $\sigma^{\boxplus_s}$ as the counterclockwise shift of $\sigma^{\boxplus}$ in the inner vertices\footnote{Note that $\sigma^{\boxplus_s}$  has the same boundary conditions as $\sigma^{\boxplus}$}. We have that
   \begin{align*}
       p^{2,b}(\tr_1^{\boxplus},\sigma^{\boxplus_s})&= p^{2, b}(\tr_1^{\boxplus},\sigma^{\boxplus}) \exp\left(
       {\Ham^b(\tr_1) - \Ham^b(\tr_2)-\frac{\beta }{2 }\sum_{i=0}^3(b(i)-(b_s(i))^2}\right)\\
       &= p^{2, b}(\tr_1^{\boxplus},\sigma^{\boxplus}) \exp\left({-\frac{\beta }{2 }\sum_{i=0}^3(b(i)-b_s(i))^2}\right) \frac{Z^b(\tr_2) }{Z^b(\tr_1) }.
   \end{align*}
   and furthermore, by doing ``properly chosen peeling steps'', meaning that we pay the price to transform $t_1$ into $t_2$ with the corresponding decoration, we also obtain that
   \begin{align*}
       p^{2,b}(\tr^{\boxplus}_1, \sigma^{\boxplus_s}) = r(b_s) \P^{2,b}(\BM=\tr_2),
   \end{align*}
   where $r(\cdot)$ is defined as follows. For a function $c:\cro{0,3} \mapsto \supp(\mu)$, we define $\sigma_c: V_i(\boxplus)\mapsto \supp(\mu)$ such that $\left(\sigma_c\right)_{\vertex}= c(i)$, where $i$ is the edge of $\vertex$ that belongs to the boundary of $\boxplus$. Then,
   \begin{align}
       r(c) = \frac{d \P^{2,b}(\boxplus \Lalg \BM, \Dec\mid_{V_i(\boxplus)} \in dx ) }{d\mu^4(dx)}(\sigma_c).
   \end{align}
  In particular $r(b)=q_4$.
   
 Note that now, to conclude, we need to verify that
   \begin{align}\label{e.needed_for_r}
       r(b_s)= q_4\exp\left({-\frac{\beta }{2 }\sum_{i=0}^3(b(i)-b_s(i))^2}\right),
   \end{align}
 and that $q_4$ is not $0$.
   
   To prove \eqref{e.needed_for_r}, we compute $r(b_s)$, see \Cref{fig:rbs}.
   \begin{align*}
       r(b_s)&= \sum_{\q} \int p^{2,b}(\boxplus \glue \q,\sigma^{\boxplus_s}+\sigma) \mu^{\sharp V_i(\q)}(d\sigma)\\
       &= \sum_{\q} \int p^{2,b}(\boxplus \glue \q ,\sigma^{\boxplus}+\sigma)  \exp\left({-\frac{\beta }{2 }\sum_{i=0}^3(b(i)-b_s(i))^2 - (b(i) - \sigma_{\vertex_i} )^2  +(b_s(i)-\sigma_{\vertex_i})^2 }\right) \mu^{\sharp V_i(\q)} (d\sigma),
   \end{align*}
  here $\vertex_i$ is the vertex of $\mathfrak r \glue \q$ that appears when doing a peeling step on the edge $i$ of $\q$.
  
  \begin{figure}[h!]
       \centering
       \includegraphics[scale=0.8]{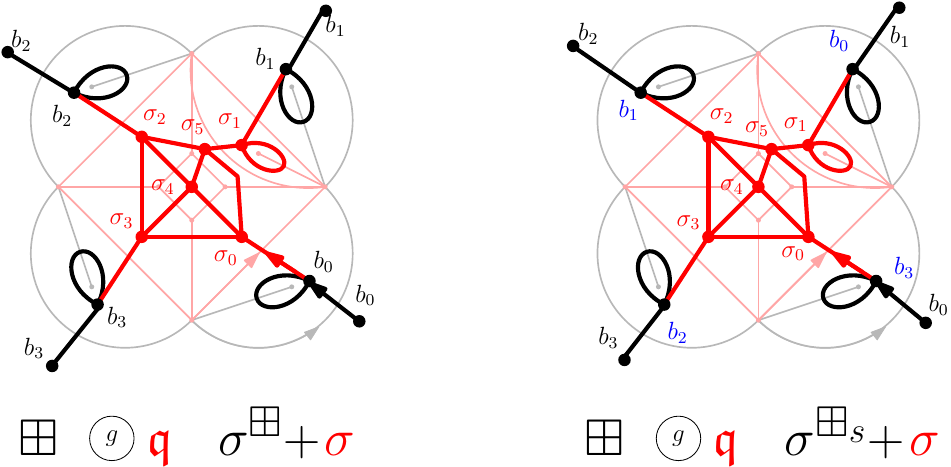}
      \caption{Examples of the configurations appearing in the computation of $r(b_s)$. The left panel shows a decorated gluing configuration, while the right panel shows the corresponding configuration after modifying the boundary-spin assignment. In both panels, the elements glued into the original map and their associated spins are shown in red. The boundary spins affected by the transformation are colored in blue in the right panel, while the remaining boundary spins are shown in black. The auxiliary spin configuration is determined by the boundary data $b$.}
       \label{fig:rbs}
   \end{figure}
   
   Now, denote $s(\q)$ as the map $\q$ but with the root edge moved one edge counter-clockwise and note that this is a bijection on the space of 4-regular maps with half-perimeter 2, and note that for a decoration $\sigma$ of $\q$ we can associate a decoration $s(\sigma)$ of $\sigma(\q)$ by properly renaming the coordinates of $\sigma$. A direct computation implies that
   \begin{align*}
       p^{2,b}(\mathfrak r \glue s(\q), \sigma^{\boxplus} + s(\sigma)) = p^{2,b}(\mathfrak r \glue \q, \sigma^{\boxplus}+\sigma) \exp\left (\frac{\beta }{2 }\sum_{i=0}^3 (b_s(i)-\sigma_{\vertex_i})^2 - (b(i)-\sigma_{\vertex_i})^2 \right ).
   \end{align*}
   This implies that $r(b_s)$ is equal to
   \begin{align*}
       &\sum_{\q} \int p^{2,b} (\boxplus \glue s(\q), \sigma^{\boxplus}+s(\sigma)) \exp\left({-\frac{\beta }{2 }\sum_{i=0}^3(b(i)-b_s(i))^2 }\right) \mu^{\sharp V_i(\q)}(d\sigma)\\
       = & \ r(b)\exp\left({-\frac{\beta }{2 }\sum_{i=0}^3(b(i)-b_s(i))^2 }\right) \\
       = & \ q_4 \exp\left({-\frac{\beta }{2 }\sum_{i=0}^3(b(i)-b_s(i))^2 }\right).
   \end{align*}

Now we only need to show that $q_4>0$. This follows from the fact that 
\begin{align*}
\P^{2,b}(\BM \text{ has at least one interior vertex}) >0.
\end{align*} This shows that for any boundary condition at least 2 edges have positive probability of having a vertex. Thanks to the fact that $\BM$ is an SDM and continuity of the measure on $\ell$, it is possible to see that if $q_4$ is $0$ there has to be an $i\in \cro{1,3}$ and at least two values $b(0)$ and $b(i)$ such that if the phantom face $0$ takes values $b(0)$ and $i$ takes values $b(i)$ then, under $\P^{\ell,b}$, $\BM$ always identifies the semi-edge $0$ with the semi-edge $i$. However, this is a contradiction as follows: assume for example that $i=1$, then take $b(0)=b(2)=b_0$ and $b(1)=b(3)=b_1$. Thanks to the shift invariance $\BM$ would a.s. have no vertices under $\P^{\ell,b}$. 
\end{proof}
   
\end{proof}

\begin{rem}
Let us point out that the assumption that $\BM$ has at least one internal vertex with positive probability is key for Theorem~\ref{t.uniqueness_decorated}. Indeed, the implication $(3)\Rightarrow(2)$ fails if a.s. $\BM$ has no interior vertices. In that case, $\BM$ is always a tree and the condition that $\P^{\ell,b}$ is a spin-decorated map is superfluous, since trees have no interior vertices on which to place the decoration. Consequently, the law $\P^{\ell,b}$ has too many degrees of freedom, and there is no mechanism forcing the factor $e^{-\Ham(\tr,\sigma)}$ to appear.

This issue arises in the proof when establishing the base case of the induction for trees. In the metric decorated setting, this is not necessary because we can stop in the middle of an edge, which allows us to change the boundary condition. When we do not have lengths, one needs to introduce a map with interior vertices in order to compare the laws of trees and recover the Hamiltonian weight.
\end{rem}

\bibliographystyle{alpha} 
\bibliography{biblio_2}
\appendix

\end{document}